\documentclass[10pt]{article}

\usepackage{amsfonts}
\usepackage{amsmath}
\usepackage{amsthm}
\usepackage{graphicx}
\usepackage{setspace}
\usepackage{makeidx}
\usepackage{setspace}
\usepackage[footnotesize]{caption2}

\makeindex

\newtheorem{thm}{Theorem}[section]
\theoremstyle{definition}

\theoremstyle{definition}

\theoremstyle{remark}

\theoremstyle{plain}
\newtheorem{lem}[thm]{Lemma}
\newtheorem*{aux}{Auxiliary lemma}
\newtheorem*{main}{Main theorem}

\title{Uniform asymptotics for the full moment conjecture of the Riemann zeta function}
\author{Ghaith A. Hiary and Michael O. Rubinstein
\footnote{Both authors are supported by the
National Science Foundation under awards DMS-0757627 (FRG grant) and
DMS-0635607. In addition, the second author is supported by an NSERC Discovery Grant.}}
\date{}

\begin{document}
\maketitle

\begin{abstract}
Conrey, Farmer, Keating, Rubinstein, and Snaith, recently conjectured formulas
for the full asymptotics of the moments of $L$-functions. In the case of the
Riemann zeta function, their conjecture states that the $2k$-th absolute moment
of zeta on the critical line is asymptotically given by a certain $2k$-fold
residue integral. This residue integral can be expressed as a polynomial of
degree $k^2$, whose coefficients are given in exact form by
elaborate and complicated formulas.

In this article, uniform asymptotics for roughly the first $k$ coefficients
of the moment polynomial are derived. Numerical data to support our
asymptotic formula are presented. An application to bounding the maximal size of the zeta
function is considered.
\end{abstract}

\tableofcontents

\section{Introduction}

The absolute moments of the Riemann zeta function  on the critical line are a
natural statistical quantity to study in connection with value distribution
questions. For example, they can be used to understand the maximal size of
the zeta function. These moments are also connected to the remainder term in
the general divisor problem~\cite{T}.

Hardy and Littlewood proved a leading-term asymptotic for the second
moment on the critical line \cite{HL}. A few years later, in 1926,
Ingham gave the full asymptotic expansion \cite{I}. In the same
article, Ingham gave a leading term asymptotic for the fourth moment. The full
asymptotic expansion for the fourth moment was obtained by Heath-Brown in 1979
\cite{HB}. In comparison, the higher moments seemed far more difficult and
mysterious. Keating and Snaith, in a breakthrough, conjectured the leading-term
asymptotic \cite{KS}.

Recently, however, based on number-theoretic considerations, Conrey, Farmer,
Keating, Rubinstein, and Snaith, conjectured~\cite{CFKRS1}~\cite{CFKRS2} the
following full asymptotic expansion for the $2k$-th absolute moment of the Riemann
zeta function $\zeta(s)$ on the critical line:
\begin{eqnarray}\label{eq:Pkx}
    \frac{1}{T}\int_0^T |\zeta(1/2+it)|^{2k}\,dt \sim 
    \frac{1}{T} \int_0^T P_k\left(\log\frac{t}{2\pi}\right)\,dt\,,
    \qquad \textrm{as $\,\,\,T\to \infty$}\,,
\end{eqnarray}
where $P_k(x)$ is a polynomial of degree $k^2$:
\begin{eqnarray}
    P_k(x)=:c_0(k) x^{k^2}+c_1(k) x^{k^2-1}+\cdots+c_{k^2}(k) \,,
\end{eqnarray}
given \textit{implicitly} by the $2k$-fold
residue
\begin{eqnarray}\label{eq:cfkrs00}
    P_k(x)=\frac{(-1)^k}{k!^2}\frac{1}{(2\pi i)^{2k}}\oint \cdots \oint 
    &&\frac{G(z_1,\ldots,z_{2k})\Delta^2(z_1,\ldots,z_{2k})}{\prod_{i=1}^{2k}z_i^{2k}} \\
    &&\times e^{\frac{x}{2}\sum_{i=1}^k z_i-z_{k+i}}\,dz_1\ldots dz_{2k}\,,\nonumber
\end{eqnarray}
where the path of integration is around small circles enclosing $z_i=0$, and 
\begin{equation}
    \Delta(z_1,\ldots,z_{2k}) :=\prod_{1\le i< j \le 2k} (z_j-z_i)
\end{equation}
is the Vandermonde determinant, and
\begin{eqnarray}
    \label{eq:G}
    G(z_1,\ldots,z_{2k}):=A(z_1,\ldots,z_{2k})\prod_{i,j=1}^k \zeta(1+z_i-z_{k+j})\,,
\end{eqnarray}
is a product of zetas and the ``arithmetic factor'' (Euler product)
\begin{eqnarray}
    \label{eq:A_def}
    &&A(z_1,\ldots,z_{2k})  \nonumber \\
    &&:= \prod_p \prod_{i,j=1}^k \left(1-p^{-1-z_i+z_{k+j}}\right)\,
    \int_0^1 \prod_{j=1}^k \left(1-\frac{e^{2\pi i \theta}}{p^{\frac{1}{2}+z_j}} \right)^{-1}
    \left(1-\frac{e^{-2\pi i\theta}}{p^{\frac{1}{2}-z_{k+j}}}\right)^{-1} d\theta \nonumber \\
    && \\
    &&= \prod_p \sum_{j=1}^k \prod_{i\ne j} \frac{\displaystyle \prod_{m=1}^k (1-p^{-1+z_{i+k}-z_m})}{1-p^{z_{i+k}-z_{j+k}}}\,.
\end{eqnarray}
As pointed out by \cite{CFKRS1},  the rhs of (\ref{eq:cfkrs00}) has an almost
identical form to an exact expression for the moment polynomial of random
unitary matrices, the difference being that $G(z_1,\ldots,z_{2k})$ is replaced
by the function $\prod_{i,j=1}^k (1-e^{z_{j+k}-z_i})^{-1}$ in the unitary case,
so there is no arithmetic factor.

The CFKRS conjecture (\ref{eq:cfkrs00}) agrees with the theorems of Hardy and
Littlewood, Ingham, and Heath-Brown, for $k=1$ and $k=2$. It has been supported
numerically; see~\cite{CFKRS1}, \cite{CFKRS2} \cite{HO} \cite{RY}. The
conjecture provides a method for computing the lower order coefficients of the
moment polynomial $P_k(x)$. It gives, in particular, a stronger asymptotic than
that of Keating and Snaith who, by carrying out an analogous computation for
random unitary matrices, first predicted the leading coefficient (see~\cite{KS}):
\begin{equation}
    c_0(k) = \frac{a_k g_k}{k^2!}\,,
\end{equation}
where
\begin{equation}
    \label{eq:ak def}
    a_k :=
    \prod_p (1-1/p)^{k^2} F(k,k;1;1/p)\,,
\end{equation}
and
\begin{equation}
    g_k := k^2!\prod_{j=0}^{k-1}\frac{j!}{(j+k)!}\,.
\end{equation}

\subsection{Results}

Our main theorem develops a uniform asymptotic for $c_r(k)$ in the region $0\le
r\le k^{\beta}$, for any fixed $\beta<1$. We expect the asymptotics can be
corrected so as to remain valid well beyond the first $k$ coefficients (i.e.
for $\beta \ge 1$), and that the methods in our paper, which are of
combinatorial nature, will be helpful in deriving uniform asymptotics for the
moments of other $L$-functions.

To state our main theorem, let us first define
\begin{equation}
    B_k:= \sum_p \frac{k\log p}{p-1} - \frac{F(k+1,k+1; 2; 1/p)}{F(k,k;1;1/p)}\,\frac{\log p}{p}\,,
\end{equation}
where $F(a,b;c;t)$ is the Gauss hypergeometric function
\begin{equation}\label{eq:2f1}
    F(a,b;c;t) := \frac{\Gamma(c)}{\Gamma(a)\Gamma(b)}
    \sum_{n=0}^{\infty} \frac{\Gamma(a+n)\Gamma(b+n)}{\Gamma(c+n)}\,\frac{t^n}{n!}\,.
\end{equation}
In the notation of \cite{CFKRS2}, $B_k$ is the same as
 $B_k(1;)$, which is given in Eqs.~(2.24) and (2.43) there.
The factor $B_k$ is arithmetic in nature. It is the coefficient of the
linear term in the following Taylor expansion of the arithmetic factor:
\begin{eqnarray}
    \log A(z_1,\ldots,z_{2k})=\log a_k+B_k\sum_{i=1}^k z_i-z_{k+i}+\cdots\,,
\end{eqnarray}
where it is known (see 2.7 of \cite{CFKRS1}) that
\begin{equation}
    a_k = A_k(0,\ldots,0)\,.
\end{equation}
Theorem~\ref{thm:lincoeff} will later furnish the following asymptotic for $B_k$:
\begin{equation} \label{eq:cfkrsb}
    B_k\sim 2\,k\,\log k\,,\qquad \textrm{as $\,\,\,k\to \infty$} \,.
\end{equation} 

\vspace{1mm}

\begin{main}
\label{thm:main}
Fix $\beta < 1$, let $0\le r\le k^{\beta}$, and let
\begin{equation}
    \label{eq:bound tau}
    \tau_k:=2B_k+2\gamma k\,,
\end{equation}
where $\gamma = 0.5772\ldots$ is the Euler constant. Notice by
(\ref{eq:cfkrsb}) we have
\begin{equation}
    \label{eq:tau_k bound}
    \tau_k\sim 4\,k\,\log k\,,
    \qquad \textrm{as $\,\,\,k\to \infty$} \,.
\end{equation}
Then as $k\to \infty$, and uniformly in $0\le r\le k^{\beta}$,
\begin{eqnarray}
    \label{eq:main1}
    c_r(k)&=& \tau_k^r \,\binom{k^2}{r}\, \frac{a_k g_k}{k^2!}
    \left[1+O\left(\frac{r^2}{k^2}\right)\right]\\
    \,\nonumber\\
    \label{eq:main2}
    &=& \tau_k^r \,\binom{k^2}{r}\, c_0(k)\,
    \left[1+O\left(k^{2(\beta-1)}\right)\right]\,.
\end{eqnarray}
Alternatively, 
\begin{eqnarray}
    c_r(k)= \frac{\tau_k^r k^{2r}}{r!}\,\, c_0(k)\, \left[1+O\left(k^{2(\beta-1)} \right)\right]\,,
\end{eqnarray}
as $k\to \infty$. Asymptotic constants depend only on $\beta$.
\end{main}

\vspace{2mm}

\noindent
Remarks: 1) The asymptotic formulas~\eqref{eq:main1} and~\eqref{eq:main2} of our theorem are actually
equalities for $r=0$, and $r=1$. The $r=0$ case is trivial, and the $r=1$ case
follows from either (2.71) of~\cite{CFKRS2} or~\eqref{eq:main identity} below. 
 2) For comparison, the corresponding asymptotic in the unitary case, provided in \cite{HR}, is:
\begin{equation}
    \label{eq:unitary_asymp}
    \tilde{c}_r(k) = k^r\binom{k^2}{r} \tilde{c}_0(k)\left[1 + O\left(\frac{r^2}{k^2}\right)\right]\,,
\end{equation}
where $\tilde{c}_r(k)$ is the coefficient of $x^{k^2-r}$ in the $2k$-th moment polynomial
 of random unitary matrices.

Although the CFKRS conjecture seems hopelessly difficult to prove, the precise
nature of the asymptotic formula allows one to gain insight into the behavior
of the zeta function. For example, by deriving an asymptotic for $c_r(k)$ that
is applicable as $r$ and $k$ both tend to infinity, one can understand the true
size of $\zeta(1/2+it)$. The results we present here are a step in this
direction.

One difficulty in extracting uniform asymptotics for the coefficients of
$P_k(x)$ from a residue like (\ref{eq:cfkrs00}) is that the coefficients are
given only implicitly. By comparison, both the coefficients and the roots of
the moment polynomials for random unitary matrices, which correspond to the
zeta-function moment polynomials according to the random matrix philosophy, are
known explicitly, via random matrix theory calculations. In fact, the proof of
Theorem~1 of \cite{HR}, which provides complete uniform asymptotics for the
coefficients in the unitary case, makes essential use of the information about
the roots via a saddle-point technique. In the case of the zeta function,
however, we do not have `simple' closed form expressions for the moment
polynomials.

We remark that if one directly applies the methods of this paper to the residue
expression for unitary moment polynomials, given in \cite{CFKRS1} Eq. (1.5.9),
then one encounters similar difficulties as in the zeta function (e.g. a
similar difficulty in deriving asymptotics beyond the first $k$ coefficients).
The main added simplicity in the unitary case is that it does not involve an
arithmetic factor.

Before delving into the careful details of the next sections, let us describe
the basic idea of the proof. To this end, define
\begin{eqnarray}
    \label{eq:R def}
    R(z_1,\ldots,z_{2k}):=G(z_1,\ldots,z_{2k}) \prod_{i,j=1}^k (z_i-z_{k+j})\,.
\end{eqnarray}
where, recall, $G(z_1,\ldots,z_{2k}) =A(z_1,\ldots,z_{2k})\prod_{i,j=1}^k
\zeta(1+z_i-z_{k+j})$. The extra product on the rhs in (\ref{eq:R def}) is
introduced in order to cancel the poles in the product of zetas in the
definition of $G(z_1,\ldots,z_{2k})$. This renders the function
$R(z_1,\ldots,z_{2k})$ analytic and non-zero in a neighborhood of the origin,
where it is equal to $a_k$. Therefore, we may write
\begin{eqnarray}
    \label{eq:Pk def G}
    P_k(x)= \frac{(-1)^k}{k!^2}\frac{1}{(2\pi i)^{2k}}\oint \cdots \oint
    &&\frac{\Delta^2(z_1,\ldots,z_{2k})\,e^{\frac{x}{2}\sum_{i=1}^k z_i-z_{k+i}}}
    {\prod_{i,j=1}^k (z_i-z_{k+j}) \prod_{i=1}^{2k}z_i^{2k}}  \nonumber \\
    && \\
    && \times \,e^{\log R(z_1,\ldots,z_{2k})} \,dz_1\ldots dz_{2k}\,,\nonumber
\end{eqnarray}
and consider the Taylor expansion of $\log R(z_1,\ldots,z_{2k})$:
\begin{eqnarray}
    \label{eq:logR expand}
    \log R(z_1,\ldots,z_{2k})=\log a_k + \frac{\tau_k}{2}\sum_{i=1}^k z_i-z_{k+i}+\cdots\,,
\end{eqnarray}
where, recall, $\tau_k =2B_k+2\gamma k \sim 4k\log k$, as $k\to \infty$.
Also, dropping the factor $\exp(\log R(z_1,\ldots,z_{2k}))$, define
\begin{eqnarray}
    \label{eq:Lk def}
    p_k(x,0):= \frac{(-1)^k}{k!^2}\frac{1}{(2\pi i)^{2k}}\oint \cdots \oint
    \frac{\Delta^2(z_1,\ldots,z_{2k})\,e^{\frac{x}{2}\sum_{i=1}^k z_i-z_{k+i}}}
    {\prod_{i,j=1}^k (z_i-z_{k+j}) \prod_{i=1}^{2k}z_i^{2k}}
    \,dz_1\ldots dz_{2k}
\end{eqnarray}
(a more general function $p_k(x,\alpha)$ will be introduced in the next section).
Our basic claim is that the approximation
\begin{eqnarray}
    \label{eq:P approx}
    P_k(x)\approx a_k \, p_k(x+\tau_k,0)\,,
\end{eqnarray}
obtained from $P_k(x)$ by truncating the Taylor expansion of $\log
R(z_1,\ldots,z_{2k})$ at the linear term, is good enough to deduce asymptotics
for the coefficients $\{c_r(k),\, 0\le r\le k^{\beta}\}$, for any fixed $\beta
< 1$, in the sense the leading term asymptotic of the coefficient of
$x^{k^2-r}$, $0\le r \le k^{\beta}$, on either side of~\eqref{eq:P approx} is
the same.

Notice the formula defining $p_k(x,0)$ does not involve the complicated arithmetic factor
$A(z_1,\ldots,z_{2k})$ present in the residue expression for $P_k(x)$.
Moreover, by the results of Conrey, Farmer, Keating, Rubinstein, and Snaith,
 the function $p_k(x+\tau_k,0)$ can be evaluated explicitly as a
polynomial in $x$ of degree $k^2$. For, by property (\ref{eq:pk prop}) later,
 and the formulas in \textsection{2.7} of~\cite{CFKRS1}, we have
\begin{eqnarray}
    \label{eq:Lk_tau}
    p_k(x+\tau_k,0)=\frac{g_k}{k^2!}(x+\tau_k)^{k^2}\,.
\end{eqnarray}

The idea that the linear term in the Taylor expansion of $\log R(z_1,\ldots,z_{2k})$ ought 
 to dominate over $0\le r\le k^{\beta}$ was inspired, in part, by the analogous 
  asymptotic (\ref{eq:unitary_asymp}), derived in \cite{HR},  
  for the moments of the characteristic polynomial of random unitary matrices.  

As mentioned earlier, the main theorem of this paper shows that
the coefficients of the polynomial  
$a_k p_k(x+\tau_k,0) = \frac{a_k g_k}{k^2!}(x+\tau_k)^{k^2}$  
provide the leading asymptotics, as $k \to \infty$,
for essentially the first $k$ coefficients of $P_k(x)$. 
The proof of this theorem will naturally
 split into two main parts. In the first part, which is presented in 
 \textsection{\ref{sec: symm}}, \textsection{\ref{sec: Nk alg section}}, and
 \textsection{\ref{sec: consq of the algs}}, 
 we obtain estimates on certain functions in $k$, later
 denoted by $p_k$. In the second part, which is presented in \textsection{\ref{sec: arith factor}}, 
  we obtain bounds on the Taylor coefficients of the logarithm of the arithmetic factor. 
The latter bounds (and in some cases asymptotics) are fairly involved but generally straightforward, 
 while the former bounds are more subtle, requiring somewhat 
more thought.  Both bounds are obtained via  essentially combinatorial arguments. 

\subsection{Numerical verifications and an application to the maximal size of $|\zeta(1/2+it)|$.}
\label{sec: numerical ver}

Table 1 provides numerical confirmation of our Main Theorem, listing values of the ratio
\begin{equation}
    \frac{c_r(k)}{c_0(k)\binom{k^2}{r}\tau_k^r}
\end{equation}
for $k=10,20,30,40,50$ and $0 \leq r \leq 7$.
Our theorem provides an estimate for this ratio of $1+O\left( (r/k)^2\right)$,
and our table is consistent with such a remainder term, agreeing, for example, to 3-4 decimal places for
$r=2$ and $k=50$, and 2-3 decimal places for $r=8$ and $k=50$.

\begin{table}
    \begin{footnotesize}
        \centerline{
        \begin{tabular}{|c|c|c|c||}
            \hline
            $k$ & $r$ & $c_r(k)$ & $c_r(k)/\left(c_0(k)\binom{k^2}{r}\tau_k^r\right)$  \\
            \hline
            10 & 0 & 3.548884925e-148 &          1 \\
            10 & 1 & 2.357691331e-144 &          1 \\
            10 & 2 & 7.702336630e-141 & 0.9934255388 \\
            10 & 3 & 1.649486344e-137 & 0.9803060865 \\
            10 & 4 & 2.604519447e-134 & 0.9608017974 \\
            10 & 5 & 3.233666778e-131 & 0.9352015310 \\
            10 & 6 & 3.287651416e-128 & 0.9039165203 \\
            10 & 7 & 2.814729470e-125 & 0.8674698258 \\
            \hline 
            20 & 0 & 9.404052083e-789 &          1 \\
            20 & 1 & 7.007560591e-784 &          1 \\
            20 & 2 & 2.600909647e-779 & 0.9986738069 \\
            20 & 3 & 6.410977573e-775 & 0.9960221340 \\
            20 & 4 & 1.180624032e-770 & 0.9920509816 \\
            20 & 5 & 1.732651855e-766 & 0.9867716274 \\
            20 & 6 & 2.110801042e-762 & 0.9802005819 \\
            20 & 7 & 2.195579847e-758 & 0.9723595087 \\
            \hline 
            30 & 0 & 2.174528185e-2019 &          1 \\
            30 & 1 & 6.409313254e-2014 &          1 \\
            30 & 2 & 9.429995281e-2009 & 0.9994621075 \\
            30 & 3 & 9.234275546e-2004 & 0.9983864033 \\
            30 & 4 & 6.770756592e-1999 & 0.9967738368 \\
            30 & 5 & 3.964993050e-1994 & 0.9946262257 \\
            30 & 6 & 1.931729883e-1989 & 0.9919462534 \\
            30 & 7 & 8.053463103e-1985 & 0.9887374636 \\
            \hline 
            40 & 0 & 1.878520688e-3887 &          1 \\
            40 & 1 & 1.450126078e-3881 &          1 \\
            40 & 2 & 5.592030026e-3876 & 0.9997132915 \\
            40 & 3 & 1.436301603e-3870 & 0.9991398909 \\
            40 & 4 & 2.764308226e-3865 & 0.9982800615 \\
            40 & 5 & 4.252265871e-3860 & 0.9971343131 \\
            40 & 6 & 5.445979160e-3855 & 0.9957034019 \\
            40 & 7 & 5.972928889e-3850 & 0.9939883295 \\
            \hline 
            50 & 0 & 3.461963190e-6425 &          1 \\
            50 & 1 & 5.605367518e-6419 &          1 \\
            50 & 2 & 4.535291006e-6413 & 0.9998231027 \\
            50 & 3 & 2.444917857e-6407 & 0.9994693125 \\
            50 & 4 & 9.879474579e-6402 & 0.9989387280 \\
            50 & 5 & 3.191850197e-6396 & 0.9982315414 \\
            50 & 6 & 8.588531004e-6391 & 0.9973480389 \\
            50 & 7 & 1.979690769e-6385 & 0.9962886003 \\
            \hline
        \end{tabular}
        }
        \caption{A comparison of our asymptotic formula for $c_r(k)$, for
                 $k=10,20,30,40,50$ and $r\leq 7$. The 1's are explained by the
                 remark following the Main Theorem that the asymptotic formula
                 is actually an identity for $r=0$ and $r=1$.
                 We expect there to be lower terms in our asymptotic expansion,
                 and will return to the problem of determining them in a future paper.
        }\label{table:k 50}
    \end{footnotesize}
\end{table}


Next, let $\beta<1$, and, as usual, $k \in \mathbb{Z}_{\geq 0}$.
While the asymptotic formula for $c_r(k)$ given in our Main
Theorem holds, as $k \to \infty$, for $r<k^\beta$, it appears, numerically,
that our asymptotic formula is, uniformly, an upper bound for $|c_r(k)|$ for
all $0 \leq r \leq k^2$.

\begin{figure}
    \centerline{
        \includegraphics[width=5.2in]{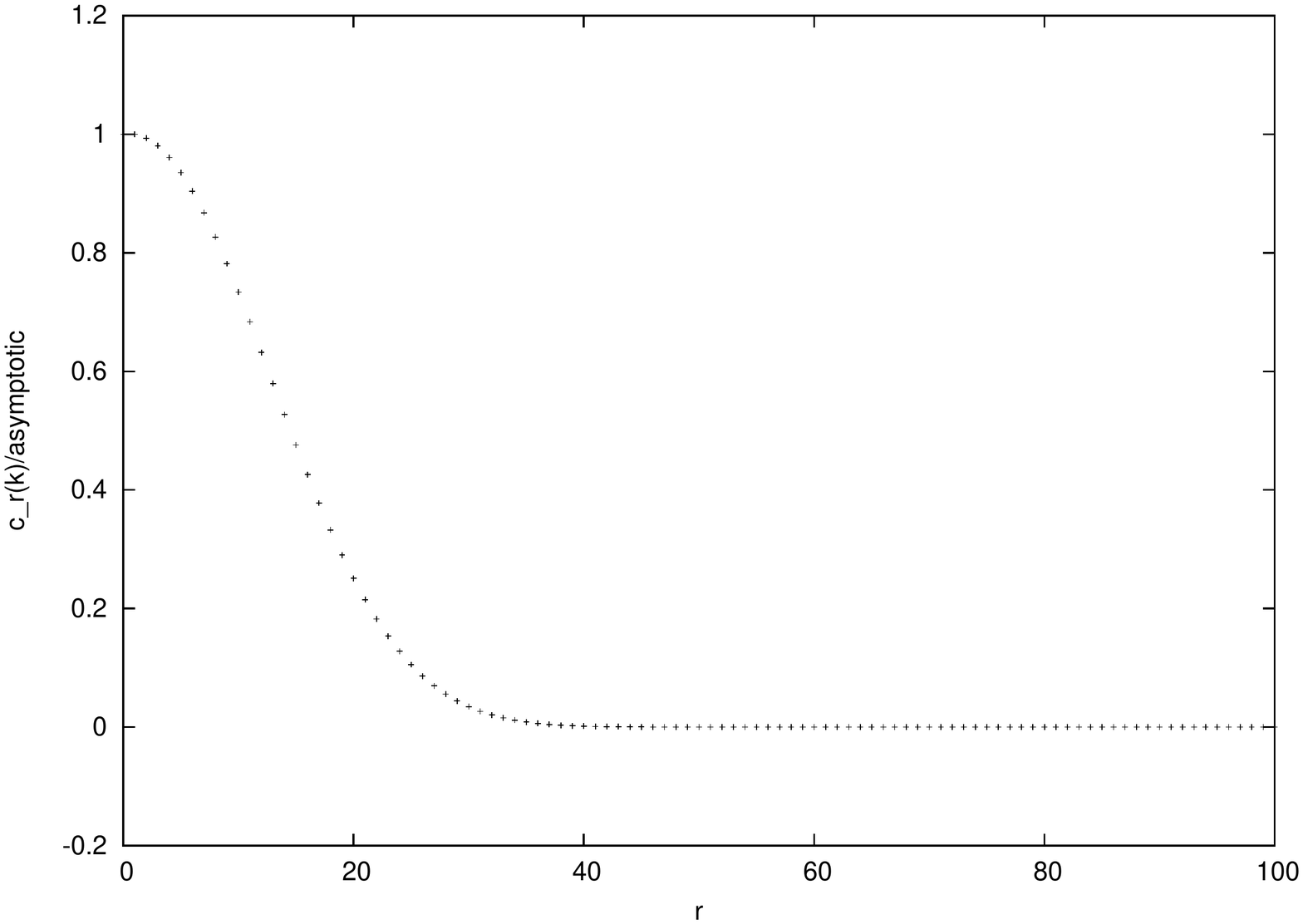}
    }
    \caption{We compare the ratio of $c_r(10)$, $0 < r < 100$, to our asymptotic formula.
    Here, $k=10$ is relatively small, and we only get reasonable agreement for the first few $r$.
    However, the graph indicates that the asymptotic formula is, uniformly, an upper bound
    for $|c_r(k)|$.
    }
    \label{fig:k 10}
\end{figure}

We therefore conjecture, for all non-negative integers $k$, and all
$0 \leq r \leq k^2$, that:
\begin{equation}
    \label{eq:c_r bound}
    |c_r(k)| \leq c_0(k)\binom{k^2}{r}\tau_k^r.
\end{equation}
We have verified this conjecture numerically for all $k \leq 13$, $0 \leq r \leq k^2$,
and all $k \leq 64$, $0 \leq r \leq 8$. The coefficients of the moment polynomials were
computed in the former case in~\cite{RY} and in the latter case using the program
developed for the computations in~\cite{CFKRS1} and~\cite{CFKRS2}.
See Figure 1 for evidence supporting this conjecture,
which depicts the ratio $c_r(k)/\left(c_0(k)\binom{k^2}{r}\tau_k^r\right)$
for $k=10$ and $0 \leq r \leq k^2$.

Assuming the bound~\eqref{eq:c_r bound},
we have, by the binomial theorem and term-wise comparison, the following upper
bound for $P_k(x)$, for all $k \in \mathbb{Z}_{\geq 0}$ and $x \in \mathbb{R}$:
\begin{equation}
    \label{eq:pkx bound}
    |P_k(x)| \leq c_0(k)(|x|+\tau_k)^{k^2}.
\end{equation}

Let $|\zeta(1/2+it_0)|=m_T:=\max_{t\in [0,T]} |\zeta(1/2+it)|$.
Lemma~3.3 of~\cite{FGH} provides:
\begin{equation}
    m_T \leq 2(CT \log T)^{1/2k}  \left(\frac{1}{T} \int_0^T
    |\zeta(1/2+it)|^{2k} dt \right)^{1/2k}
\end{equation}
for some absolute constant $C>0$.
Farmer, Gonek and Hughes use this inequality, combined with the Keating and Snaith leading term
conjecture for the moments of zeta to estimate $m_T$.
However, the leading term does a poor job at bounding the true size of the moments
if we allow $k$ to grow with $T$.

However, using our conjectured bound~\eqref{eq:pkx bound} for $P_k(x)$, we have,
in whatever range of $k$ that~\eqref{eq:Pkx} remains valid asymptotically, that
\begin{equation}
    \label{eq:m_T bound}
    m_T \leq 2(c_0(k) C_2 T \log T)^{1/2k}  \left(\frac{1}{T} \int_0^T
    (|\log(t/2\pi)| +\tau_k)^{k^2} dt \right)^{1/2k}
\end{equation}
for some absolute constant $C_2>0$. Following the argument in~\cite{FGH}, we
will, at the end, apply the above with $k$ proportionate to
$(\log(T)/\log\log(T))^{1/2}$.

The portion of the integral, $t \in (0,2\pi)$ where $\log(t/2\pi)$ is negative
contributes $O\left( (k^2)^{k^2}\right)$, on using: $\int_0^{2\pi} |\log(t/2\pi)
|^{k^2}dt = 2\pi\, k^2!$, the binomial expansion, Stirling's formula for
$k^2!$, and also $\sum_0^{k^2} \tau_k^r/r! < \exp(\tau_k)$ combined
with~\eqref{eq:tau_k bound}. (We could also slightly modify the argument
in~\cite{FGH} and ignore this interval outright.)

Next, by~\eqref{eq:tau_k bound}, we have $\tau_k=O(k\log k)$. Thus,
if $k \leq  C_3 \log(T)/\log\log(T)$, for some absolute constant $C_3$,
the contribution to the integral for
$t \in [2\pi,T]$ is $O\left(T (C_4 \log(T))^{k^2} \right)$, for some absolute
constant $C_4$.

Therefore, if $k=O(\log(T)^{1/2})$, we can ignore the portion of the integral from
$0$ to $2\pi$, and get:
\begin{equation}
    m_T \ll 2(c_0(k) C_5 T \log T)^{1/2k}  (C_4 \log{T})^{k/2}
\end{equation}
for some absolute constant $C_5$, i.e.
\begin{equation}
    \log{m_T} \ll \frac{\log{c_0(k)}}{2k} + \frac{\log(T) + \log\log(T)}{2k} +\frac{k}{2} \log\log{T} + O(k).
\end{equation}
Combining  Conrey and Gonek's estimate \cite{CG}:
\begin{equation}
    \log a_k \sim -k^2\log(2e^\gamma\log k) + o(k^2) \textrm{ for } k \rightarrow \infty,
\end{equation}
with the asymptotics of the Barnes $G$-function, see (3.17) and (3.18)
of~\cite{FGH}, gives:
\begin{equation}
    \frac{\log{c_0(k)}}{2k} = -\frac{k \log{k}}{2} +O(k \log\log{k}).
\end{equation}
Hence,
\begin{equation}
    \log{m_T} \ll \frac{\log(T) + \log\log(T)}{2k} +\frac{k}{2} \log\log{T} - \frac{k\log{k}}{2} + O(k \log\log{k}),
\end{equation}
i.e. bound (3.20) of~\cite{FGH} continues to hold even when we use our upper bound for the moment polynomials, rather
than the much smaller and less precise (as $k$ grows) leading term.

Taking, as in~\cite{FGH}, $k \sim c (\log(T)/\log\log(T))^{1/2}$, and choosing the optimal $c=2^{1/2}$, thus
gives the identical upper bound (3.9) of~\cite{FGH}:
\begin{equation}
    m_T \ll \exp\left(  (\tfrac12 \log T \log\log T)^{1/2} +
    O\left(\frac{(\log T)^{1/2} \log\log\log T}{(\log\log T)^{1/2}}\right)\right).
\end{equation}

Table 2 compares values of  $\int_0^T |\zeta(1/2+it)|^{2k}dt$, for $T=100000000.643$, $k\leq 13$,
to: the Keating and Snaith leading term $c_0(k) T \log(T)^{k^2}$ prediction, the full
asymptotics $\int_0^T P_k(\log(t/2\pi))dt$, and, finally, using
our upper bound for $P_k(x)$, i.e. to $c_0(k) \int_0^T (|\log(t/2\pi)| + \tau_k)^{k^2}dt$.

The values for the third column in Table 2 come from~\cite{RY}, and the lower
accuracy for $k=11,12,13$ reflects the precision to which we computed,
in~\cite{RY}, the coefficients of the moment polynomials.
The numerical integration of the moments of zeta was carried out in~\cite{RY}
using tanh-sinh quadrature, integrating the humps between successive zeros of
zeta on the critical line, hence we stopped at $100000000.643$ rather than
$10^8$.

The values in the 4th and 5th columns are given with more precision as they
only rely on $c_0(k)$ and $c_1(k)$ which have been computed to higher accuracy.
The table shows, first, that the full moment conjecture successfully captures,
here, the moments well beyond $k=4 \approx (2\log(T)/\log\log(T))^{1/2}$. It
also shows that the leading term alone quickly (for example, at $k=4$) fails to
capture the true size of the moments, whereas, our upper bound for the moment
polynomials seems to give an upper bound for the moments of zeta valid for a large
range of $k$, hence justifying its use in bounding the maximum size of zeta,
$m_T$.

\begin{table}
    \begin{footnotesize}
        \centerline{
        \begin{tabular}{|c|c|c|c|c|}
            \hline
            $k$ &
            $\int_0^T |\zeta(1/2+it)|^{2k}dt$ &
            $\int_0^T P_k(\log(t/2\pi))dt$ &
            $c_0(k) T \log(T)^{k^2}$ &
            $c_0(k) \int_0^T (|\log(t/2\pi)| + \tau_k)^{k^2}dt$ \\
            \hline
            1 &  1.6737236904e+09 & 1.6737234985e+09 & 1.8420680869e+09 & 1.6737235247e+09 \\
            2 &  6.3738834341e+11 & 6.3738992350e+11 & 5.8330132790e+11 & 6.7489927655e+11 \\
            3 &  8.0458531434e+14 & 8.0458140334e+14 & 1.3940397179e+14 & 1.2999952534e+15 \\
            4 &  1.7376480696e+18 & 1.7374512576e+18 & 4.3322247610e+15 & 8.5349032584e+18 \\
            5 &  5.0837678819e+21 & 5.0816645028e+21 & 6.0772270922e+15 & 1.8070544717e+23 \\
            6 &  1.8153019937e+25 & 1.8136396872e+25 & 1.8242195930e+14 & 1.2033327456e+28 \\
            7 &  7.4805129691e+28 & 7.4688841259e+28 & 6.5819531631e+10 & 2.4552753344e+33 \\
            8 &  3.4385117285e+32 & 3.4309032713e+32 & 1.7844629682e+05 & 1.4940783176e+39 \\
            9 &  1.7238857795e+36 & 1.7191846566e+36 & 2.4462083265e-03 & 2.6420504382e+45 \\
            10 &  9.2785048601e+39 & 9.2517330046e+39 & 1.2040915381e-13 & 1.3256809885e+52 \\
            11 &  5.2991086420e+43 & 5.28630715e+43 & 1.5747149879e-26 & 1.8471999998e+59 \\
            12 &  3.1825481927e+47 & 3.17945e+47 & 4.1820123844e-42 & 7.0111752824e+66 \\
            13 &  1.9956246380e+51 & 2.00e+51 & 1.7694787451e-60 & 7.1249837060e+74 \\
            \hline
        \end{tabular}
        }
        \caption{A comparison of three estimates for the moments of zeta,
        with $T=100000000.643$, and $k\leq 13$. The second and third columns
        are taken from~\cite{RY}.
        }\label{table:comparison}
    \end{footnotesize}
\end{table}

\section{Proof of the main theorem}

In the remainder of the paper, asymptotic constants are always absolute, and are taken as $k\to \infty$,
unless otherwise is stated. 

\begin{proof}[Proof of the main theorem]
    Let $\alpha:=(\alpha_1,\ldots,\alpha_{2k})$ be a $2k$-tuple in $\mathbb{Z}_{\geq 0}^{2k}$, and let $|\alpha|:=\alpha_1+\cdots+\alpha_{2k}$ denote its weight. Write
\begin{eqnarray}\label{eq:a_alpha_def}
    \log A(z_1,\ldots,z_{2k})=:\log a_k+B_k\sum_{i=1}^k z_i-z_{k+i}+\sum_{|\alpha|>1} a_{\alpha}\, z_1^{\alpha_1}\ldots z_{2k}^{\alpha_{2k}}\,,
\end{eqnarray}
the second sum being over tuples with weight greater than 1. Also, write
\begin{eqnarray} \label{eq:b_alpha_def}
    \log\left(\prod_{i,j=1}^k (z_i-z_{k+j})\zeta(1+z_i-z_{k+j})\right)=:
    \gamma k\sum_{i=1}^k z_i-z_{k+i}+\sum_{|\alpha|>1} b_{\alpha}\,
    z_1^{\alpha_1}\ldots z_{2k}^{\alpha_{2k}}\,.
\end{eqnarray}
The linear term in the Taylor expansion (\ref{eq:b_alpha_def}) is $\gamma k$, 
 which is an easy consequence of the expansion $z\zeta(1+z) = 1 + \gamma z +\cdots$. 
Lastly, define
\begin{eqnarray}
    \label{eq:p_k def}
    p_k(x,\alpha):=\frac{(-1)^k}{k!^2}\frac{1}{(2\pi i)^{2k}}\oint \cdots \oint
    \frac{\Delta^2(z_1,\ldots,z_{2k})\,e^{\frac{x}{2}\sum_{i=1}^k
    z_i-z_{k+i}}}{\prod_{i,j=1}^k (z_i-z_{k+j})\,\prod_{i=1}^{2k}z_i^{2k}}
    \,z_1^{\alpha_1}\ldots z_{2k}^{\alpha_{2k}}\,dz_1\ldots dz_{2k}\,, \nonumber \\
    \,
\end{eqnarray}
and let $c_{\alpha}$ be the Taylor coefficients determined by 
\begin{eqnarray}
    \label{eq:c_alpha_def}
    e^{\sum_{|\alpha|>1} (a_{\alpha}+b_\alpha)\,z_1^{\alpha_1}\ldots
    z_{2k}^{\alpha_{2k}}}=:1+\sum_{|\alpha|>1}c_{\alpha}\,z_1^{\alpha_1}\ldots
    z_{2k}^{\alpha_{2k}}\,.
\end{eqnarray}
So, on recalling $\tau_k =2B_k+2\gamma k$, the $c_{\alpha}$'s satisfy:
\begin{eqnarray}
    A(z_1,\ldots,z_{2k}) \prod_{i,j=1}^k (z_i-z_{k+j})\zeta(1+z_i-z_{k+j}) = a_k e^{\frac{\tau_k}{2}\sum_{i=1}^k z_i-z_{k+i}} 
    \left( 1+ \sum_{|\alpha|>1} c_{\alpha} z_1^{\alpha_1}\ldots z_{2k}^{\alpha_{2k}}\right)\,,
\end{eqnarray}
where, as before, $\tau_k \sim 4k\log k$ as $k\to \infty$. Therefore, we have
\begin{eqnarray}
    \label{eq:Pk to pk}
    P_k(x) = a_k\, p_k(x + \tau_k, 0) + a_k \sum_{|\alpha|>1} c_{\alpha} \,p_k(x +\tau_k,\alpha)\,,
\end{eqnarray}
where the second argument in $p_k(x+\tau_k,0)$ stands for the zero $2k$-tuple.

Notice the sum in (\ref{eq:Pk to pk}) is actually finite,
because  if $|\alpha|>k^2$ (or if $\alpha_j \ge 2k$ for some $j$), 
then $p_k(x,\alpha) = 0$, because by degree considerations
the integrand in the residue (\ref{eq:p_k def}) defining $p_k(x,\alpha)$  
will have no poles. Also, by the change of variables, $z_j \leftarrow x z_j$, we have
\begin{eqnarray}
    \label{eq:pk prop}
    p(x,\alpha)=x^{k^2-|\alpha|}\, p(1,\alpha)\,,
\end{eqnarray}
which, along with the formulas in \textsection{2.7} of~\cite{CFKRS1}, yields 
\begin{eqnarray} \label{eq:ksb}
    p_k(x,0)= x^{k^2} \, p_k(1,0) = x^{k^2} \,\frac{g_k}{k^2!}\,.
\end{eqnarray}
(We used formulas (\ref{eq:pk prop}) and (\ref{eq:ksb}) to evaluate $p_k(x+\tau_k,0)$ in (\ref{eq:Lk_tau}) earlier).
In light of property (\ref{eq:pk prop}),  it is convenient to set
\begin{eqnarray}\label{eq:p_k def1}
    p_k(\alpha):=p_k(1,\alpha)\,. 
\end{eqnarray}
Combining (\ref{eq:Pk to pk}), the observation made thereafter, and (\ref{eq:pk prop}), we arrive at
\begin{eqnarray}
    \label{eq:equate coeff}
    P_k(x) =a_k\, \left(x+\tau_k \right)^{k^2} \, p_k(0) + a_k \, \sum_{n=2}^{k^2}
    \left(x+\tau_k\right)^{k^2-n} \sum_{|\alpha|=n} c_{\alpha}\,
    p_k(\alpha)\,.
\end{eqnarray}
In particular, observing $a_k \,p_k(0) = c_0(k)$, and equating the coefficient
 of $x^{k^2-r}$ on both sides of (\ref{eq:equate coeff}), we obtain
\begin{eqnarray}
    \label{eq:crk identity}
    c_r(k)&=& \tau_k^r\, \binom{k^2}{r}\,c_k(0)+ a_k\,\sum_{n=2}^{r}
    \tau_k^{r-n}\,\binom{k^2-n}{r-n} \sum_{|\alpha|=n} c_{\alpha} p_k(\alpha)\nonumber\\
     &=& \tau_k^r\,\binom{k^2}{r}\, c_k(0) \left[1+\sum_{n=2}^{r}
    \frac{r!\,(k^2-n)!}{(r-n)! \,k^2!} \, \frac{1}{\tau_k^n} \sum_{|\alpha|=n}
    c_{\alpha} \frac{p_k(\alpha)}{p_k(0)} \right]\,. \label{eq:main identity}
\end{eqnarray}
The above is an identity, valid for any $0\le r\le k^2$. Also, notice
the double sum in \eqref{eq:crk identity} is empty if 
$r=0,1$, so $c_r(k)= \tau_k^r\,\binom{k^2}{r}\,c_0(k)$
 for $r=0,1$.

Our plan is to show, for $0\le r\le k^{\beta}$, $c_r(k)\approx \tau_k^r\binom{k^2}{r} c_0(k)$.
 To do so, we will show that the term 1 preceding the
 double sum in (\ref{eq:main identity}) dominates. 
 This will follow from the following three bounds, 
as we soon explain:

\begin{itemize}
\item First bound: By Theorem~\ref{thm:nkbd}, as $k\to \infty$ and uniformly in $|\alpha|<k/2$, we have 
\begin{eqnarray} \label{eq:pk bound}
    \frac{p_k(\alpha)}{p_k(0)}  \ll (\lambda_1 k \log (|\alpha| + 10))^{|\alpha|}\,,
\end{eqnarray}
where $\lambda_1$ is some absolute constant. 
This is proved in \textsection{\ref{sec: consq of the algs}} as a by-product of the ``symmetrization algorithm''
(see \textsection{\ref{sec: symm}}), and the algorithm 
to compute a certain  ``symmetrized version'' of $p_k(\alpha)$,
which we denote $N_k^0(\alpha)$ (see \textsection{\ref{sec: the alg to compute N}}).
 The notation $N_k^0(\alpha)$ is chosen to distinguish it from the related function
 $N_k(\alpha)$, defined in \cite{CFKRS2}.
 The said algorithms are essentially combinatorial recursions. In the case 
 of $N^0_k(\alpha)$, the recursion stops 
  much earlier than what is obvious, due to a certain anti-symmetry relation,  
  which is the reason algorithm is able to produce a non-trivial 
   bound on $N_k^0(\alpha)$, essentially by counting the number of terms involved in it.   
We remark the bound \eqref{eq:pk bound} is sharp in the power of $k$, as the second
example in \textsection{\ref{sec: Nk alg examples}} illustrates.  

\item Second bound: By Theorem~\ref{arithmfac}, the coefficients 
$a_{\alpha}$ in the Taylor expansion of $\log A(z_1,\ldots,z_{2k})$, 
which were defined in (\ref{eq:a_alpha_def}), satisfy:
\begin{eqnarray} \label{eq:a alpha bound}
    a_{\alpha}\ll  \lambda_2^{|\alpha|}\,(\log k)^{|\alpha|}\,\left[
    m(\alpha)^{|\alpha|} \,k^{2-\min\{m(\alpha),2\}} +  |\alpha|!\,
    k^{2-m(\alpha)}\right]\,,
\end{eqnarray}
where $m(\alpha)$ denotes the number of non-zero entries in $\alpha$, and 
 $\lambda_2$ is some absolute constant. 
This is proved in \textsection{\ref{sec: arith factor}} by an elementary, though lengthy,
counting of the terms that contribute. It will transpire that, for $0\le r\le k^{\beta}$, most of the contribution 
 to $c_r(k)$ comes from  ``the combinatorial sum for the small primes'', 
 see \textsection{\ref{sec: small primes comb}}.

\item Third bound: By lemma~\ref{lem:prodz}, the Taylor coefficients $b_{\alpha}$
 of the product of zetas, which were defined in (\ref{eq:b_alpha_def}), satisfy:

\begin{eqnarray}\label{eq:b alpha bound}
    b_{\alpha}\ll \lambda_3^{|\alpha|}\,k^{2-m(\alpha)}\,.
\end{eqnarray}
This is proved in \textsection{\ref{sec: zeta product}} by means of Cauchy's estimate. 
\end{itemize}

We now appeal to the auxiliary lemma stated later in this section. 
Specifically, by (\ref{eq:a alpha bound}) and (\ref{eq:b alpha bound}), the coefficients
 $a_{\alpha} + b_{\alpha}$ still satisfy the conditions of that lemma. 
 So on applying the lemma we obtain the following bound on the Taylor coefficients $c_{\alpha}$, 
 which were defined in (\ref{eq:c_alpha_def}): 
As $k\to \infty$, and uniformly in $n<k/e$, 
\begin{eqnarray}
    \label{eq:c alpha bound}
    \sum_{|\alpha|=n}|c_{\alpha}|\ll (\lambda_4\,k\log k)^n\,.
\end{eqnarray}
Notice the number of summands on the lhs above is not
 far off from the upper bound, so, on average, the 
  $|c_{\alpha}|$'s are not large when $|\alpha|<k/e$.

Substituting \eqref{eq:pk bound} and \eqref{eq:c alpha bound} 
directly into identity \eqref{eq:main identity}, and recalling $r\le k^{\beta}$, yields
\begin{eqnarray}
    \label{eq:another bound1}
    \sum_{n=2}^{r} \frac{r!\,(k^2-n)!}{(r-n)! \,k^2!} \, \frac{1}{\tau_k^n}
    \sum_{|\alpha|=n} \left|c_{\alpha} \frac{p_k(\alpha)}{p_k(0)}\right|&\ll&
    \sum_{n=2}^{r} \frac{r^n}{k^{2n}\tau_k^n} (\lambda_1\, k \, \log k)^n (\lambda_4\,k\log (n+10))^n \nonumber\\
    &\ll& \sum_{n=2}^r \frac{(\lambda\,r\, \log n)^n }{k^n}\,, \label{eq:bound1}
\end{eqnarray}
for some absolute constant $\lambda$.
Here, we used the following elementary bound
\begin{eqnarray}
    \label{eq:easy bound}
    \frac{r!\,(k^2-n)!}{(r-n)! \,k^2!}\le \frac{r^n}{k^{2n}}\,,
\end{eqnarray}
which follows from $(r-j)/(k^2-j)=(r/k^2)(1-j/r)/(1-j/k^2) \leq r/k^2$ with $j\leq(n-1) < r$, and
$r< k^2$ (in fact, $r<k$ in this proof).

Finally, summing the series in (\ref{eq:bound1}), 
and using the assumed bound on $r$, shows that the lhs of (\ref{eq:another bound1}) 
is bounded by $O_{\beta}\left( (r /k)^2 \right)$, completing the proof.
\end{proof}

\begin{aux}\label{aux}
Let $f$ be a multi-variate series in $2k$ variables
\begin{eqnarray}
    \label{eq:f_l}
    f(x_1,\ldots,x_{2k}):=\sum_{n=2}^{\infty}
    \sum_{\substack{\,\,\,\alpha\in \mathbb{Z}_{\geq 0}^{2k}\\|\alpha|=n}}
    a_{\alpha}\, x_1^{\alpha_1}\ldots x_{2k}^{\alpha_{2k}}\,.
\end{eqnarray}
Assume the coefficients $a_{\alpha}$ satisfy bounds (\ref{eq:a alpha bound}).
Then the coefficients $c_{\alpha}$ in the Taylor expansion
\begin{eqnarray}
    e^{f(x_1,\ldots,x_{2k})}=:1+\sum_{n=2}^{\infty}
    \sum_{|\alpha|=n}
    c_{\alpha}\,x_1^{\alpha_1}\ldots x_{2k}^{\alpha_{2k}}
\end{eqnarray}
satisfy
\begin{eqnarray}
    \label{eq:c alpha bound 2}
    \sum_{\substack{|\alpha|=n}} |c_{\alpha}|\ll (\lambda_5\, \log k)^n\,
    k^n\,,\qquad \textrm{for $n<k/e$} \,, 
\end{eqnarray}
for some absolute constant $\lambda_5$.
\end{aux}
Remarks: i) This lemma applies as well if we replace $a_\alpha$ by $a_\alpha+b_\alpha$,
with $b_\alpha$ satisfying~\eqref{eq:b alpha bound}, because $a_\alpha+b_\alpha$ together
satisfy a bound of the same form as~\eqref{eq:a alpha bound}, but with $\lambda_2$
replaced by the maximum of $\lambda_2$ and $\lambda_3$. ii) We are using this lemma
 in~\eqref{eq:c alpha bound}.

\begin{proof}

Define
\begin{eqnarray}
    C(n):=\sum_{\substack{|\alpha|=n}} |c_{\alpha}|\,,\qquad
    A(q):=\sum_{\substack{|\alpha|=q}} |a_{\alpha}|\,.
\end{eqnarray}
We plan to obtain a bound on $C(n)$ in terms of an expression involving 
 $A(q)$, then we will bound $A(q)$ with the aid of estimate (\ref{eq:a alpha bound})
 for the $a_{\alpha}$'s, which is assumed in the statement of the lemma.

To this end, exponentiate (\ref{eq:f_l}), turning the outer sum into a product,
and writing, for the inner sum,
\begin{equation}
    \label{eq:inner sum}
    \exp\left(
        \sum_{|\alpha|=n}
        a_{\alpha}\, x_1^{\alpha_1}\ldots x_{2k}^{\alpha_{2k}}
    \right)
    = \sum_{d=0}^\infty \frac{1}{d!}
           \left( \sum_{|\alpha|=n}
                a_{\alpha}\, x_1^{\alpha_1}\ldots x_{2k}^{\alpha_{2k}}
           \right)^d
    \,,
\end{equation}
we get, on multiplying out the product, that
\begin{equation}
    1+\sum_{n=2}^{\infty} \sum_{|\alpha|=n}
    c_{\alpha}\,x_1^{\alpha_1}\ldots x_{2k}^{\alpha_{2k}}
    = \prod_{n=2}^{\infty} 
     \sum_{d_n=0}^\infty \frac{1}{d_n!}
           \left( \sum_{|\alpha|=n}
                a_{\alpha}\, x_1^{\alpha_1}\ldots x_{2k}^{\alpha_{2k}}
           \right)^{d_n}
   \,. 
\end{equation}
By choosing which of the sums in the above infinite product contribute (i.e., 
 which of the sums has a term chosen from it different from 1), we obtain
\begin{eqnarray}
    \label{eq:Cn bound}
    C(n)\le \sum_{\substack{q_1d_1+\cdots+q_rd_r=n,\, r\geq 1\\ q_r>\cdots>q_2>q_1\ge 2,\,
    d_i\ge1}} \frac{1}{d_1!d_2!\ldots d_r!}
    \,A\left(q_1\right)^{d_1}A\left(q_2\right)^{d_2}\ldots
    A\left(q_r\right)^{d_r}\,.
\end{eqnarray}
We now derive a bound on the $A(q_j)$'s. Given an integer $2\le q\le n$, write
\begin{eqnarray}
    A(q)= \sum_{j=1}^{q} \sum_{\substack{|\alpha|=q\\ m(\alpha)=j}}
    |a_{\alpha}|=\sum_{\substack{|\alpha|=q\\ m(\alpha)=1}}
    |a_{\alpha}|+\sum_{j=2}^{q} \sum_{\substack{|\alpha|=q\\ m(\alpha)=j}}
    |a_{\alpha}|\,,
\end{eqnarray}
where, recall, $m(\alpha)$ is equal to the number of non-zero $\alpha_i$'s.
Substituting the bounds (\ref{eq:a alpha bound}) for the $|a_{\alpha}|$'s, we get
\begin{eqnarray}
    \label{eq:Aq sums}
    A(q)\ll \sum_{\substack{|\alpha|=q\\ m(\alpha)=1}} \left(\lambda_2\right)^q\,
    q!\, (\log k)^q\,k + \sum_{j=2}^q \sum_{\substack{|\alpha|=q\\ m(\alpha)=j}}
    \left(\lambda_2\right)^q\, j^q\,(\log k)^q  + \sum_{j=2}^q
    \sum_{\substack{|\alpha|=q\\ m(\alpha)=j}} \frac{\left(\lambda_2\right)^q\,q!\,
    (\log k)^q}{k^{j-2}}\,.\nonumber \\
    \,
\end{eqnarray}
But
\begin{eqnarray}
    \label{eq:Aq sums1}
    \sum_{\substack{|\alpha|=q\\ m(\alpha)=j}} 1= \binom{2k}{j}
    \binom{q-1}{j-1}\,,
\end{eqnarray}
as there are $\binom{2k}{j}$ ways to select $j$ of the $z_i$'s and $\binom{q-1}{j-1}$
ways to sum to $q$ using precisely $j$ positive (ordered) integers. The latter fact can be seen
by arranging $q$ `dots' in a row and breaking them into $j$ summands by
selecting $j-1$ out of $q-1$ barriers between the dots.

Therefore, for $q<k/2$ (for later purposes, we actually assume $q\le n < k/e$ in this proof), we have
generously,
\begin{eqnarray}
    \label{eq:Aq sums2}
    \sum_{j=2}^q \binom{2k}{j} \binom{q-1}{j-1}\,j^q\le \sum_{j=2}^q
    \frac{(2k)^j j^q q^j}{(j!)^2}\le k^q\,(100)^q\,.\label{eq:binomial sum}
\end{eqnarray}
The first inequality follows by expanding the binomial coefficients as ratios
of factorials, and noting that: i) $(2k)!/(2k-j)!\leq (2k)^j$.
ii) $j (q-1)!/(q-j)! \leq j q^{j-1} \leq q^j$. The second inequality
in~\eqref{eq:binomial sum} follows by
noticing that the terms of the sum are, in our range, increasing (consider the ratio of
two successive terms), hence an upper bound for sum is $q$ times the last term, which
can be estimated by Stirling's formula.
Similarly,
\begin{eqnarray}
    \label{eq:Aq sums3}
    \sum_{j=2}^q \binom{2k}{j} \binom{q-1}{j-1}\,q!\,k^{2-j}\le 2^q\,k^2
    \sum_{j=2}^q \frac{q^j q^q}{(j!)^2} \le k^2\,q^q\,(100)^q\,.
\end{eqnarray}

Using \eqref{eq:Aq sums2} to bound the second sum in \eqref{eq:Aq sums},
 using (\ref{eq:Aq sums3}) to bound the third sum, 
 and noting that the number of terms in the first sum there is
\begin{eqnarray}
    \label{eq:Aq sums0}
    \sum_{\substack{|\alpha|=q\\ m(\alpha)=1}} 1 = 2k\,,
\end{eqnarray}
which follows since there are $2k$ choices for the $z_j$'s, together yields
\begin{eqnarray}
    A(q)&\ll& k^2\, \left(\lambda_2\, q\, \log k\right)^q+ k^q\,
    \left(100\,\lambda_2\, \log k\right)^q+k^2\,\left(100\, \lambda_2\, q\, \log
    k\right)^q \\ 
    \, \nonumber \\ 
    &\ll& k^q\, \left(100\, \lambda_2 \, \log
    k\right)^q\, \left[1 + k^2 \left(\frac{q}{k}\right)^q\right]\,.
\end{eqnarray}
Substituting the above into (\ref{eq:Cn bound}), we obtain for some absolute constant $\lambda_6$,
\begin{eqnarray}
    \label{eq:Cn product bound}
    C(n)\ll k^n (\lambda_6  \log k)^{n}  \sum_{\substack{q_1d_1+\cdots+q_rd_r=n,\,r\ge 1\\
    q_r>\cdots>q_2>q_1\ge 2,\, d_i\ge 1}} \frac{1}{d_1!d_2!\ldots d_r!} \,
    \prod_{i=1}^r \left[1 + k^2 \left(\frac{q_i}{k}\right)^{q_i}\right]^{d_i}\,.
\end{eqnarray}

Since the function $(x/k)^x$ is monotonically decreasing for $x\in [1,k/e)$, it follows
\begin{equation}
    k^2 \left(\frac{q_i}{k}\right)^{q_i} \le  
    4\,, \quad \textrm{if $2\le q_i<k/e$}\,.
\end{equation}
Thus,
\begin{equation}
    \label{eq:product bound c}
    \prod_{i=1}^r \left[1 + k^2 \left(\frac{q_i}{k}\right)^{q_i}\right]^{d_i} \le 
    5^n,\qquad \textrm{if $2\le q_i<k/e$}\,.
\end{equation}
Here we have used $\sum d_i \leq n$.
Also,
\begin{eqnarray}
    \label{eq:e n bound}
    \sum_{\substack{q_1d_1+\cdots+q_rd_r=n,\, r\ge 1\\ q_r>\cdots>q_2>q_1\ge 1\,,\,
    d_i\ge1}} \frac{1}{d_1!d_2!\ldots d_r!}< e^n,. 
\end{eqnarray}
because the lhs is the coefficient of $x^n$ in $\prod_{m=1}^n \sum_{d=1}^\infty x^{md}/d!$
(we truncate the product at $m=n$ since each $q_i\leq n$). But that coefficient is less than
the sum total of all the coefficients, i.e. $<\prod_{m=1}^n \sum_{d=1}^\infty 1/d! < e^n$.

Substitute (\ref{eq:product bound c}) and (\ref{eq:e n bound}) into (\ref{eq:Cn product bound}), 
 we have, for $n<k/e$, 
\begin{eqnarray}
    C(n)&\ll& (5\,\lambda_6\, \log k)^n\,k^n\,
    \sum_{\substack{q_1d_1+\cdots+q_rd_r=n\\q_i\ge 2\,,\, d_i\ge1\,, r\ge 1}}
    \frac{1}{d_1!d_2!\ldots d_r!}  \\ &\ll&  (15\, \lambda_6\, \log k)^n\,
    k^n\,, 
\end{eqnarray}
as claimed.

\end{proof}

\section{An algorithm to reduce to the first half}
\label{sec: symm}
We show that the residue expression for $p_k(\alpha)$, given by
(\ref{eq:p_k def}) and (\ref{eq:pk prop}), can be
 reduced to variables in the first half only; 
 i.e., involving $z_1,\ldots,z_k$ only. 
 To do so, we will need the following two lemmas.

\begin{lem}\label{lem:sym1}
Suppose $H(z_1,\ldots,z_{2n})$ is regular in
\mbox{$\mathcal{D}:=\{|(z_1,\ldots,z_{2n})|<\delta\}$}. For
\mbox{$(\alpha_1,\ldots,\alpha_{2n})\in \mathcal{D}$}, such that the
$\alpha_i$'s are distinct, define
\begin{equation}
    \mathcal{K}(\alpha_1,\ldots,\alpha_{2n}):=\sum_{\sigma\in S_{2n}}
    \frac{H(\alpha_{\sigma(1)},\ldots,\alpha_{\sigma(2n)})}{\prod_{i,j=1}^n
    (\alpha_{\sigma(i)}-\alpha_{\sigma(n+j)})}\,,
\end{equation}
where $S_{2n}$ be the permutation group of $2n$ elements. Then, it holds
\begin{equation}
    \mathcal{K}(\alpha_1,\ldots,\alpha_{2n})= \frac{(-1)^n}{(2\pi i)^{2n}}\oint
    \cdots \oint
    \frac{H(z_1,\ldots,z_{2n})\,\Delta^2(z_1,\ldots,z_{2n})}{\prod_{i,j=1}^n
    (z_i-z_{n+j})\,\prod_{i,j=1}^{2n} (z_i-\alpha_j)} \,dz_1\ldots dz_{2k}\,,
\end{equation}
where the integration contour consists of circles contained in $\mathcal{D}$ around the $\alpha_i$'s.
In particular, if the integration contour is chosen so each circle encloses 0 as well, then the limit
\begin{equation}
    \lim_{\substack{\alpha_i\to 0\\1\le i\le 2n}}
    \,\mathcal{K}(\alpha_1,\ldots,\alpha_{2n})= \frac{(-1)^n}{(2\pi i)^{2n}}\oint
    \cdots \oint
    \frac{H(z_1,\ldots,z_{2n})\,\Delta^2(z_1,\ldots,z_{2n})}{\prod_{i,j=1}^n
    (z_i-z_{n+j})\,\prod_{i=1}^{2n} z_i^{2n}} \,dz_1\ldots dz_{2k}\,,
\end{equation}
exists, and is finite.
\end{lem}

\begin{proof}
    This lemma is a slight variant of lemmas 2.5.1 and 2.5.3 in~\cite{CFKRS1}.
\end{proof}

\begin{lem}\label{lem:sym2}
Let $H(z_1,\ldots,z_{2n})$ and $f(z_1,\ldots,z_{2n})$ be two regular functions
in $\mathcal{D}$. Suppose also $f$ is symmetric with respect to all its
arguments (so $f$ is invariant under the action of $S_{2n}$). Define
\begin{eqnarray}
    I(f):=\frac{(-1)^n}{(2\pi i)^{2n}}\oint \cdots \oint
    \frac{H(z_1,\ldots,z_{2n})\,
    f(z_1,\ldots,z_{2n})\,\Delta^2(z_1,\ldots,z_{2n})}{\prod_{i,j=1}^n
    (z_i-z_{n+j})\,\prod_{i=1}^{2n} z_i^{2n}} \,dz_1\ldots dz_{2k}\,,\nonumber \\
\end{eqnarray}
where the integration contour consists of circles in $\mathcal{D}$ around 0. Then
\begin{equation}
    I(f)=f(0,\ldots,0)\,I(1)\,.
\end{equation}

\end{lem}

\begin{proof}

Define
\begin{equation}
    \mathcal{K}_f(\alpha_1,\ldots,\alpha_{2n}):=\sum_{\sigma\in S_{2n}}
    \frac{H(\alpha_{\sigma(1)},\ldots,\alpha_{\sigma(2n)})\,f(\alpha_{\sigma(1)},\ldots,\alpha_{\sigma(2n)})}{\prod_{i,j=1}^n
    (\alpha_{\sigma(i)}-\alpha_{\sigma(n+j)})}\,.
\end{equation}
Then,
\begin{eqnarray}
    I(f)&=& \lim_{\substack{\alpha_i\to 0\\1\le i\le 2n}}
    \,\mathcal{K}_f(\alpha_1,\ldots,\alpha_{2n})\\ &=&
    \lim_{\substack{\alpha_i\to 0\\1\le i\le 2n}}
    f(\alpha_1,\ldots,\alpha_{2n})\,\lim_{\substack{\alpha_i\to 0\\1\le i\le
    2n}}\,\mathcal{K}_1(\alpha_1,\ldots,\alpha_{2n})\\ &=&
    f(0,\ldots,0)\, I(1)\,.  
\end{eqnarray}
\end{proof}

\subsection{The first step: from $p_k(\alpha)$ to $p_k(\lambda;0)$}

Recall, for a tuple $\alpha=(\alpha_1,\ldots,\alpha_{2k})\in \mathbb{Z}_{\geq 0}^{2k}$ we defined
\begin{eqnarray}
    \label{eq:pk sym}
    p_k(\alpha):=\frac{(-1)^k}{k!^2}\frac{1}{(2\pi i)^{2k}}\oint \cdots \oint
    \frac{\Delta^2(z_1,\ldots,z_{2k})\,e^{\frac{1}{2}\sum_{i=1}^k
    z_i-z_{k+i}}}{\prod_{i,j=1}^k (z_i-z_{k+j})\,\prod_{i=1}^{2k}z_i^{2k}}
    \,z_1^{\alpha_1}\ldots z_{2k}^{\alpha_{2k}}\,dz_1\ldots dz_{2k}\,.\nonumber \\
    \,
\end{eqnarray}
In this subsection we show that $p_k(\alpha)$ can always be written as a 
relatively short (for purposes of our analysis) linear combination
of functions of the form $p_k(\beta_1,\ldots,\beta_k,0,\ldots,0)$, where
$\beta_i\in \mathbb{Z}_{\geq 0}$ for all $1\le i\le k$.
So consider a $2k$-tuple $\alpha=(\alpha_1,\ldots,\alpha_{k+d},0,\ldots,0)$ where 
$1\le d\le k$, and such that $\alpha_{k+i}>0$ for $1\le i\le d$.
Since the integral (\ref{eq:pk sym}) is then symmetric in $z_{k+d},\ldots,z_{2k}$, it follows 
\begin{eqnarray}
    p_k(\alpha)=\frac{(-1)^k}{k!^2}\frac{1}{(2\pi i)^{2k}}\oint \cdots \oint
    \frac{\Delta^2(z_1,\ldots,z_{2k})\,e^{\frac{1 }{2}\sum_{i=1}^k
    z_i-z_{k+i}}}{\prod_{i,j=1}^k (z_i-z_{k+j})\,\prod_{i=1}^{2k}z_i^{2k}}
    \,z_1^{\alpha_1}\ldots z_{k+d-1}^{ \alpha_{k+d-1}} \,\times  \nonumber \\
    \,  \\ 
    \frac{1}{k-d+1}\,\left(\sum_{j=k+d}^{2k}
    z_j^{\alpha_{k+d}}\right)\,\,dz_1\ldots
    dz_{2k}\,,\qquad\qquad\qquad\qquad\nonumber
\end{eqnarray}
and by lemma~\ref{lem:sym2},
\begin{eqnarray}
    p_k(\alpha)=\frac{(-1)^k}{k!^2}\frac{1}{(2\pi i)^{2k}}\oint \cdots \oint
    \frac{\Delta^2(z_1,\ldots,z_{2k})\,e^{\frac{1 }{2}\sum_{i=1}^k
    z_i-z_{k+i}}}{\prod_{i,j=1}^k (z_i-z_{k+j})\,\prod_{i=1}^{2k}z_i^{2k}}
    \,z_1^{\alpha_1}\ldots z_{k+d-1}^{ \alpha_{k+d-1}} \,\times  \nonumber \\
    \,  \\ 
    \frac{1}{k-d+1}\,\left(\sum_{j=k+d}^{2k}
    z_j^{\alpha_{k+d}}-\sum_{j=1}^{2k} z_j^{\alpha_{k+d}}\right)\,\,dz_1\ldots
    dz_{2k}\,.\qquad\qquad\qquad\qquad\nonumber
\end{eqnarray}
This can be seen
from lemma~\ref{lem:sym2} by pulling out the second sum in brackets in front of
the integral, evaluated at all $z_j=0$, to give 0.
For $1\le j\le 2k$, let us thus define
\begin{eqnarray}
    \eta^{(j)}&:=&(\overbrace{0,\ldots,0}^{\textrm{$j-1$
    zeros}},\alpha_{k+d},0,\ldots,0)\,, \\ 
    \, \nonumber \\
    \alpha^{(j)}&:=&\alpha-\eta^{(k+d)}+\eta^{(j)}\,,
\end{eqnarray}
where the addition and subtraction in the definition of $\alpha^{(j)}$ is
 done component-wise. Then we have
\begin{eqnarray}
    p_k(\alpha)= \frac{-1}{k-d+1}\,\sum_{j=1}^{k+d-1} p_k(\alpha^{(j)})\,.
\end{eqnarray}
In particular, we have expressed $p_k(\alpha)$ as the sum of $k+d-1$
functions of the form $p_k(\beta)$, where each tuple $\beta$ has its last
possibly non-zero entry in position $k+d-1$ (instead of position $k+d$, as was the case
for $\alpha$ itself), and each $\beta$ satisfies $|\beta|=|\alpha|$.
By iterating this procedure several times, we obtain the following lemma.

\begin{lem}\label{lem:sym lem}
Let \mbox{$\alpha=(\alpha_1,\ldots,\alpha_{2k})\in \mathbb{Z}_{\geq 0}^{2k}$},
 and let $d$ be the number of non-zero entries in the second half of $\alpha$
(i.e. among the entries $\alpha_{k+1},\ldots,\alpha_{2k}$). 
Further, given $\lambda=(\lambda_1,\ldots,\lambda_k)\in \mathbb{Z}_{\geq 0}^k$, define
$p_k(\lambda;0):=p_k(\lambda_1,\ldots,\lambda_k,0,\ldots,0)$.
Then the function $p_k(\alpha)$ can be written in the form
\begin{equation}
    p_k(\alpha)=\frac{(-1)^d}{\prod_{j=1}^d (k-d+j)}\,\sum_{\lambda\in
    \mathcal{S}_{\alpha}} p_k(\lambda;0)\,,
\end{equation}
where $\mathcal{S}_{\alpha}$ is a certain set of tuples $\lambda \in
\mathbb{Z}_{\geq 0}^k$, with $|\lambda|=|\alpha|$, of cardinality
$|\mathcal{S}_{\alpha}|= \prod_{j=1}^d (k+d-j)$.
\end{lem}

\subsection{An example}

Given a tuple of the form 
\begin{eqnarray}
    (\alpha_1,\ldots,\alpha_l,0,\ldots,0,\alpha_{k+1},\ldots,\alpha_{k+d},0,\ldots,0)
    \in \mathbb{Z}_{\geq 0}^{2k}\,,
\end{eqnarray}
where the $\alpha_j$'s are possibly non-zero, let us write it, for notational convenience, 
in the form $(\alpha_1,\ldots,\alpha_l; \,\alpha_{k+1},\ldots,\alpha_{k+d})$.
Now suppose we wish to symmetrize $p_k(2,2,1;\,2,1)$. By independent means,
using the determinantal identities in~\cite{CFKRS2} for specific values of $k$
and polynomial interpolation, one can compute
\begin{eqnarray}
    p_k(2,2,1;\,2,1)=6 (k+2) (k^2-10) (k+1)^2\,p_k(0)\,.
\end{eqnarray}
On the other hand,
the first iteration of the  symmetrization algorithm applied to $p_k(2,2,1\,;\,2,1)$ produces
\begin{eqnarray}
    p_k(2,2,1;\,2,1)&=&\frac{1}{k-1}\,[\,-\,p_k(3,2,1;\,2)-p_k(2,3,1;\,2)-p_k(2,2,2;\,2)
     \nonumber \\ 
     &&-\sum_{r=1}^{k-3}
    p_k(2,2,1,\overbrace{0,\ldots,0}^{\textrm{$r-1$
    zeros}},1;\,2) -p_k(2,2,1;\,3)\,]\,. 
\end{eqnarray}
Therefore, by routine symmetry considerations,
\begin{eqnarray}
    p_k(2,2,1;\,2,1)&=&\frac{1}{1-k}\,[\,2\,p_k(3,2,1;\,2)+p_k(2,2,2;\,2)  \nonumber\\
    &&+(k-3)\,p_k(2,2,1,1;\,2)+p_k(2,2,1;\,3)\,]\,.
\end{eqnarray}
We verify the two sides of the above equality are equal. By independent means,
\begin{eqnarray}
    p_k(3,2,1;\,2)&=& 2 (k+2) (k+1) (k^4-58 k^2 + 417)\,p_k(0)\\
    p_k(2,2,2;\,2)&=& -72 (k+2) (k+1) (k^2-11) \,p_k(0)\\
    p_k(2,2,1;\,3)&=& 6 (k-3) (k-4) (k+4) (k+3) (k+2) (k+1)\,p_k(0) \\
    p_k(2,2,1,1;\,2)&=& -8 (k+3) (k+2) (k+1) (2k^2 -47) \,p_k(0)\,.
\end{eqnarray}
Using some algebraic manipulations, we thus obtain
\begin{eqnarray}
    &&2\,p_k(3,2,1;\,2)+p_k(2,2,2;\,2)+p_k(2,2,1;\,3) \\
    &&+(k-3)\,p_k(2,2,1,1;\,2) = -6 (k-1) (k+2) (k^2-10)
    (k+1)^2\,p_k(0)\,.\nonumber
\end{eqnarray}
Upon dividing the above by $1-k$, we arrive at $p_k(2,2,1;\,2,1)$, as claimed.

\subsection{The second step: from $p_k(\lambda;0)$ to $N_k^0(\lambda)$}

According to the lemma~\ref{lem:sym lem}, the function $p_k(\alpha)$, where
$\alpha\in \mathbb{Z}_{\geq 0}^{2k}$, can be written in terms of functions of
the form
\begin{eqnarray}
    \label{eq:p_k lambda}
    p_k(\lambda;0):=\frac{(-1)^k}{k!^2}\frac{1}{(2\pi i)^{2k}}\oint \cdots
    \oint \frac{\Delta^2(z_1,\ldots,z_{2k})\,e^{\frac{1}{2}\sum_{i=1}^k
    z_i-z_{k+i}}}{\prod_{i,j=1}^k (z_i-z_{k+j})\,\prod_{i=1}^{2k}z_i^{2k}}
    \,z_1^{\lambda_1}\ldots z_{k}^{\lambda_k}\,dz_1\ldots dz_{2k}\,,
\end{eqnarray}
where  $\lambda=(\lambda_1,\ldots,\lambda_k)\in \mathbb{Z}_{\geq 0}^k$, and
 $p_k(\lambda;0)=p_k(\lambda_1,\ldots,\lambda_k,0,\ldots,0)$.
We now show that the variables $z_{k+1},\ldots,z_{2k}$, can be completely 
 eliminated from the above expression for $p_k(\lambda;0)$. 
 That is, the integral~\eqref{eq:p_k lambda} 
 can be made to involve variables in the first half only
 (so the ``cross-terms'' are eliminated).

\begin{lem}\label{lem:red0}

Let $\lambda=(\lambda_1,\ldots,\lambda_k)\in \mathbb{Z}_{\geq 0}^k$, $k \ge 2$, and define
\begin{eqnarray}
    \label{eq:N_k^0 lambda}
    N_k^0(\lambda):=\frac{(-1)^{\binom{k}{2}}}{k!}\,\frac{1}{(2\pi
    i)^k}\oint \cdots \oint \frac{\Delta^2(z_1,\ldots,z_k)\,e^{\sum_{i=1}^k
    z_i}}{\prod_{i=1}^k z_i^{2k}} \,z_1^{\lambda_1}\ldots
    z_{k}^{\lambda_k}\,dz_1\ldots dz_k \,.
\end{eqnarray}
Then $p_k(\lambda;0) = N_k^0(\lambda)$. 
\end{lem}

\begin{proof}

    Applying lemma~\ref{lem:sym2} to (\ref{eq:p_k lambda}) with 
    $f(z_1,\ldots,z_{2k})= \exp(\frac{1}{2} \sum_1^{2k} z_i)$,
    so that $f(0,\ldots,0)=1$,
\begin{eqnarray}
    p_k(\lambda;0)=\frac{(-1)^k}{k!^2}\frac{1}{(2\pi i)^{2k}}\oint \cdots \oint
    \frac{\Delta^2(z_1,\ldots,z_{2k})\,e^{\sum_{i=1}^k z_i}}{\prod_{i,j=1}^k
    (z_i-z_{k+j})\,\prod_{i=1}^{2k}z_i^{2k}} \,z_1^{\lambda_1}\ldots
    z_{k}^{\lambda_k}\,dz_1\ldots dz_{2k}\,. 
\end{eqnarray}
Also,
\begin{eqnarray}
    \Delta^2(z_1,\ldots,z_{2k})=
    \Delta^2(z_1,\ldots,z_k)\,\Delta^2(z_{k+1},\ldots,z_{2k})\,\prod_{i,j=1}^k
    (z_i-z_{k+j})^2\,.
\end{eqnarray}
Therefore,
\begin{eqnarray}
    p_k(\lambda;0)=\frac{(-1)^k}{k!^2}\frac{1}{(2\pi i)^{2k}}\oint \cdots \oint
    \frac{\Delta^2(z_1,\ldots,z_k)\,e^{\sum_{i=1}^k z_i}}{\prod_{i=1}^k
    z_i^{2k}} \,z_1^{\lambda_1}\ldots z_{k}^{\lambda_k} \, \times  \nonumber \\
    \, \\ 
    \oint \cdots \oint
    \frac{\Delta^2(z_{k+1},\ldots,z_{2k})\,\prod_{i,j=1}^k
    (z_i-z_{k+j})}{\prod_{i=1}^k z_{k+i}^{2k}}\,dz_{k+1}\ldots dz_{2k}
    \,dz_1\ldots dz_k\,.\qquad\qquad\nonumber
\end{eqnarray}
The polynomial $\Delta^2(z_{k+1},\ldots,z_{2k})$ is homogeneous of degree $2\binom{k}{2}=k^2-k$. 
Also, the polynomial $\prod_{i,j=1}^k (z_i-z_{k+j})$ is homogeneous of degree $k^2$.
Note that the coefficient of $z_{k+1}^{k-1}\ldots z_{2k}^{k-1}$ in
     $\Delta^2(z_{k+1},\ldots,z_{2k})$ is
     $(-1)^{\binom{k}{2}}\,k!$,
and the coefficient of $z_{k+1}^k\ldots z_{2k}^k$ in 
    $\prod_{i,j=1}^k (z_i-z_{k+j})$ is $(-1)^{k^2}=(-1)^k$.
So, computing the residue at $z_{k+1}=\ldots = z_{2k}=0$ gives
\begin{eqnarray}
    \frac{(-1)^k}{(2\pi i )^k}\,\oint \cdots \oint
    \frac{\Delta^2(z_{k+1},\ldots,z_{2k})\,\prod_{i,j=1}^k
    (z_i-z_{k+j})}{\prod_{i=1}^k z_{k+i}^{2k}} \,dz_{k+1}\ldots
    dz_{2k}=(-1)^{\binom{k}{2}}\,k!\,.\nonumber \\
    \,
\end{eqnarray}
The lemma follows.
\end{proof}

\section{An algorithm to compute $N_k^0(\lambda)$}
\label{sec: Nk alg section}

Given a multivariate formal power series $Q(z_1,\ldots,z_k)$, define
\begin{eqnarray}
    [\lambda_1,\ldots,\lambda_k]_Q:=\textrm{ Coefficient of $\prod_{j=1}^k
    z_j^{2k-\lambda_j-1}$ in $Q(z_1,\ldots,z_k)$}\,. 
\end{eqnarray} 
Let 
\begin{eqnarray}
    F(z_1,\ldots,z_k):= \Delta^2(z_1,\ldots,z_k)\,e^{\sum_{i=1}^k
    z_i}\,.
\end{eqnarray}
Then,
\begin{eqnarray}
    \label{eq:Nk again}
    \frac{1}{(2\pi i)^k}\oint \cdots \oint
    \frac{F(z_1,\ldots,z_k)}{\prod_{i=1}^k
    z_i^{2k}} \,z_1^{\lambda_1}\ldots z_{k}^{\lambda_k}\,dz_1\ldots dz_k =
    [\lambda_1,\ldots,\lambda_k]_F \,. \nonumber \\
    \,
\end{eqnarray}
Also, by its definition,
\begin{eqnarray}
    \label{eq:N_k^0 b}
    p_k(\lambda_1,\ldots,\lambda_k,0,\ldots,0)=N_k^0(\lambda)
    =\frac{(-1)^{\binom{k}{2}}}{k!}\,[\lambda_1,\ldots,\lambda_k]_F\,.
\end{eqnarray}

The purpose of this section is to derive an algorithm to compute the
coefficients $[\lambda_1,\ldots,\lambda_k]_F$. As an easy by-product of the
algorithm, sharp enough upper bounds on the magnitude of these coefficients
are obtained. The algorithm comes in the form of a recursion that
  dissipates the entries of a given tuple $\lambda$, while also
  decreasing its weight.

Notice since $F$ is symmetric with respect 
to the all of the $z_j$'s, then $[\lambda_1,\ldots,\lambda_k]_F$ and $N_k^0(\lambda)$
 are symmetric with respect to all of the $\lambda_j$'s.

To help get used to the notation, note for instance, for $k\ge 2$, 
\begin{eqnarray}
    \frac{(-1)^{\binom{k}{2}}}{k!}\,
    [0,\ldots,0]_F&=&\frac{(-1)^{\binom{k}{2}}}{k!}\,\times \,\textrm{
    Coefficient of $z_1^{2k-1}\ldots z_k^{2k-1}$ in $F(z_1,\ldots,z_k)$
    }\,\nonumber\\ \, \nonumber \\ &=& 
    N_k^0(0)=p_k(0)=
    \frac{g_k}{k^2!}
    \,.
\end{eqnarray}
The last step is equation~\eqref{eq:ksb}.

We will need several lemmas, and we will make use of the function
\begin{eqnarray}
    \label{eq:G_j}
    G_j(z_1,\ldots,z_k):=\frac{F(z_1,\ldots,z_k)}{z_1-z_j}\,.
\end{eqnarray}
Notice $z_1-z_j$ divides the Vandermonde determinant in $F$, so $G_j(z_1,\ldots,z_k)$ is a polynomial.
In the lemmas to follow, we consider tuples $(\lambda_1,\ldots,\lambda_k)\in \mathbb{Z}_{\geq 0}^k$.
 Although the restriction $\lambda_j\ge 0$ is what is relevant to our problem, it is often not necessary.

\begin{lem} \label{lem:iter0}
Let $(\lambda_1,\ldots,\lambda_k)\in \mathbb{Z}_{\geq 0}^k$. Then,
\begin{eqnarray}
    \label{eq:recur F}
    [\lambda_1,\lambda_2,\ldots,\lambda_k]_F=(2k-\lambda_1)\,[\lambda_1-1,\lambda_2,\ldots,\lambda_k]_F-2\sum_{j=2}^k
    [\lambda_1,\lambda_2,\ldots,\lambda_k]_{G_j}\,.\nonumber \\
    \,
\end{eqnarray}
\end{lem}

\begin{proof}

By logarithmic differentiation, we have
\begin{eqnarray}
    \frac{\frac{\partial}{\partial z_1}
    F(z_1,\ldots,z_k)}{F(z_1,\ldots,z_k)}=1+2 \sum_{j=2}^k
    \frac{1}{z_1-z_j}\,.
\end{eqnarray}
So
\begin{eqnarray}
    \frac{\partial}{\partial z_1}
    F(z_1,\ldots,z_k)&=&F(z_1,\ldots,z_k)+2\sum_{j=2}^k
    \frac{F(z_1,\ldots,z_k)}{z_1-z_j} \nonumber \\
    &=&F(z_1,\ldots,z_k)+2\sum_{j=2}^k
    G_j(z_1,\ldots,z_k)\,. \label{eq:diff F}
\end{eqnarray}
Equating the coefficient of $\prod_{j=1}^k z_j^{2k-\lambda_j-1}$ on both sides above,
 we have
\begin{eqnarray}
    \label{eq:recur 1}
    [\lambda_1,\ldots,\lambda_k]_{\frac{\partial}{\partial z_1}
    F}=[\lambda_1,\ldots,\lambda_k]_F+2\sum_{j=2}^k
    [\lambda_1,\ldots,\lambda_k]_{G_j} \,.
\end{eqnarray}
By differentiating the power series of $F$ with respect to $z_1$, the lhs also equals
\begin{eqnarray}
    \label{eq:recur 2}
    [\lambda_1,\ldots,\lambda_k]_{\frac{\partial}{\partial z_1}
    F}=(2k-\lambda_1)\,[\lambda_1-1,\lambda_2,\ldots,\lambda_k]_F\,.
\end{eqnarray}
By substituting \eqref{eq:recur 2} into \eqref{eq:recur 1}, the lemma follows. 
\end{proof}

It is actually more convenient to rewrite the recursion \eqref{eq:recur F} in the form
\begin{eqnarray}
    \label{eq:recur F 1}
    [\lambda_1+1,\lambda_2,\ldots,\lambda_k]_F=
    (2k-\lambda_1-1)\,[\lambda_1,\lambda_2,\ldots,\lambda_k]_F-2\sum_{j=2}^k
    [\lambda_1+1,\lambda_2,\ldots,\lambda_k]_{G_j}\,.\nonumber\\
\end{eqnarray}
Also, for better readability, let us  drop entries $\lambda_j$ unaltered
from their ``original values'' in a {\em{reference tuple
$\lambda=(\lambda_1,\ldots,\lambda_k)$}}, except for the first entry
$\lambda_1$, which will always be displayed. For example, if
\mbox{$\lambda=(\lambda_1,\ldots,\lambda_k)$} is the reference tuple, then the
expressions
\begin{eqnarray}
    [\lambda_1,\lambda_j+1] \qquad \textrm{and}\qquad
    [\lambda_1+3,\lambda_k+9]\,,
\end{eqnarray}
will now stand for
\begin{eqnarray}
    [\lambda_1,\ldots,\lambda_{j-1},\lambda_j+1,\lambda_{j+1},\ldots,\lambda_k] \qquad
    \textrm{and}\qquad
    [\lambda_1+3,\lambda_2,\ldots,\lambda_{k-1},\lambda_k+9]\,,\nonumber \\
    \,
\end{eqnarray}
So now the recursion \eqref{eq:recur F 1} can be expressed more simply as
\begin{eqnarray}
    \label{eq:recur F simple}
    [\lambda_1+1]_F=(2k-\lambda_1-1)[\lambda_1]_F-2\sum_{j=2}^k
    [\lambda_1+1]_{G_j}\,.
\end{eqnarray}

\begin{lem}\label{lem:iter1}
Let $(\lambda_1,\ldots,\lambda_k)\in \mathbb{Z}_{\geq 0}^k$ be the reference tuple. Then
\begin{eqnarray}
    \label{eq:recur G}
    [\lambda_1+1]_{G_j}=[\lambda_1]_F+[\lambda_1,\lambda_j+1]_{G_j}\,.
\end{eqnarray}
In particular, for any integer $\Delta \ge -1$, and $2\le j\le k$, we have
\begin{eqnarray}
    \label{eq: recur G multi}
    [\lambda_1+1]_{G_j}=\sum_{l=0}^{\Delta} [\lambda_1-l,\lambda_j+l]_F +
    [\lambda_1-\Delta,\lambda_j+\Delta+1]_{G_j}\,.
\end{eqnarray}
\end{lem}

\begin{proof}

The relation \eqref{eq:recur G} is symmetric in the $z_j$'s, $j\ge 2$. 
So we may as well take $j=2$. Write
\begin{eqnarray}
    G_2(z_1,\ldots,z_k)&=&c_1\,
    z_1^{2k-\lambda_1-2}z_2^{2k-\lambda_2-1} z_3^{2k-\lambda_3-1} \ldots
    z_k^{2k-\lambda_k-1} \nonumber \\ \,\\ && +\, c_2\,
    z_1^{2k-\lambda_1-1}z_2^{2k-\lambda_2-2} z_3^{2k-\lambda_3-1}\ldots
    z_k^{2k-\lambda_k-1} +\cdots \,.\nonumber
\end{eqnarray}
Thus, $c_1=[\lambda_1+1]_{G_2}$, and $c_2=[\lambda_1,\lambda_2+1]_{G_2}$.
Notice
\begin{eqnarray}
    \label{eq:recur G 1}
    (z_1-z_2)\,G_2(z_1,\ldots,z_k)=
    (c_1-c_2)\,z_1^{2k-\lambda_1-1}\,z_2^{2k-\lambda_2-1}\ldots z_k^{2k-\lambda_k-1} +
    \cdots\,.
\end{eqnarray}
Since, by definition, $F(z_1,\ldots,z_k)=(z_1-z_2)\,G_2(z_1,\ldots,z_k)$,
it follows from \eqref{eq:recur G 1} that
\begin{eqnarray}
    [\lambda_1]_F=c_1-c_2=[\lambda_1+1]_{G_2}-[\lambda_1,\lambda_2+1]_{G_2}\,.
\end{eqnarray}
Equivalently,
    $[\lambda_1+1]_{G_2}=[\lambda_1]_F+[\lambda_1,\lambda_2+1]_{G_2}$.
The last part of the lemma follows by applying the recursion \eqref{eq:recur G} a total of $\Delta+1$ times.
\end{proof}

\begin{lem} \label{lem:iter2}
Let $(\lambda_1,\ldots,\lambda_k) \in \mathbb{Z}_{\geq 0}^k$ be the reference tuple. Assume
$\lambda_1\ge \lambda_j$ for $j\le k$, and define
\begin{eqnarray}
\Delta_j:=\left\lfloor \frac{\lambda_1-\lambda_j}{2}\right\rfloor\,.
\end{eqnarray}
Then,
\begin{displaymath}
    [\lambda_1-\Delta_j,\lambda_j+\Delta_j+1]_{G_j}=\left\{\begin{array}{cl}
    -\,\frac{1}{2}\,[\lambda_1-\Delta_j,\lambda_j+\Delta_j]_F & \textrm{if
    $\lambda_1-\lambda_j$ is even}\,,\\ &\\ 0 &\textrm{if
    $\lambda_1-\lambda_j$ is odd} \,.\\ \end{array}\right.
\end{displaymath}
\end{lem}

\begin{proof}

Since $F(z_1,\ldots,z_k)$ is symmetric with respect to all of the $z_j$'s,
it follows that $G_j(z_1,\ldots,z_k) = F(z_1,\ldots,z_k)/(z_1-z_j)$ is anti-symmetric with respect to $z_1$
and $z_j$; i.e.:
\begin{eqnarray}
    \label{eq:G anti sym}
    G_j(z_1,\ldots,z_j,\ldots)=-\,G_j(z_j,\ldots,z_1,\ldots)\,.
\end{eqnarray}

In particular, if we view $G_j$ as a polynomial in $z_1$ and $z_j$, and write 
\begin{eqnarray}
    G_j(z_1,\ldots,z_k)=\sum_{m,n\in \mathbb{Z}_{\geq 0}} c_{m,n} \,z_1^m z_j^n\,,
\end{eqnarray}
so the coefficients $c_{m,n}$ are now polynomials in $\{z_i\,:\, i\ne 1,j\}$, then
by the anti-symmetry of $G_j$, in \eqref{eq:G anti sym}, we have $c_{m,n}= -\,c_{n,m}$,
 and so
\begin{eqnarray}
    \label{eq:coeff relations}
    c_{m,m}=0\,,\qquad 
    c_{m+1,m}=-\,c_{m,m+1}\,.
\end{eqnarray}

Next, note
\begin{displaymath}
    (\lambda_1 -\Delta_j)-(\lambda_j+\Delta_j+1) =\left\{\begin{array}{cl} -1 &
    \textrm{if $\lambda_1-\lambda_j$ is even}\,,\\ &\\ 0 & \textrm{if
    $\lambda_1-\lambda_j$ is odd}\,. \end{array}\right.
\end{displaymath}
If $\lambda_1-\lambda_j$ is odd,  so $\lambda_1
-\Delta_j=\lambda_j+\Delta_j+1$, it follows from the first relation in \eqref{eq:coeff relations}
, with $m=2k- (\lambda_1-\Delta_j) - 1= 2k - (\lambda_j+\Delta_j+1) - 1$, that
\begin{eqnarray}
    [\lambda_1-\Delta_j,\lambda_j+\Delta_j+1]_{G_j} = 0\,.
\end{eqnarray}
On the other hand, if $\lambda_1-\lambda_j$ is even, so $\lambda_1
-\Delta_j=\lambda_j+\Delta_j$,  then  the identity
\begin{eqnarray}
    [\lambda_1-\Delta_j+1,\lambda_j+\Delta_j]_{G_j}
    =[\lambda_1-\Delta_j,\lambda_j+\Delta_j]_F+[\lambda_1-\Delta_j,\lambda_j+
    \Delta_j+1]_{G_j}\,,\nonumber \\
    \,
\end{eqnarray}
readily deducible from the recursion 
$[\lambda_1+1]_{G_j}=[\lambda_1]_F+[\lambda_1,\lambda_j+1]_{G_j}$ of
lemma~\ref{lem:iter1}, together with the second
relation in~\eqref{eq:coeff relations} applied with
$m + 1=2k-(\lambda_1-\Delta_j)-1$ and $m=2k-(\lambda_j+\Delta_j+1)-1$, imply
\begin{eqnarray}
    [\lambda_1-\Delta_j,\lambda_j+\Delta_j+1]_{G_j} = -\,\frac{1}{2}\,
    [\lambda_1-\Delta_j,\lambda_j+\Delta_j]_F\,,
\end{eqnarray}
as required.
\end{proof}

\subsection{An algorithm to compute $N_k^0(\lambda)$}
\label{sec: the alg to compute N}

We show how to compute $[\lambda_1,\ldots,\lambda_k]_F$ via a recursion. Since by relation
\eqref{eq:N_k^0 b} we have $N_k^0(\lambda) = \frac{(-1)^{\binom{k}{2}}}{k!}\,[\lambda_1,\ldots,\lambda_k]_F$,
then the said recursion can be directly used to compute $N_k^0(\lambda)$ as well. We will employ
this recursion in  \textsection{\ref{sec: consq of the algs}} to bound $N_k^0(\lambda)$. 

\begin{lem}\label{lem:theiteration}

Let $(\lambda_1+1,\lambda_2,\ldots,\lambda_k)\in \mathbb{Z}_{\geq 0}^k$. 
Assume $\lambda_1+1 \ge \lambda_j$ for $j\le k$.
Define
\begin{eqnarray}
    \Delta_j:=\left\lfloor \frac{\lambda_1-\lambda_j}{2}\right\rfloor\,,\qquad
    \delta_j:=\left\{\begin{array}{cl} -\,\frac{1}{2}\,,& \textrm{if
    $\lambda_1-\lambda_j$ is even}\\ \,\\ 0\,,& \textrm{if
    $\lambda_1-\lambda_j$ is odd.}\end{array}\right.
\end{eqnarray}
Then, with $\lambda=(\lambda_1,\ldots,\lambda_k)$ as the reference tuple, we have

\begin{eqnarray}
    \label{eq:Nk alg}
    [\lambda_1+1]_F=(2k -\lambda_1-1)\,[\lambda_1]_F-2\sum_{j=2}^k \left[
    \delta_j\,[\lambda_1-\Delta_j,\lambda_j+\Delta_j]_F +\sum_{l=0}^{\Delta_j}
    [\lambda_1-l,\lambda_j+l]_F\right] \,.\nonumber \\
    \,
\end{eqnarray}
In other words, the coefficient corresponding to the tuple
$(\lambda_1+1,\lambda_2,\ldots,\lambda_k)$, which has weight $|\lambda|+1$,
can be expressed as a linear combination involving tuples
of weight $|\lambda|$ only. 
\end{lem}
\noindent
Remark: if $\lambda_1 = \lambda_j - 1$, so $\Delta_j=-1$, then the sum over $k$ in \eqref{eq:Nk alg} vanishes, since $\delta_j=0$ in that case.
\begin{proof}

By lemma~\ref{lem:iter0},
\begin{eqnarray}
    [\lambda_1+1]_F=(2k-\lambda_1-1)[\lambda_1]_F-2\sum_{j=2}^k [\lambda_1+1]_{G_j}\,.
\end{eqnarray}
And by lemma~\ref{lem:iter1}, applied with $\Delta=\Delta_j$, we have 
\begin{eqnarray}
    [\lambda_1+1]_{G_j}=\sum_{l=0}^{\Delta_j} [\lambda_1-l,\lambda_j+l]_F +
    [\lambda_1-\Delta_j,\lambda_j+\Delta_j+1]_{G_j}\,.
\end{eqnarray}
Therefore,
\begin{eqnarray}
    [\lambda_1+1]_F=(2k -\lambda_1-1)[\lambda_1]_F-2\sum_{j=2}^k
    \left[\sum_{l=0}^{\Delta_j} [\lambda_1-l,\lambda_j+l]_F +
    [\lambda_1-\Delta_j,\lambda_j+\Delta_j+1]_{G_j}\right]\,.\nonumber \\
    \,
\end{eqnarray}
The result now follows from lemma~\ref{lem:iter2}.
\end{proof}

\subsection{Examples}
\label{sec: Nk alg examples}

Say we wish to compute $N_k^0(4,2,1,0,\ldots,0)$.  For notational convenience,
given a tuple $(\lambda_1,\ldots,\lambda_l,0,\ldots,0)\in \mathbb{Z}_{\geq
0}^k$, let us define
\begin{eqnarray}
    N_k^0(\lambda_1,\ldots,\lambda_l,0,\ldots,0)=:N_k^0(\lambda_1,\ldots,\lambda_l)\,.
\end{eqnarray}
Using this notation, the function to be computed is $N_k^0(4,2,1)$.
Lemma~\ref{lem:theiteration} and \eqref{eq:N_k^0 b} provides, on collecting terms,
\begin{eqnarray}
    N_k^0(4,2,1)&=& (2k - 4)\,N_k^0(3,2,1)- 2(k-1)\, N_k^0(3,2,1)- N_k^0(2,2,2) - \nonumber\\
    && 2\,(k-3)\, N_k^0(2,2,1,1)\nonumber \\ 
    &=& -\,2\,N_k^0(3,2,1)-N_k^0(2,2,2) - 2\,(k-3)\, N_k^0(2,2,1,1)\,.
\end{eqnarray}
Note the lhs involves a tuple of weight 7, whereas the rhs
involves tuples of weight 6 only, as should be.
By independent means, using determinantal identities in~\cite{CFKRS2} for specific values of $k$
and polynomial interpolation, we computed
\begin{eqnarray}
    N_k^0(3,2,1)&=& -3 k (k-3) (k+3) (k+2) (k+1)\,N_k^0(0)\,,\\
    N_k^0(2,2,2)&=& 24 k (k+2) (k+1)\, N_k^0(0)\,,\\
    N_k^0(2,2,1,1)&=& 12 k (k+3) (k+2) (k+1)\, N_k^0(0)\,,\\
    N_k^0(4,2,1)&=&-6 k (k+2)(k+1) (3k^2-23)\, N_k^0(0)\,.\label{eq:N 4 2 1}
\end{eqnarray}
Let us check that lemma~\ref{lem:theiteration} does in fact yield the correct $N_k^0(4,2,1)$. The rhs is 
\begin{eqnarray}
    &&[\,6 k (k-3) (k+3) (k+2) (k+1)-24 k (k+2) (k+1) \nonumber\\
    &&-24 k (k-3)(k+3) (k+2)(k+1)\,]\,N_k^0(0)\,.
\end{eqnarray}
The above can be simplified to
\begin{eqnarray}
    &&6 k (k+2)(k+1)\, [\,(k-3) (k+3) -4 - 4 (k-3)(k+3) \,] \nonumber\\
    &&=6 k (k+2)(k+1) (-3 k^2 +23)\,,
\end{eqnarray}
which agrees with~\eqref{eq:N 4 2 1}

As another example, let 
\begin{eqnarray}
    1_n:=(\overbrace{1,\ldots,1}^{\textrm{$n$ entries}},0,\ldots,0)\,.
\end{eqnarray}
Then one computes, by directly using \eqref{eq:Nk alg} and the symmetry of 
 $N_k^0(1_n)$ with respect to the $\lambda_j$'s with $j> n$,  
\begin{eqnarray}
    \label{eq:N111 1}
    N_k^0(1_n) &=& (2k - 1) N_k^0(1_{n-1}) - \sum_{j=n+1}^k N_k^0(1_{n-1}) \nonumber\\
    &=& (k + n - 1) N_k^0(1_{n-1})\,.
\end{eqnarray}
From which it follows
\begin{eqnarray}
    \label{eq:N111 2}
    N_k^0(1_n) = N_k^0(0)\prod_{j=0}^{n-1} (k+j)\,.
\end{eqnarray}
One can obtain similar simple expressions 
 for other special choices of $\lambda$.

\section{Applications of the algorithms}
\label{sec: consq of the algs}

As a consequence of the recursions in \textsection{\ref{sec: symm}} and \textsection{\ref{sec: the alg to compute N}}, 
we show that $p_k(\alpha)/p_k(0)$ grows at most polynomially in $k$, and at most exponentially in $|\alpha|$, for $|\alpha| < k/2$.
We need the following lemma.

\begin{lem}\label{lem:ni1}

Let $\lambda=(\lambda_1,\ldots,\lambda_k)\in \mathbb{Z}_{\geq 0}^k$,
such that $|\lambda|<k$. Then, 
\begin{equation}
    \label{eq:Nk simple bound to prove}
    \frac{N_k^0(\lambda)}{N_k^0(0)} \le  \frac{16^{|\lambda|}\, (\log(|\lambda|+10))^{|\lambda|} \, k^{|\lambda|}}{\lambda_1\lambda_2\ldots \lambda_{m(\lambda)}}\,.
\end{equation}
\end{lem}
\begin{proof}

Consider a tuple \mbox{$(\lambda_1+1,\lambda_2,\ldots,\lambda_k)$}, which has weight
$|\lambda|+1$. By the symmetry of $N_k^0(\lambda)$ with respect to all of the
 $\lambda_j$'s (see the remark at the beginning of \textsection{\ref{sec: the alg to compute N}}),
 we may assume $\lambda_1+1\ge \lambda_2 \ge \ldots\ge \lambda_k$.
 Without loss of generality, we may make a similar 
 assumption on the ordering of all the tuples that occur in the present proof.

Maintaining the convention whereby entries unchanged from their values in
the reference tuple $\lambda=(\lambda_1,\ldots,\lambda_k)$ are dropped, we
have by lemma~\ref{lem:theiteration}, after some simple manipulations, that
\begin{eqnarray}
    \label{eq:Nk bound 0}
    |N_k^0(\lambda_1+1)| &\le& (2k-1) \,|N_k^0(\lambda_1)|  +2\sum_{j=2}^k
    \sum_{l=0}^{\Delta_j} |N_k^0(\lambda_1-l,\lambda_j+l)|\,, 
\end{eqnarray}
where $\Delta_j = \lfloor (\lambda_1-\lambda_j)/2\rfloor$.
Note the term $\delta_j\,[\lambda_1-\Delta_j,\lambda_j+\Delta_j]_F$
that appears in the lemma is dropped because in the event
$\delta_j=-1/2$ it simply reduces the $l=\Delta_j$ term of the inner sum 
in the lemma by a factor of $1/2$, which is smaller than the stated bound.

The rhs in \eqref{eq:Nk bound 0}  involves tuples 
of weight $|\lambda|$ only, while the lhs involves a tuple of 
weight $|\lambda|+1$. This suggests inducting on $|\lambda|$.
So assume we have verified the following induction hypothesis 
for all tuples $\lambda'$ of weight $\le |\lambda|$:
\begin{eqnarray}
    \label{eq:Nk bound ind}
    \frac{|N_k^0(\lambda')|}{N_k^0(0)}
    \le \frac{16^{|\lambda'|} \, (\log(|\lambda'|+10))^{|\lambda'|}\,k^{|\lambda'|}}{\lambda_1'\lambda_2'\ldots\lambda'_{m(\lambda')}}\,.
\end{eqnarray}
We now wish to show it holds for $N_k^0(\lambda_1+1)$;
 that is, we wish to show it for tuples of weight $|\lambda|+1$.
 
By identity \eqref{eq:N111 2}, and the assumption $|\lambda|<k$,
the induction hypothesis holds for all  
$k$-tuples $\lambda'=(1,\ldots,1,0,\ldots,0)$. 
 So we may take tuples of this form
  as the base cases for the induction.
Also, notice if $\lambda_1=0$, then given our assumption 
$\lambda_1+1\ge \lambda_2\ge \ldots\ge \lambda_k$,
 the tuple $(\lambda_1+1,\lambda_2,\ldots,\lambda_k)$ must be
 of the form $(1,\ldots,1,0,\ldots,0)$, 
  and this falls within the base cases of the induction.
    Therefore, we may assume $\lambda_1>0$, so that 
 $m(\lambda_1+1,\lambda_2,\ldots,\lambda_k) = m(\lambda_1,\ldots,\lambda_k)$. 
 What we wish to show then is
\begin{eqnarray}
    \label{eq:Nk bound to prove}
    \frac{|N_k^0(\lambda_1+1)|}{N_k^0(0)} \le 
    \frac{16^{|\lambda|+1} (\log(|\lambda|+4))^{|\lambda|+1}\,k^{|\lambda|+1}}{(\lambda_1+1)\lambda_2\ldots\lambda_{m(\lambda)}}\,.
\end{eqnarray}

Consider the first term on the rhs of \eqref{eq:Nk bound 0}, 
 as well as the terms with $l=0$ in the inner sum there.
By the induction hypothesis,
\begin{eqnarray}
    \label{eq:Nk proof 1}
    \frac{|(2k-1)N_k^0(\lambda_1)|}{N_k^0(0)} + 2\sum_{j=2}^k \frac{|N_k^0(\lambda_1,\lambda_j)|}{N_k^0(0)} \le 
    4\,\frac{16^{|\lambda|}\,(\log(|\lambda|+10))^{|\lambda|+1}\,k^{|\lambda|+1}}
    {(\lambda_1+1)\lambda_2\ldots\lambda_{m(\lambda)}}\,,
\end{eqnarray}
where we used that the above sum involves $\le 4k$ tuples of weight $|\lambda|$, and
$(\lambda_1 +1)/\lambda_1 \le 2 \le \log(|\lambda|+10)$, which is valid since $\lambda_1>0$.
Also by the induction hypothesis, 
\begin{eqnarray}
    2\sum_{j=2}^{m(\lambda)} \sum_{l=1}^{\Delta_j} \frac{|N_k^0(\lambda_1- l,\lambda_j+l)|}{N_k^0(0)}
    &\le& \frac{16^{|\lambda|}\,(\log(|\lambda|+10))^{|\lambda|}\,k^{|\lambda|}}{(\lambda_1+1)\lambda_2\ldots\lambda_{m(\lambda)}}
    2\sum_{j=2}^{m(\lambda)} \sum_{l=1}^{\Delta_j} \frac{(\lambda_1+1) \lambda_j}{(\lambda_1 - l)(\lambda_j+l)}\nonumber\,.\\
\end{eqnarray}
Therefore, since $\lambda_1-l \ge (\lambda_1 + 1)/2$ for $1\le l\le \Delta_j$ and $j\le m(\lambda)$, we have
\begin{eqnarray}
    2\sum_{j=2}^{m(\lambda)} \sum_{l=1}^{\Delta_j} \frac{(\lambda_1+1)\lambda_j}{(\lambda_1 - l)(\lambda_j+l)}
    \le 4\sum_{j=2}^{m(\lambda)} \sum_{l=1}^{\Delta_j} \frac{\lambda_j}{\lambda_j+l} \le 4|\lambda| \log(|\lambda|+10)\,,
\end{eqnarray}
where we used 
$\sum_{l=1}^{\Delta_j} 1/(\lambda_j+l) \le \log(|\lambda|+10)$,
and $\sum_{j=2}^{m(\lambda)} \lambda_j \le |\lambda|$. 
Combined with $|\lambda|<k$, we obtain
\begin{eqnarray}
    \label{eq:Nk proof 2}
    2\sum_{j=2}^{m(\lambda)} \sum_{l=1}^{\Delta_j} \frac{|N_k^0(\lambda_1- l,\lambda_j+l)|}{N_k^0(0)}
    &\le& 4\,\frac{16^{|\lambda|} (\log(|\lambda|+10))^{|\lambda|+1}\,k^{|\lambda|+1}}{(\lambda_1+1)\lambda_2\ldots\lambda_{m(\lambda)}}\,.
\end{eqnarray}

Last, since by definition $\lambda_j = 0$ for $j>m(\lambda)$, and since
$N(\lambda_1-l,\lambda_j+l)$ is symmetric with respect to the $\lambda_j$'s,
we have
\begin{eqnarray}
    \label{eq:Nk proof 3}
    2\sum_{j=m(\lambda)+1}^{k} \sum_{l=1}^{\Delta_j} \frac{|N_k^0(\lambda_1- l,\lambda_j+l)|}{N_k^0(0)}
    &=& 2\,(k-m(\lambda))\,\sum_{1\le l \le \lambda_1/2} \frac{|N_k^0(\lambda_1-l,l)|}{N_k^0(0)} \nonumber \\
    &\le& 2\,\frac{16^{|\lambda|}\,(\log(|\lambda|+10))^{|\lambda|}\,k^{|\lambda|+1}}{(\lambda_1+1)\lambda_2\ldots\lambda_{m(\lambda)}}
    \sum_{1\le l\le \lambda_1/2} \frac{\lambda_1+1}{(\lambda_1 - l)l}\nonumber\\
    &\le&8\, \frac{16^{|\lambda|}\,(\log(|\lambda|+10))^{|\lambda|+1}\,k^{|\lambda|+1}}{(\lambda_1+1)\lambda_2\ldots\lambda_{m(\lambda)}}\,,
\end{eqnarray}
where we used $(\lambda_1+1)/(\lambda_1-l) \le 4$ for $l\le \lambda_1/2$,
 and $\sum_{1\le l\le \lambda_1/2} 1/l \le \log(|\lambda|+10)$.
Assembling the bounds \eqref{eq:Nk proof 1}, \eqref{eq:Nk proof 2}, and \eqref{eq:Nk proof 3},
the claim follows.
\end{proof}

\begin{thm}\label{thm:nkbd}
Let $\alpha=(\alpha_1,\ldots,\alpha_{2k})\in \mathbb{Z}_{\geq 0}^{2k}$. Then,
there exists an absolute constant $\eta$ such that as $k\to \infty$,
and uniformly in $|\alpha|< k/2$, 
\begin{equation}
    \frac{p_k(\alpha)}{p_k(0)}\ll \eta^{|\alpha|}\,(k\log (|\alpha|+10))^{|\alpha|}\,.
\end{equation}
Note, from the residue~\eqref{eq:p_k def} defining $p_k(\alpha)$,
if $\alpha_j\ge 2k$ for any $1\le j\le k$, then $p_k(\alpha)=0$.
\end{thm}

\begin{proof}
By lemma~\ref{lem:sym lem},
\begin{equation}
    |p_k(\alpha)| \le \frac{1}{\prod_{j=1}^d (k-d+j)}\,\sum_{\lambda\in
    \mathcal{S}_{\alpha}} |N_k^0(\lambda)|\,,
\end{equation}
where $d$ is the number of non-zero entries in the second half of $\alpha$
(i.e. among $\alpha_{k+1},\ldots,\alpha_k$), and $S_{\alpha}$ is a set of
tuples $\lambda \in \mathbb{Z}_{\geq 0}^k$ satisfying $|\lambda|=|\alpha|$,
 of size $|\mathcal{S}_{\alpha}|= \prod_{j=1}^d (k+d-j)$.
Since $|\lambda| = |\alpha| < k/2$, we can apply lemma~\ref{lem:ni1}
 to the $N_k^0(\lambda)$'s, which yields
\begin{eqnarray}
    |p_k(\alpha)| &\le&  \frac{|S_{\alpha}|}{\prod_{j=1}^d (k-d+j)} \,
    16^{|\alpha|} (k\log (k+10))^{|\alpha|} \,N_k^0(0) \nonumber \\
    &\ll& (48)^{|\alpha|}\,(k\log (|\alpha|+10))^{|\alpha|}\,p_k(0)\,,
\end{eqnarray}
where we used $N_k^0(0) = p_k(0)$ and the estimate
\begin{eqnarray}
    \frac{\prod_{j=1}^{d} (k+d-j)}{\prod_{j=1}^d (k-d+j)}
    = \prod_{j=0}^{d-1} \frac{1+ j/k}{1-j/k}
    \le 3^{|\alpha|}\,,
\end{eqnarray}
which holds since $d\le |\alpha| < k/2$ and so 
$(1+j/k)/(1-j/k) \le 3$ for $j<d$.
\end{proof}

Another, more precise, consequence is that $p_k(\lambda;0)/p_k(0)$ 
is a polynomial in $k$ of degree at most $|\lambda|$. This is not
 specifically used in the proof of the main theorem in this
 paper, but it is an important fact that the ideas developed 
 so far can prove fairly straightforwardly.

\begin{thm} \label{thm:poly}
Fix a positive integer $m$. Fix
$\lambda=(\lambda_1,\ldots,\lambda_m,0,\ldots,0) \in \mathbb{Z}_{\geq 0}^k$.
Then, $p_k(\lambda;0)/p_k(0)$ is a polynomial in $k$ of degree $\le |\lambda|$. 
\end{thm}

\begin{proof}
We induct on $|\lambda|$. The base case is trivial. Assume that we have verified
the theorem for all tuples of weight $\le |\lambda|$ and consider the case of $|\lambda|+1$.
By symmetry, we may assume that
\begin{eqnarray}
    \lambda_1+1\ge \lambda_2\ge \cdots\ge \lambda_m\,.
\end{eqnarray}
And by the recursion in lemma~\ref{lem:theiteration}, applied with
$(\lambda_1,\ldots,\lambda_m,0,\ldots,0)$ as the reference tuple, we have
\begin{eqnarray}
    p_k(\lambda_1+1)=(2k -\lambda_1-1)\,p_k(\lambda_1)-2\sum_{j=2}^k \left[
    \delta_j\,p_k(\lambda_1-\Delta_j,\lambda_j+\Delta_j) +\sum_{l=0}^{\Delta_j}
    p_k(\lambda_1-l,\lambda_j+l)\right] \,.\nonumber \\
    \,
\end{eqnarray}

First, observe, by the induction hypothesis, $p_k(\lambda_1)/p_k(0)$ is a polynomial in $k$ of degree
at most $|\lambda|$. Therefore, $(2k-\lambda_1-1)\,p_k(\lambda_1)/p_k(0)$ is a
polynomial in $k$ of degree at most $|\lambda|+1$.

Second, since $\lambda_{m(\alpha)+1} = \ldots =\lambda_k=0$, we can collect the terms
$j=m(\alpha)+1,\ldots,k$ together in the above sum over $j$, and using
$\Delta_{m(\alpha)+1}=\ldots=\Delta_k$, we obtain
\begin{eqnarray}
    \sum_{j=2}^k \sum_{l=0}^{\Delta_j}
    p_k(\lambda_1-l,\lambda_j+l)&=&\sum_{j=2}^{m(\alpha)} \sum_{l=0}^{\Delta_j}
    p_k(\lambda_1-l,\lambda_j+l)+\sum_{j=m(\alpha)+1}^k \sum_{l=0}^{\Delta_j}
    p_k(\lambda_1-l,\lambda_j+l) \nonumber \\ 
    &=&\sum_{j=2}^{m(\alpha)}
    \sum_{l=0}^{\Delta_j}
    p_k(\lambda_1-l,\lambda_j+l)+(k-m(\alpha))\sum_{l=0}^{\Delta_{m(\alpha)+1}}
    p_k(\lambda_1-l,l)\,.\nonumber\\
    \,
\end{eqnarray}
Again, by the induction hypothesis, $p_k(\lambda_1-l,\lambda_j+l)/p_k(0)$ is a
polynomial in $k$ of degree at most $|\lambda|$, for all $2\le j\le m(\alpha)$. Also,
$m(\alpha)$ and $\Delta_j$ are independent of $k$.  Hence, the right hand side above, divided by $p_k(0)$, is a
polynomial in $k$ of degree at most $|\lambda|+1$.

Last, since $\delta_j$ is also independent of $k$, and since 
\begin{eqnarray}
    \sum_{j=2}^k \delta_j\,p_k(\lambda_1-\Delta_j,\lambda_j+\Delta_j)&=&
    \sum_{j=2}^{m(\alpha)} \delta_j\,p_k(\lambda_1-\Delta_j,\lambda_j+\Delta_j) 
    \\ 
    && +\,(k-m(\alpha))\,\delta_{m(\alpha)+1}\,
    p_k(\lambda_1-\Delta_{m(\alpha)+1},\Delta_{m(\alpha)+1})\,,\nonumber 
\end{eqnarray}
it follows by another application of the induction hypothesis that
the rhs above is  a polynomial in $k$ of degree
at most $|\lambda|$, completing the proof.

\end{proof}

\section{The arithmetic factor}
\label{sec: arith factor}

The function $A(z_1,\ldots,z_{2k})$ is analytic and does not vanish in a
neighborhood of the origin (where it is equal to $a_k$). So, one may consider
the Taylor expansion,
\begin{eqnarray}
    \log A(z_1,\ldots,z_{2k})=: \log a_k+B_k\sum_{i=1}^k z_i-z_{k+i}+
    \sum_{\substack{\,\,\,\alpha\in \mathbb{Z}_{\geq 0}^{2k}\\|\alpha|>1}}
    a_{\alpha}\, z_1^{\alpha_1}\ldots z_{2k}^{\alpha_{2k}}\,.
\end{eqnarray}

The goal of this section is to produce upper bounds on the coefficients
$a_{\alpha}$ (in fact, we give an asymptotic when $m(\alpha)=1$).

Before doing so, let us introduce some notation. Let
$\lambda:=(\lambda_1,\ldots,\lambda_k)$ and $\rho:=(\rho_1,\ldots,\rho_k)$
denote tuples in $\mathbb{Z}_{\geq 0}^k$. Further, for primes $p$, define
\begin{eqnarray}
    \label{eq:S and A}
    S_{n,p}:=\sum_{|\lambda|=|\rho|=n} p^{\sum_{i=1}^k \rho_i z_{k+i}-\lambda_i z_i}
    ,\qquad
    A_p:=\prod_{i,j=1}^k \left(1-
    \frac{p^{z_{k+j}-z_i}}{p}\right)\,\sum_{n=0}^{\infty} \frac{S_{n,p}}{p^n}
    \,,\nonumber \\
    \,
\end{eqnarray}
where dependencies of $S_{n,p}$ and $A_p$ on $(z_1,\ldots,z_{2k})$
are suppressed to avoid notational clutter.

With the above notation, the arithmetic factor can be expressed as
\begin{eqnarray}
    A(z_1,\ldots,z_{2k}):=\prod_p A_p\,.
\end{eqnarray}
For any absolute constant $c>1$ say, one may write
\begin{eqnarray}
    \label{eq:log A}
    \log A(z_1,\ldots,z_{2k})&=&\overbrace{\sum_{p\le ck^2} \log
    A_p}^{\textrm{``Small primes''}} +\overbrace{\sum_{p> ck^2} \log
    A_p}^{\textrm{``Large primes''}}\,.
\end{eqnarray}

We will bound the contributions of ``the small primes'' and ``the large primes'' to a
coefficient $a_{\alpha}$, separately. To this end, split the ``the small primes''
sum into
\begin{eqnarray}
    \overbrace{\sum_{p\le ck^2} \sum_{i,j=1}^k \log \left(1-
    \frac{p^{z_{k+j}-z_i}}{p}\right)}^{\textrm{Convergence factor sum}}
    +\overbrace{\sum_{p\le ck^2} \log \left(1+\sum_{n=1}^{\infty}
    \frac{S_{n,p}}{p^n}\right)}^{\textrm{Combinatorial sum}}\,. \nonumber \\
    \,
\end{eqnarray}
(Here, we used the fact $S_{0,p}=1$.) Similarly, split the ``the large primes'' sum into
\begin{eqnarray}
    \overbrace{\sum_{p>ck^2} \left[\frac{S_{1,p}}{p}+\sum_{i,j=1}^k\log
    \left(1- \frac{p^{z_{k+j}-z_i}}{p}\right)\right]}^{\textrm{Convergence
    factor sum}} +\overbrace{\sum_{p>ck^2} \left[\log
    \left(1+\sum_{n=1}^{\infty} \frac{S_{n,p}}{p^n}\right)
    -\frac{S_{1,p}}{p}\right]}^{\textrm{Combintorial sum}}. 
\end{eqnarray}

So, the sum (over primes) has been separated into four pieces. In the next few
subsections, the contribution to $a_{\alpha}$  of each of piece is bounded, or,
in some cases, an asymptotic is provided. In the last subsection, the various
bounds are collected, then presented as a theorem.

Before we proceed, let us make two remarks. First, the symmetry
\begin{eqnarray}
    \log A(z_1,\ldots,z_{2k})=\log
    A(-z_{k+1},\ldots,-z_{2k},-z_1,\ldots,-z_k)\,,
\end{eqnarray}
implies
\begin{eqnarray}
    a_{(\alpha_1,\ldots,\alpha_k,\alpha_{k+1},\ldots,\alpha_{2k})}=
    (-1)^{|\alpha|}
    a_{(\alpha_{k+1},\ldots,\alpha_{2k},\alpha_1,\ldots,\alpha_k)}\,.
\end{eqnarray}

Second, the symmetry
\begin{eqnarray}
    \log A(z_1,\ldots,z_{2k})=\log A(z_{\sigma(1)},\ldots,z_{\sigma(k)},z_{k+\tau(1)},\ldots,z_{k+\tau(k)})\,,
\end{eqnarray}
where $\sigma$ and $\tau$ are any members of the permutation group of $\{1,\ldots,k\}$, implies
\begin{eqnarray}
    a_{(\alpha_1,\ldots,\alpha_{2k})}
    =a_{(\alpha_{\sigma(1)},\ldots,\alpha_{\sigma(k)},
    \alpha_{k+\tau(1)},\ldots,\alpha_{k+\tau(k)})}\,.
\end{eqnarray}

In particular, to understand the Taylor coefficients of $\log
A(z_1,\ldots,z_{2k})$, it is enough to understand $a_{\alpha}$ for tuples
$\alpha$ of the form
\begin{eqnarray}
    \alpha=(\alpha_1,\ldots,\alpha_l,0,\ldots,0,\alpha_{k+1},\ldots,\alpha_{k+d},0,\ldots,0)\,,\qquad
    0\le d\le l\le k\,,\qquad \alpha_i>0 \,.\nonumber\\
    \,
\end{eqnarray}
We will use the convention where if $d=0$, then $\alpha_{k+1}=\cdots=\alpha_{2k}=0$.

Throughout this section, it is assumed $k$ and $c$ (in~\eqref{eq:log A})
are large enough. For the sake of definiteness, let us require
\begin{eqnarray}
    k>1000\,,\qquad  \textrm{and}  \qquad 10<c<1000\,,
\end{eqnarray}
which will suffice.

\subsection{Contribution of ``the small primes'': via Cauchy's estimate}

\subsubsection{The combinatorial sum}
\label{sec: small primes comb}

We wish to estimate the Taylor coefficients (about zero) of
\begin{eqnarray}
    \sum_{p\le ck^2} \log \left(1+\sum_{n=1}^{\infty}
    \frac{S_{n,p}}{p^n}\right)=:\sum_{p\le ck^2} C_p \,. 
\end{eqnarray}

Fix a prime $p$. We consider the coefficient of $z_1^{\alpha_1}\ldots
z_{2k}^{\alpha_{2k}}$ in the Taylor expansion of a local factor $C_p$,
and denote it by $a_{\alpha,p}$. Since $p$ is fixed, we may drop
 the dependency on it in $S_{n,p}$. So, let us write
\begin{eqnarray}
    C_p=\log \left(1+\sum_{n=1}^{\infty} \frac{S_n}{p^n}\right)\,.
\end{eqnarray}

We consider two possibilities: $m(\alpha)=1$ or $m(\alpha)>1$. Let us first
handle the case $m(\alpha)>1$.

As explained earlier, it may be assumed $\alpha$ is of the form
\begin{eqnarray}
    \alpha=(\alpha_1,\ldots,\alpha_l,0,\ldots,0,\alpha_{k+1},\ldots,\alpha_{k+d},0,\ldots,0)\,,\qquad
    0\le d\le l\le k\,,\qquad \alpha_i>0 \,.\nonumber \\
    \,
\end{eqnarray}
By symmetry, it may be further assumed $\alpha_1\ge \cdots\ge \alpha_l$ and
\mbox{$\alpha_{k+1}\ge \cdots \ge \alpha_{k+d}$}.

There are two possibilities, either $\alpha_2=0$ or not. Assume $\alpha_2\ne
0$. A quick review of the argument to follow should show that the case
$\alpha_2=0$ is completely analogous (one will need to differentiate with respect to $z_{k+1}$ instead of $z_2$, noting the fact that since $m(\alpha)>1$ then if $\alpha_2=0$, then
$\alpha_{k+1}\ne 0$). Given the assumption $\alpha_2\ne 0$, define
\begin{eqnarray}
    C_p^{''}:=\left.\frac{\partial^2}{\partial z_1 \partial z_2}\,
    C_p\right|_{\substack{z_i=0\,,\,z_{k+j}=0\\ l< i\le k\,,\, d<j\le
    k}}\,.
\end{eqnarray}
Then
\begin{eqnarray}
    a_{\alpha,p}= \frac{1}{\alpha_1\alpha_2}
    \textrm{ Coefficient of
    $z_1^{\alpha_1-1} z_2^{\alpha_2-1}
    z_3^{\alpha_3}\ldots z_l^{\alpha_l} z_{k+1}^{\alpha_{k+1}}\ldots
    z_{k+d}^{\alpha_{k+d}}$ in $C_p^{''}$}\,. \label{eq:a alpha}
\end{eqnarray}

Define
\begin{eqnarray}
    Q:=\left.1+\sum_{n=1}^{\infty} \frac{S_n}{p^n}
    \right|_{\substack{z_i=0\,,\, z_{k+j}=0\\l<i\le k\,,\,d<j\le
    k}}\,,& Q_1:=\left.\sum_{n=1}^{\infty}
    \frac{1}{p^n}\,\frac{\partial}{\partial z_1}\,S_n
    \right|_{\substack{z_i=0\,,\, z_{k+j}=0\\l<i\le k\,,\,d<j\le
    k}}\,,\nonumber \\  \,\\ Q_2:=\left.\sum_{n=1}^{\infty}
    \frac{1}{p^n}\,\frac{\partial}{\partial z_2}\,S_n
    \right|_{\substack{z_i=0\,,\, z_{k+j}=0\\l<i\le k\,,\,d<j\le k}}\,,&
    Q_{12}:=\left.\sum_{n=1}^{\infty} \frac{1}{p^n}\,\frac{\partial^2}{\partial
    z_1\partial z_2}\,S_n\right|_{\substack{z_i=0\,,\, z_{k+j}=0\\l<i\le
    k\,,\,d<j\le k}}\,. \nonumber \\
    \,
\end{eqnarray}

By a straightforward calculation,
\begin{eqnarray}
    C_p^{''}= \frac{Q_{12}}{Q}-\frac{Q_1 Q_2}{Q^2}\,.
\end{eqnarray}
Letting
\begin{eqnarray}
    \Omega:=\left\{|z_1|=\frac{\delta}{10^6\,l}\,,\ldots\,,|z_l|
    =\frac{\delta}{10^6\,l}\,,|z_{k+1}|=\frac{\delta}{10^6\,l}\,,\ldots\,,|z_{k+d}|
    =\frac{\delta}{10^6\,l}\right\}\,,\nonumber \\
    \,
\end{eqnarray}
with $\delta>0$ chosen so that $Q\ne 0$ on or inside $\Omega$ (such a $\delta$
exists), it follows from~\eqref{eq:a alpha} and Cauchy's estimate that
\begin{eqnarray} \label{corr5}
    |a_{\alpha,p}| \le \left(\frac{\delta}{10^6\,l}\right)^{2-|\alpha|}
    \left[\frac{\max_{\Omega} |Q_{12}|}{\min_{\Omega} |Q|}+\frac{\max_{\Omega}
    |Q_1 |^2}{\min_{\Omega} |Q|^2}\right]\,. 
\end{eqnarray}

Now, set
\begin{eqnarray}
    \delta=\frac{1}{1000\,\log (ck^2)}\,.
\end{eqnarray}
We do not know this is a valid choice of $\delta$ a priori, but we will know
this a posteriori. 

\vspace{2mm}

\noindent
{\bf{The Denominator.}}$\,\,\,$We first estimate $\min_{\Omega} |Q|$. So, let
\begin{eqnarray}
    \mu:=(\mu_1\,,\,\ldots\,,\,\mu_l)\,,\qquad
    \tau:=(\tau_1\,,\,\ldots\,,\tau_d)\,, \qquad \mu\in \mathbb{Z}_{\geq
    0}^l\,,\,\tau \in \mathbb{Z}_{\geq 0}^d\,.
\end{eqnarray}
Then, define
\begin{eqnarray}
    Q^{(\mu,\tau)}:=\left.\frac{\partial^{|\mu|+|\tau|}\,Q}{\partial
    z_1^{\mu_1}\ldots \partial z_l^{\mu_l}\,\partial z_{k+1}^{\tau_1}\ldots
    \partial z_{k+d}^{\tau_d}} \, \right|_{\substack{z_i=0\,,\,z_{k+j}=0\\ 1\le
    i\le l\,,\,1\le j\le d}}\,. 
\end{eqnarray}
It follows
\begin{eqnarray}
    Q=Q^{(0)}+\sum_{|\mu|+|\tau|\ge 1} \frac{Q^{(\mu,\tau)}}{\mu_1!\ldots
    \mu_l!\,\tau_1!\ldots\tau_d!}\,z_1^{\mu_1}\ldots
    z_l^{\mu_l}\,z_{k+1}^{\tau_1} \ldots z_{k+d}^{\tau_d}\,,
\end{eqnarray}
where by definition,
\begin{eqnarray}
    Q^{{(0)}} = \sum_{n=1}^{\infty} \frac{1}{p^n}\,\binom{k+n-1}{n}^2\,.
\end{eqnarray}
Let
\begin{equation}
    \label{eq:mathcal D}
    \mathcal{D} :=
    \sum_{|\mu|+|\tau|\ge 1} \frac{|Q^{(\mu,\tau)}|}{\mu_1!\ldots
    \mu_l!\,\tau_1!\ldots\tau_d!}\,|z_1^{\mu_1}\ldots z_l^{\mu_l}\,
    z_{k+1}^{\tau_1}\ldots z_{k+d}^{\tau_d}|
\end{equation}
We shall show there exists an absolute constant $\eta_1\in (0,1)$ such that
\begin{eqnarray}
    \mathcal{D} \le \eta_1\, Q^{(0)}
\end{eqnarray}
for
\begin{eqnarray}
    (z_1\,,\ldots\,,z_l\,,z_{k+1}\,,\ldots\,,z_{k+d})\in \Omega\,.
\end{eqnarray}
From that it follows
\begin{eqnarray}
    \min_{\Omega} |Q|\ge (1-\eta_1)
    Q^{(0)}=(1-\eta_1)\left[1+\sum_{n=1}^{\infty}
    \frac{1}{p^n}\,\binom{k+n-1}{n}^2\right]\,,
\end{eqnarray}
because by setting all $z_j=0$ in~\eqref{eq:S and A} we have
\begin{equation}
    \sum_{|\lambda|=|\rho|=n} 1 = \binom{k+n-1}{n}^2. 
\end{equation}
The latter can be seen by arranging $k+n-1$ `dots' in a row and breaking
them into $k$ non-negative summands by selecting $k-1$ of the dots as barriers.

Now, bounding the rhs of~\eqref{eq:mathcal D} on $\Omega$ gives
\begin{eqnarray}
    \mathcal{D}\le  \sum_{\substack{h+g\ge 1\\ h\le l\,,\, g\le
    d}}\,\frac{1}{(10^6\,l)^{h+g}} \sum_{\substack{m(\mu)=h\\m(\tau)=g}}
    \frac{|Q^{(\mu,\tau)}|\,\,\delta^{|\mu|+|\tau|}}{\mu_{i_1}!\ldots
    \mu_{i_h}!\,\tau_{j_1}!\ldots\tau_{j_g}!}\,.
\end{eqnarray}
Here we have used $h \leq |\mu|$ and $g\leq |\tau|$ so that
$(10^6\,l)^{h+g} \leq (10^6\,l)^{|\mu|+|\tau|}$

Let us examine the inner sum above. For $h$ and $g$ any non-negative integers
satisfying $h+g\ge 1$, $h\le l$, $g\le d$, we have
\begin{eqnarray}
    \label{eq:Q expansion}
    \left.Q^{\,}\right|_{\substack{z_i=0\,,\,z_{k+j}=0\\ h<i\le l\,,\, g<j\le
    d}}  =1+\sum_{n=1}^{\infty} \frac{1}{p^n}\,\sum_{a=0}^n \sum_{b=0}^n
    \,\,\binom{k+n-h-a-1}{n-a} \binom{k+n-g-b-1}{n-b}\, \times \nonumber \\
    \,\\
    \sum_{\substack{\lambda=(\lambda_1,\ldots,\lambda_h)\,,\,\lambda_i\ge 0\\
    \rho=(\rho_1,\ldots,\rho_g)\,,\,\rho_i\ge 0\\|\lambda|= a\,,\,|\rho|= b}}
    p^{\rho_1 z_{k+1}+\cdots+\rho_g z_{k+g} -\lambda_1 z_1-\cdots-\lambda_h
    z_h}\,.\nonumber
\end{eqnarray}
In the above, the binomial coefficient $\binom{k+n-h-a-1}{n-a}$, for example,
represents the number of ways to write $n-a$ as the sum of $k-h$ non-negative
summands. Notice if $h=0$ then the inner-most sum vanishes unless $a=0$, and if $h=k$ then $\binom{k+n-h-a-1}{n-a}$ is 0 unless $a=n$, in which case it is 1; analogously if $g=0,k$.

So, for $\mu=(\mu_1,\ldots,\mu_h,0,\ldots,0)\in \mathbb{Z}_{\geq 0}^l$, and
$\tau=(\tau_1,\ldots,\tau_g,0,\ldots,0)\in \mathbb{Z}_{\geq 0}^d$, such that $|\mu|+|\tau| \ge 1$,
\begin{eqnarray}
    |Q^{(\mu,\tau)}| \le \sum_{n=h}^{\infty} \frac{1}{p^n}
    \sum_{a=h}^{n}\sum_{b=g}^{n} \binom{k+n-h-a-1}{n-a} \binom{k+n-g-b-1}{n-b}
    \times \qquad\qquad\qquad\qquad\qquad \nonumber \\ \, \\
    \sum_{\substack{\lambda=(\lambda_1,\ldots,\lambda_h)\,,\,\lambda_i\ge 1\\
    \rho=(\rho_1,\ldots,\rho_g)\,,\,\rho_i\ge 1\\ |\lambda|=a\,,\,|\rho|=b}}
    \,\, (\lambda_1 \log p)^{\mu_{1}}\ldots (\lambda_h \log p)^{\mu_{h}}
    (\rho_1 \log p)^{\tau_{1}}\ldots (\rho_g \log
    p)^{\tau_{g}}\,.\qquad\qquad\qquad\qquad  \nonumber
\end{eqnarray}
The sums over $a,b$ start at $h,g$ respectively because the partial derivatives
of~\eqref{eq:Q expansion} vanish if the exponent in the innermost sum
has fewer than $h$ of $z_1,\ldots,z_h$ or fewer than $g$ of $z_{k+1},\ldots,z_{k+g}$.
For the same reason, we can start the sum over $n$ at $\max(h,g)$, and choose
$h$.

Therefore, by symmetry  of $Q$ with respect to $z_1,\ldots,z_l$, and, separately, with respect to $z_{k+1},\ldots,z_{k+d}$,
\begin{eqnarray}
    \sum_{\substack{m(\mu)=h\\m(\tau)=g}}
    \frac{|Q^{(\mu,\tau)}|\,\,\delta^{|\mu|+|\tau|}}{\mu_{i_1}!\ldots
    \mu_{i_h}!\,\tau_{j_1}!\ldots\tau_{j_g}!} \le \sum_{n=h}^{\infty}
    \frac{1}{p^n}  \sum_{a=h}^{n}\sum_{b=g}^{n} \binom{k+n-h-a-1}{n-a}
    \binom{k+n-g-b-1}{n-b} \times \qquad \\ \,\nonumber\\
    \sum_{\substack{\lambda=(\lambda_1,\ldots,\lambda_h)\,,\,\lambda_i\ge 1\\
    \rho=(\rho_1,\ldots,\rho_g)\,,\,\rho_i\ge 1\\ |\lambda|=a\,,\,|\rho|=b}}
    \,\,\,\sum_{\substack{m(\mu)=h\\m(\tau)=g}} \frac{(\delta \lambda_1 \log
    p)^{\mu_{i_1}}\ldots (\delta \lambda_h \log p)^{\mu_{i_h}} (\delta \rho_1
    \log p)^{\tau_{j_1}}\ldots (\delta \rho_g \log
    p)^{\tau_{j_g}}}{\mu_{i_1}!\ldots \mu_{i_h}! \tau_{j_1}!\ldots
    \tau_{j_g}!}\,. \nonumber \label{eq:bound 1}
\end{eqnarray}
Summing over $h+g\ge 1$, $h\le l$, $g\le d$, we obtain
\begin{eqnarray} \label{corr1}
    \mathcal{D}\le \sum_{\substack{h+g\ge 1\\ h\le l\,,\, g\le d}}
    \frac{1}{(10^6\,l)^{h+g}}\,\binom{l}{h}\binom{d}{g}\,\sum_{n=h}^{\infty}
    \frac{1}{p^n} \,\times
    \qquad\qquad\qquad\qquad\qquad\qquad\qquad\qquad\qquad\qquad
     \\ \,\nonumber \\ \qquad\qquad\sum_{a=h}^{n}\sum_{b=g}^{n}
    \binom{k+n-h-a-1}{n-a} \binom{k+n-g-b-1}{n-b} \binom{a-1}{h-1}
    \binom{b-1}{g-1} \,p^{\delta (a+b)}\nonumber
    \,.
\end{eqnarray}

In the above sum, the binomial coefficients $\binom{l}{h}$ and $\binom{d}{g}$
represent the number of ways to select the $\mu_i$'s and $\tau_i$'s so that
$m(\mu)=h$ and $m(\tau)=g$. Also, the factor $p^{\delta (a+b)}$ arises from
$\exp(\log(p) (\lambda_1+\dots +\lambda_h + \rho_1+\ldots +\rho_g))$, writing
this as a product of $\exp$'s and using the Taylor series about 0 for $\exp(x)$
to produce the terms in the innermost sum of~\eqref{eq:bound 1}. There are two
special cases: When $g=0$, the quantity $\binom{b-1}{g-1}$ is defined to be
zero unless $b=0$, where it is defined to be 1, and when $g=k$, the quantity
$\binom{k+n-g-b-1}{n-b}$ is 0, unless $b=n$, in which case it is 1. Similar considerations apply to
special values of $h$.

For $n< 8k$ say, use the following estimates. First, notice that
$\binom{k+n-h-a-1}{n-a}$ is the number of ways to write $n-a$ as the sum of
exactly $k-h$ non-negative integers, and $\binom{a-1}{h-1}$ is equal to the
number of ways to write $a$ as the sum of exactly $h$ positive integers.
Therefore, $\binom{k+n-h-a-1}{n-a}\binom{a-1}{h-1}$ is at most the number of
ways to write $n$ as the sum of exactly $k$ non-negative integers, where the
first $k-h$ parts sum to $n-a$ and the last $h$ parts sum to $a$. So by summing
over $a$, we see
\begin{eqnarray}
   \sum_{a=h}^{n} \binom{k+n-h-a-1}{n-a} \binom{a-1}{h-1} \le \binom{k+n-1}{n}\,, 
\end{eqnarray}
where $\binom{k+n-1}{n}$ is the number of ways to 
write $n$ as the sum of exactly $k$ non-negative integers.
In the range $100h \le n$, we thus obtain
\begin{eqnarray}
    \sum_{a=h}^{100 h-1} \binom{k+n-h-a-1}{n-a} \binom{a-1}{h-1} \,p^{\delta a}
    \le \binom{k+n-1}{n}\,p^{100 \delta h} \,.
\end{eqnarray}
In the range $100h \le a \le n$, estimate (205) is no longer good enough for our purposes. Instead, we note
\begin{eqnarray}
    \frac{\binom{k+n-h-a-1}{n-a}}{\binom{k+n-1}{n} } 
    = \frac{\prod_{j=0}^{a-1} (n-j)\,\prod_{j=1}^h (k-j)}{\prod_{j=1}^{a+h} (k + n -j)}
    \le (1+k/n)^{-a}\,(1+n/k)^{-h}\,,
\end{eqnarray}
it follows
\begin{eqnarray}
    \label{eq:comb sum small primes est 1}
    \sum_{a=100h}^{n} \binom{k+n-h-a-1}{n-a} \binom{a-1}{h-1} \,p^{\delta a}
    \le \binom{k+n-1}{n} \sum_{a=100h}^n \frac{\binom{a-1}{h-1}\,p^{\delta
    a}}{\left(1+\frac{k}{n}\right)^a\,\left(1+\frac{n}{k}\right)^h}\,.\nonumber \\
    \,
\end{eqnarray}
Recalling $\delta=\frac{1}{1000\,\log (ck^2)}$ and $p\le ck^2$,
 we have $p^{\delta} \le 1.001$. 
Writing $a=100h+m$, one deduces 
\begin{eqnarray}
    \frac{\binom{100h + m -1}{h-1}}{\binom{100h-1}{h-1}}
    = \frac{\prod_{j=0}^{m-1} (100h+j)}{\prod_{j=0}^{m-1} (99h+j+1)} 
    \le (1+1/99)^m\,. 
\end{eqnarray}
Also, for $n< 8k$, it holds $1+k/n \ge 9/8$.
So it is seen that the sum  \eqref{eq:comb sum small primes est 1} 
is bounded by
\begin{eqnarray}
    \le 100\,\binom{k+n-1}{n}\,\frac{\binom{100h-1}{h-1} p^{100\delta
    h}}{\left(\frac{9}{8}\right)^{100h}}\le 100\,\binom{k+n-1}{n}\,,
\end{eqnarray}
where, in the last inequality, we used $\binom{100h-1}{h-1} \le (100h)^h/h! \le 300^h$, $p^{100 \delta h} \le (1.2)^h$, and $(9/8)^{100h} \ge 1000^h$. Put together, we have
\begin{eqnarray}
\sum_{n=h}^{8k-1} \frac{1}{p^n} \sum_{a=h}^{n}\sum_{b=g}^{n} \binom{k+n-h-a-1}{n-a} \binom{k+n-g-b-1}{n-b}\,\times\qquad\qquad\qquad \nonumber\\
\, \\
\binom{a-1}{h-1} \binom{b-1}{g-1} \,p^{\delta (a+b)}\le 10000\,p^{100 \delta(h+g)}\,Q^{(0)} \,.\nonumber 
\end{eqnarray}
For $n\ge 8k$, use the estimate 
\begin{eqnarray}
\sum_{a=h}^{n} \binom{k+n-h-a-1}{n-a} \binom{a-1}{h-1} \,p^{\delta a} \le \binom{k+n-1}{n} p^{\delta n} \,,  
\end{eqnarray}
which, again, is deducible via a combinatorial interpretation of the sum. This estimate yields
\begin{eqnarray}
\sum_{n=8k}^{\infty} \frac{1}{p^n} \sum_{a=h}^{n}\sum_{b=g}^{n} \binom{k+n-h-a-1}{n-a} \binom{k+n-g-b-1}{n-b}\,\times\qquad\qquad\qquad\\
\,\nonumber \\
\binom{a-1}{h-1} \binom{b-1}{g-1} \,p^{\delta (a+b)}\le \sum_{n=8k}^{\infty} \frac{p^{2 \delta n}}{p^n} \binom{k+n-1}{n}^2 \,.\nonumber 
\end{eqnarray}
Collecting the bounds so far, and using some straightforward manipulations, we have by \eqref{corr1} that $\mathcal{D}$ is bounded by
\begin{eqnarray}
    \sum_{\substack{h+g\ge 1\\ h\le l\,,\, g\le d}}
    \frac{1}{(10^6\,l)^{h+g}}\,\binom{l}{h}\binom{d}{g}\,\left[10000\,p^{\delta 100
    (h+g)}\, Q^{(0)}+ \sum_{n=8k}^{\infty} \frac{p^{2\delta n}}{p^n}\,
    \binom{k+n-1}{n}^2 \right]\le  \nonumber \\ \,\nonumber\\
    \sum_{\substack{h+g\ge 1\\ h\le l\,,\, g\le d}} \frac{10000\, p^{100 \delta
    (h+g)} l^{h+g}}{(10^6\,l)^{h+g}} \left[Q^{(0)}+ \frac{p^{16\delta
    k}}{p^{8k}}\,\binom{9k-1}{8k}^2 \,\sum_{j=0}^{\infty} \frac{p^{2\delta
    j}}{p^j}\, \left(\frac{9}{8}\right)^{2j} \right]\le \frac{Q^{(0)}}{2}\,.
    \nonumber \\ \,
\end{eqnarray}
Here we have used the assumption that $d\leq l$ in the inequality
$\binom{l}{h}\binom{d}{g}\leq l^{h+g}$. Also,
note in the last inequality we used the following observation: since $Q^{(0)}$
contains the term $\frac{1}{p^k}\,\binom{2k-1}{k}^2$, and since
\begin{eqnarray}
    \frac{\frac{1}{p^{8k}}\,\binom{9k-1}{8k}^2}{\frac{1}{p^k}\,\binom{2k-1}{k}^2}=\frac{1}{p^{7k}}\,
    \prod_{l=1}^{k-1} \left( 1+\frac{7k}{k+l}\right)^2\le \frac{8^{2k}}{2^{7k}} = \frac{1}{2^k}\,
\end{eqnarray}
(the above uses $(1+7k/(k+l))<8$ and $p\geq 2$), it follows
\begin{eqnarray}
    \frac{p^{16\delta k}}{p^{8k}}\,\binom{9k-1}{8k}^2 \le
    \left(\frac{p^{16\delta}}{2}\right)^k\,Q^{(0)}\le \frac{Q^{(0)}}{10}\,.
\end{eqnarray}
In sum, we have shown
\begin{eqnarray}\label{corr6}
\max_{\Omega} \mathcal{D}\le \frac{1}{2}\,Q^{(0)}\qquad \Rightarrow\qquad \min_{\Omega}|Q|\ge \frac{1}{2}\,Q^{(0)}\,.
\end{eqnarray}
\vspace{2mm}
\noindent
{\bf{The Numerator.}}
Having disposed of $\min_{\Omega} |Q|$, we direct our attention to $\max_{\Omega} |Q_{12}|$ and $\max_{\Omega} |Q_1|^2$. 
We deal with $\max_{\Omega} |Q_{12}|$ first. We will show there exists an absolute constant $\eta_2$ such that
\begin{eqnarray}
    \max_{\Omega} |Q_{12}|\le \eta_2 \,l^3\, \frac{(\log p)^2}{p}\,Q^{(0)}\,.
\end{eqnarray}

First, note over $\Omega$,
\begin{eqnarray}
    \frac{|Q_{12}|}{(\log p)^2}\le \sum_{n=2}^{\infty} \frac{1}{p^n}\,\sum_{a=2}^n\sum_{b=0}^n \binom{k+n-l-a-1}{n-a} \binom{k+n-d-b-1}{n-b}\,\times \nonumber \\
    \,\nonumber \\
    p^{\delta\,\frac{a+b}{l}} \sum_{\substack{\lambda=(\lambda_1,\ldots,\lambda_l)\,,\,\lambda_i\ge 0\\ \rho=(\rho_1,\ldots,\rho_d)\,,\,\rho_i\ge 0\\|\lambda|=a-2\,,\,|\rho|=b}}(\lambda_1+1)(\lambda_2+1)\,.
\end{eqnarray}
(Note the sum over $a$ starts at 2 instead of 0 because. otherwise, either the derivative with respect to $z_1$ or $z_2$ will vanish.) Therefore, since $(\lambda_1+1)(\lambda_2+1) \le a^2$,
\begin{eqnarray}\label{corr2}
    \frac{|Q_{12}|}{(\log p)^2} \le \sum_{n=2}^{\infty} \frac{1}{p^n}\,\sum_{a=2}^n\sum_{b=0}^n \binom{k+n-l-a-1}{n-a} \binom{k+n-d-b-1}{n-b}\,p^{\delta\,\frac{a+b}{l}} a^2 \sum_{\substack{\lambda=(\lambda_1,\ldots,\lambda_l)\,,\,\lambda_i\ge 0\\ \rho=(\rho_1,\ldots,\rho_d)\,,\,\rho_i\ge 0\\|\lambda|=a-2\,,\,|\rho|=b}} 1\,.
\nonumber \\
    \,
\end{eqnarray}
When $n< 8k $, it follows by considering the ranges $b < 100d$ and $100d \ge b \le n$ separately as before, while noting that $d\le l$ by hypothesis, that
\begin{eqnarray}\label{corr3}
    &&\sum_{b=0}^n  \binom{k+n-d-b-1}{n-b}\,p^{\frac{\delta b}{l}} \sum_{\substack{\rho=(\rho_1,\ldots,\rho_d)\\ \rho_i\ge 0\,,\,|\rho|=b}} 1= \sum_{b=0}^n  \binom{k+n-d-b-1}{n-b}\,\binom{d+b-1}{b}\,p^{\frac{\delta b}{l}} \le  \nonumber \\
    && \qquad \binom{k+n-1}{n}\,p^{\frac{100 \delta d}{l}}+ \binom{k+n-1}{n}\sum_{b=100 d}^n \frac{\binom{d+b-1}{b}\,p^{\frac{\delta b}{l}}}{\left(1+\frac{k}{n}\right)^b\,\left(1+\frac{n}{k}\right)^d} \le 100\,\binom{k+n-1}{n}\,.\nonumber\\
    \,
\end{eqnarray}
When $n\ge 8k$, we have
\begin{eqnarray}\label{corr4}
    \sum_{b=0}^n  \binom{k+n-d-b-1}{n-b}\,p^{\frac{\delta b}{l}} \sum_{\substack{\rho=(\rho_1,\ldots,\rho_d)\\ \rho_i\ge 0\,,\,|\rho|=b}} 1&\le& \binom{k+n-1}{n}\,p^{\frac{\delta n}{l}}\,.
\end{eqnarray}
In the above expressions, when $d=0$, the quantity $\binom{d+b-1}{b}$ is interpreted as 0 unless $b=0$. Similar care should be taken in interpreting expressions when $l$ or $d$ equals $k$. In any case, if we define
\begin{eqnarray}
    \mathcal{N}:=\sum_{n=2}^{8k} \frac{1}{p^n}\,\binom{k+n-1}{n} \sum_{a=2}^{n} \binom{k+n-l-a-1}{n-a} \,p^{\frac{\delta a}{l}}\,a^2 \sum_{\substack{\lambda=(\lambda_1,\ldots,\lambda_l)\\\lambda_i\ge 0\,,\,|\lambda|=a-2}} 1\,, 
\end{eqnarray}
then, after a little bit of work combining \eqref{corr2}, \eqref{corr3}, and \eqref{corr4}, we have generously
\begin{eqnarray}
    \frac{|Q_{12}|}{(\log p)^2}\le 100\,\mathcal{N} +100\,(8k)^2\,\frac{p^{16 \delta k }}{p^{8k}}\,\binom{9k-1}{8k}^2\le 100\,\mathcal{N}+\frac{1}{p}\,Q^{(0)}\,.
\end{eqnarray}
So, we just need to bound $\mathcal{N}$. To this end, note
\begin{eqnarray}
    \sum_{a=2}^n \binom{k+n-l-a-1}{n-a} \,p^{\frac{\delta a}{l}}\,a^2\,\sum_{\substack{\lambda=(\lambda_1,\ldots,\lambda_l)\\\lambda_i\ge 0\,,\,|\lambda|=a-2}} 1 = \sum_{a=2}^n \binom{k+n-l-a-1}{n-a} \binom{l+a-3}{a-2} \,p^{\frac{\delta a}{l}}\,a^2\,. \nonumber \\
\,
\end{eqnarray}
Define
\begin{equation}
    M:= \left\lceil \frac{c\, k}{\sqrt{p}}\right\rceil\,.
\end{equation}
Further define
\begin{eqnarray}
    \Sigma_1&:=&\sum_{n=2}^{M-1} \frac{1}{p^n}\,\binom{k+n-1}{n} \sum_{a=2}^n \binom{k+n-l-a-1}{n-a}\binom{l+a-3}{a-2}\,p^{\frac{\delta a}{l}}\,a^2 \nonumber\\
\,\\
    \Sigma_2&:=&\sum_{n=M}^{\infty} \frac{1}{p^n}\,\binom{k+n-1}{n} \sum_{a=2}^n \binom{k+n-l-a-1}{n-a}\binom{l+a-3}{a-2}\,p^{\frac{\delta a}{l}}\,a^2\,.\nonumber\\
\end{eqnarray}
In particular,
\begin{eqnarray}
    \mathcal{N}\le \Sigma_1+\Sigma_2\,.
\end{eqnarray}

We bound $\Sigma_1$. Observe that
\begin{eqnarray}
    \sum_{a=2}^{100 l-1} \binom{k+n-a-l-1}{n-a}\binom{l+a-3}{a-2} \,a^2\,p^{\frac{\delta a}{l}}
    \le \sum_{a=2}^{100 l} \binom{k+n-3}{n-2} \,a^2\,p^{\frac{\delta a}{l}}\le \binom{k+n-3}{n-2}\,(100 l)^3\,.\nonumber \\
\,
\end{eqnarray}
Also, for $n<M$, 
\begin{eqnarray}
    \frac{k}{n-2}\ge \frac{k}{M}\ge \frac{\sqrt{p}}{c}\,.
\end{eqnarray}
Therefore,
\begin{eqnarray}
    \sum_{a=100 l}^n \binom{k+n-a-l-1}{n-a}\binom{l+a-3}{a-2} \,a^2\,p^{\frac{\delta a}{l}}
    \le \binom{k+n-3}{n-2}\sum_{a=100 l}^n \frac{\binom{l+a-3}{a-2}\,a^2\,p^{\frac{\delta a}{l}}}{\left(1+\frac{k}{n-2}\right)^{a-2}}
    \le 100 \binom{k+n-3}{n-2}\,.\nonumber \\
\,
\end{eqnarray}
In summary, 
\begin{eqnarray}
    \Sigma_1\le (100 l )^3\, \sum_{n=2}^{M-1} \frac{1}{p^n}\,\binom{k+n-1}{n}^2\,\frac{ \binom{k+n-3}{n-2}}{\binom{k+n-1}{n}}
    \le \frac{(100 l)^3}{\left(1+\frac{k}{n}\right)^2}\,Q^{(0)}
    \le \frac{ \eta_3\,l^3}{p}\,Q^{(0)}\,, \nonumber \\
\end{eqnarray}
where $\eta_3$ is some absolute constant. As for $\Sigma_2$, note
\begin{eqnarray}
    \sum_{a=2}^n \binom{k+n-a-l-1}{n-a}\binom{l+a-3}{a-2}\,a^2\,p^{\frac{\delta a}{l}}
    \le \binom{k+n-3}{n-2}\,n^3\,p^{\delta n}\,.
\end{eqnarray}
Therefore,  using the change of variable $n = M + j$, we have
\begin{eqnarray}
    \Sigma_2&\le& \sum_{n=M}^{\infty} \frac{n^3\,p^{\delta n}}{p^n}\,\binom{k+n-1}{n} \binom{k+n-3}{n-2} \\
    \,\nonumber \\
    &\le& \frac{M^3\,p^{\delta M}}{p^M} \binom{k+M-1}{M}\binom{k+M-3}{M-2} 
    \sum_{j=0}^{\infty} \frac{(1+j/M)^3\,p^{\delta j}}{p^j}\,\left(1+\frac{2\,k}{M}\right)^{2j}\,.\nonumber
\end{eqnarray}
Since
\begin{eqnarray}
    \sum_{j=0}^{\infty} \frac{(1+j/M)^3\,p^{\delta j}}{p^j}\,\left(1+\frac{2\,k}{M}\right)^{2j}
    \le \sum_{j=0}^{\infty} \left(1+\frac{j}{M}\right)^3 \,\left(\frac{p^{\delta/2}}{\sqrt{p}}+\frac{2}{c}\right)^{2j}
    \le \eta_4\,,
\end{eqnarray}
where $\eta_4$ is some absolute constant, it follows
\begin{eqnarray}
    \Sigma_2\le \eta_4\,\frac{M^3\,p^{\delta M}}{p^M} \binom{k+M-1}{M}\binom{k+M-3}{M-2}\,,
\end{eqnarray}
Now, define
\begin{equation}
    M_1:= \left\lfloor \frac{5 k}{\sqrt{p}}\right\rfloor\,.
\end{equation}
Note $Q^{(0)}$ contains the term $\frac{1}{p^{M_1}}\,\binom{k+M_1-1}{M_1}^2$. Thus,
\begin{eqnarray}
    \frac{\frac{1}{p^M}\, \binom{k+M-1}{M} \binom{k+M-3}{M-2}}{Q^{(0)}} 
    \le \left(\frac{M}{k}\right)^2\,\frac{1}{p^{M-M_1}}\,\left(1+\frac{k}{M_1+1}\right)^{2(M-M_1)}\,.
\end{eqnarray}
Note,
\begin{eqnarray}
    \frac{M}{k}\le \frac{4\,c}{\sqrt{p}}\,,
\end{eqnarray}
and
\begin{eqnarray}
    \frac{1}{p^{M-M_1}}\,\left(1+\frac{k}{M_1+1}\right)^{2(M-M_1)}
    \le \left(\frac{1}{\sqrt{p}}+\frac{1}{5}\right)^{2(M-M_1)} \,.
\end{eqnarray}
Therefore, for some absolute constant $\eta_5$, we have
\begin{eqnarray}
    \Sigma_2\le \frac{\eta_5}{p}\,  M^3\,p^{\delta M} \,\left(\frac{1}{\sqrt{p}}+\frac{1}{5}\right)^{2(M-M_1)}\,Q^{(0)}\,.
\end{eqnarray}
Since $M-M_1\ge \frac{c\,k}{2\,\sqrt{p}}-1$, we have
\begin{eqnarray}
    p^{\delta M}\,\left(\frac{1}{\sqrt{p}}+\frac{1}{5}\right)^{2(M-M_1)} 
    \le 2\, e^{\frac{c\, k}{500\,\sqrt{p}}}\, (0.91)^{\frac{c\,k}{\sqrt{p}}}\le 2\,(0.92)^{\frac{c\,k}{\sqrt{p}}}\,. 
\end{eqnarray}
Hence,
\begin{eqnarray}
    M^3\,p^{\delta M} \,\left(\frac{1}{\sqrt{p}}+\frac{1}{5}\right)^{2(M-M_1)} 
    \le \left(\frac{c\,k}{\sqrt{p}}\right)^3\,(0.92)^{\frac{c\,k}{\sqrt{p}}}\le \eta_6\,, 
\end{eqnarray}
for some absolute constant $\eta_6$. So, there exists an absolute constant $\eta_7$ such that
\begin{eqnarray}
    \Sigma_2\le \frac{\eta_7}{p}\,Q^{(0)}\,.
\end{eqnarray}

Assembling previous bounds together, we thus obtain
\begin{eqnarray} \label{corr7}
    \max_{\Omega} |Q_{12}|\ll l^3\,\frac{(\log p)^2}{p}\,Q^{(0)}\,, 
\end{eqnarray}
as claimed. The case $\max_{\Omega} |Q_1|^2$ is similar. There, we obtain
\begin{eqnarray} \label{corr8}
    \max_{\Omega} |Q_1|^2\ll \frac{(\log p)^2}{p}\,\left[Q^{(0)}\right]^2\,. 
\end{eqnarray}
\vspace{2mm}
\noindent
{\bf{Summary.}}$\,\,\,$Combining \eqref{corr5}, \eqref{corr6}, \eqref{corr7}, \eqref{corr8}, and the fact $l \le m(\alpha)$, we have therefore shown the existence of an absolute constant $\eta_8$ such that
\begin{eqnarray}
    |a_{\alpha,p}|\ll (\eta_8\,m(\alpha))^{|\alpha|}\,\left(\log
    k\right)^{|\alpha|-2}\,\frac{(\log p)^2}{p}\,\qquad \textrm{ for
    $m(\alpha)>1$}\,.
\end{eqnarray}
Thus, when $m(\alpha)>1$, the contribution to $a_{\alpha}$ of the combinatorial
sum corresponding to ``the small primes'' is
\begin{eqnarray}
    \label{eq:small primes comb m alpha > 1}
    \ll (\eta_8\,m(\alpha))^{|\alpha|} \, (\log k)^{|\alpha|-2} \sum_{p\le
    ck^2} \frac{(\log p)^2}{p}\ll (\eta_8\,m(\alpha)\,\log
    k)^{|\alpha|}\,.
\end{eqnarray}

Finally, the case $m(\alpha)=1$ can be handled analogously. In that case, we
obtain for some absolute constant $\eta_9$,
\begin{eqnarray}
    |a_{\alpha,p}|\ll \left(\eta_9\, \log k\right)^{|\alpha|-2}\,\frac{(\log
    p)^2}{\sqrt{p}}\,\qquad \textrm{ for $m(\alpha)=1$}\,.
\end{eqnarray}
Thus, when $m(\alpha)=1$, the contribution to $a_{\alpha}$ of the combinatorial
sum corresponding to ``the small primes'' is
\begin{eqnarray}
    \label{eq:small primes comb m alpha = 1}
    \ll \left(\eta_9\, \log k\right)^{|\alpha|-2} \sum_{p\le ck^2} \frac{(\log
    p)^2}{\sqrt{p}}\ll \left(\eta_9\, \log k\right)^{|\alpha|-1} k\,.
\end{eqnarray}

\subsubsection{The convergence factor sum}

In this subsection, we redefine, for convenience, $C_p$ and $a_{\alpha,p}$ 
of the previous subsection.

We wish to bound the Taylor coefficients (about zero) of
\begin{eqnarray}\label{eq:thm62proof0}
    \sum_{p\le ck^2} \sum_{i,j=1}^k
    \log\left(1-\frac{p^{z_{k+j}-z_i}}{p}\right)=:\sum_{p\le ck^2}
    C_p\,,
\end{eqnarray}
where, again, we redefined $C_p$ to avoid notational clutter.
 Because only two $z_i$'s appear in each term of the inner sum on the lhs,
the Taylor coefficients $a_{\alpha,p}$ of a local factor $C_p$ are zero except for the
coefficients of monomials of the type $z_i^{u}$, with $1\le i\le 2k$ (case
$m(\alpha)=1$), or $z_i^u z_{k+j}^v$, with $1\le i,j\le k$ (case
$m(\alpha)=2$). Here $u,v\in \mathbb{Z}_{\geq 0}$. By symmetry, it is
enough to consider the monomials $z_1^u$ and $z_1^u z_{k+1}^v$.

We deal with the case $m(\alpha)=1$ first. So, let $a_{\alpha,p}$ denote the
coefficient of $z_1^u$ in $C_p$, where
\begin{eqnarray}
    \alpha=(u,0,\ldots,0)\,,\qquad u\in \mathbb{Z}_{\geq 0}\,.
\end{eqnarray}

Consider the derivative
\begin{eqnarray} \label{eq:thm62proof1}
    C_p':=\left.\frac{\partial}{\partial
    z_1}\,C_p\right|_{\substack{z_i=0\\2\le i\le 2k}}= \frac{k\,\log
    p}{p}\,\frac{p^{-z_1}}{1-\frac{p^{-z_1}}{p}}\,.\nonumber\\
\end{eqnarray}
Let $\Omega:=\{|z_1|=\delta\}$, where $\delta$ is sufficiently small
(to be specified shortly).
By Cauchy's estimate,
\begin{eqnarray}
    |a_{\alpha,p}|\le \delta^{1-u}\,\max_{\Omega} C_p'\le
    \delta^{1-u}\,\frac{k\,\log
    p}{p}\,\frac{p^{\delta}}{1-\frac{p^{\delta}}{p}}\,.
\end{eqnarray}
Choosing $\delta=1/(10\,\log ck^2)$, we obtain,
\begin{eqnarray}
    |a_{\alpha,p}|\le (50\log k)^{u-1}\, \frac{50\,k\,\log p}{p}\,.
\end{eqnarray}
This uses our assumption that $k \geq 1000$,
$10 \leq c\leq 1000$, and, here, $p \leq ck^2$, so
that, with plenty of room to spare, $10 \log(c k^2) < 50 \log(k)$,
and $p^{\delta}/(1-p^{\delta-1})<50$.

Therefore, when $m(\alpha)=1$, the contribution to $a_{\alpha}$ of the
convergence factor sum corresponding to ``the small primes'' is
\begin{eqnarray}
    \label{eq:small primes conv m alpha = 1}
    \ll (50\log k)^{|\alpha|-1}\,50\,k\,\sum_{p\le ck^2} \frac{\log p}{p}\ll
    (50\,\log k)^{|\alpha|}\,k\,.
\end{eqnarray}

The case $m(\alpha)=2$ can be handled similarly. Let $a_{\alpha,p}$ now denote
the coefficient of $z_1^u z_{k+1}^v$, where
\begin{eqnarray}
    \alpha=(u,0,\ldots,0,v,0,\ldots,0)\,,\qquad u,v\in \mathbb{Z}_{\geq
    0}\,.
\end{eqnarray}

Consider the derivative
\begin{eqnarray}
    C_p'':=\left.\frac{\partial^2}{\partial z_1\partial
    z_{k+1}}\,C_p\right|_{\substack{z_i=0\,,\,z_{k+i}=0\\2\le i\le k}}=
    \frac{(\log
    p)^2}{p}\,\frac{p^{z_{k+1}-z_1}}{1-\frac{p^{z_{k+1}-z_1}}{p}}
    \left[1+\frac{1}{p}\,\frac{p^{z_{k+1}-z_1}}{1-\frac{p^{z_{k+1}-z_1}}{p}}
    \right]\,.
\end{eqnarray}

Let $\Omega:=\{|z_1|=\delta\,,\,|z_{k+1}|=\delta\}$, with $\delta$ chosen as
before. By Cauchy's estimate,
\begin{eqnarray}
    |a_{\alpha,p}|\le \delta^{2-|\alpha|}\,\max_{\Omega} C_p''\le
    \delta^{2-|\alpha|}\,\frac{50\,(\log p)^2}{p}\le (50\,\log
    k)^{|\alpha|-2}\, \frac{50\,(\log p)^2}{p}\,.
\end{eqnarray}

Therefore, when $m(\alpha)=2$, the contribution to $a_{\alpha}$ of the
convergence factor sum corresponding to ``the small primes'' is
\begin{eqnarray}
    \label{eq:small primes conv m alpha = 2}
    \ll (50\,\log k)^{|\alpha|-2}\, 50\,\sum_{p\le ck^2}\frac{(\log p)^2}{p}\ll
    (50\,\log k)^{|\alpha|}\,.
\end{eqnarray}

\subsection{Contribution of ``the large primes'': via Taylor expansions}

\subsubsection{The combinatorial sum}

Next we bound the Taylor coefficients (about zero) of
\begin{eqnarray}
   \sum_{p>ck^2} \left[\log \left(1+\sum_{n=1}^{\infty} \frac{S_{n,p}}{p^n}\right)
   -\frac{S_{1,p}}{p}\right]=:\sum_{p>ck^2} C_p\,,
\end{eqnarray}
again redefining $C_p$.
Fix a prime $p$. Since $p$ is fixed, we may drop dependency on it in $S_{n,p}$. Applying Taylor expansions to the local factor $C_p$, we obtain
\begin{eqnarray}
    \label{eq:C_p a}
    C_p= \sum_{n=2}^{\infty}
    \frac{S_n}{p^n}+\sum_{m=2}^{\infty}\frac{(-1)^{m+1}}{m}\left(
    \sum_{n=1}^{\infty} \frac{S_{n}}{p^n}\right)^m \,,
\end{eqnarray}
again redefining $C_p$.
Next, write
\begin{equation}
    \sum_{m=2}^{\infty}\frac{(-1)^{m+1}}{m}
    \left(
        \sum_{n=1}^{\infty} \frac{S_{n}}{p^n}
    \right)^m =
    \sum_{m=2}^{\infty}\frac{(-1)^{m+1}}{m}
    \sum_{n_1,n_2,\ldots,n_m \geq 1} \frac{S_{n_1}S_{n_2}\ldots
    S_{n_m}}{p^{n_1+\ldots n_m}}\,, 
\end{equation}
sort the $n_i$'s, and count them according to their multiplicity, i.e. let
$S_{n_1}S_{n_2}\ldots S_{n_m} = S_1^{\lambda_1} S_2^{\lambda_2} \ldots
S_r^{\lambda_r}$, where each $\lambda_i\geq 0$, and $\lambda_r \geq 1$ with $r$
the largest integer amongst $n_1,\ldots,n_m$. Notice that
$\lambda_1+2\lambda_2+\cdots+r\lambda_r = n_1+\ldots +  n_m$, and that
$m= \lambda_1+\cdots+\lambda_r$. The above thus equals
\begin{equation}
    \sum_{n=2}^{\infty} \frac{1}{p^n}
    \sum_{\substack{\lambda_1+2\lambda_2+\cdots+r\lambda_r=n\\\lambda_1+\cdots+\lambda_r\geq 2\\\lambda_i\ge 0\,,\,r\ge 1}}
    \frac{(-1)^{\lambda_1+\cdots+\lambda_r+1}}{\lambda_1+\cdots+\lambda_r}\,
    \frac{(\lambda_1+\cdots+\lambda_r)!}{\lambda_1 ! \cdots \lambda_r !}\,
    S_1^{\lambda_1} S_2^{\lambda_2} \ldots S_r^{\lambda_r}\,.
\end{equation}
Next, we can absorb the first sum in~\eqref{eq:C_p a} into this by changing the
condition $\lambda_1+\cdots+\lambda_r\geq 2$ to include the case
$\lambda_1+\cdots+\lambda_r=1$. But, because $\lambda_r=1$ we then have
$\lambda_1=\ldots=\lambda_{r-1}=0$. And because
$\lambda_1+2\lambda_2+\cdots+r\lambda_r=n$, we thus have $r=n$, i.e., if we
extend the sum to include $\lambda_1+\cdots+\lambda_r=1$, it introduces
precisely the terms $\sum_{n=2}^{\infty} \frac{S_n}{p^n}$. Therefore, we have
arrived at
\begin{eqnarray}
    C_p=
    \sum_{n=2}^{\infty} \frac{1}{p^n}
    \sum_{\substack{\lambda_1+2\lambda_2+\cdots+r\lambda_r=n\\\lambda_i\ge 0\,,\,r\ge 1}}
    \frac{(-1)^{\lambda_1+\cdots+\lambda_r+1}}{\lambda_1+\cdots+\lambda_r}\,
    \frac{(\lambda_1+\cdots+\lambda_r)!}{\lambda_1 ! \cdots \lambda_r !}\,
    S_1^{\lambda_1} S_2^{\lambda_2} \ldots S_r^{\lambda_r}\,.\nonumber \\
    \,
\end{eqnarray}

We consider the coefficient of $z_1^{\alpha_1}\ldots z_{2k}^{\alpha_{2k}}$ in
the Taylor expansion of $C_p$. Let us overload notation again and
denote the said coefficient by
$a_{\alpha,p}$. As noted at the beginning of the current section, it may be
assumed $\alpha$ is of the form
\begin{eqnarray}
    \alpha=(\alpha_1,\ldots,\alpha_l,0,\ldots,0,\alpha_{k+1},\ldots,\alpha_{k+d},0,\ldots,0)\,,\qquad
    0\le d\le l\le k\,,\qquad \alpha_i >0 \,.\nonumber \\
    \,
\end{eqnarray}

In particular, as far as $a_{\alpha,p}$ is concerned, it is equivalent to
consider the series
\begin{eqnarray}
    \sum_{n=\max\{l,2\}}^{\infty} \frac{1}{p^n}
    \sum_{\substack{\lambda_1+2\lambda_2+\cdots+r\lambda_r=n\\\lambda_i\ge
    0\,,\,r\ge 1}}
    \frac{(-1)^{\lambda_1+\cdots+\lambda_r+1}}{\lambda_1+\cdots+\lambda_r}\,
    \frac{(\lambda_1+\cdots+\lambda_r)!}{\lambda_1 ! \cdots \lambda_r !}\,
    S_1^{\lambda_1} S_2^{\lambda_2} \ldots S_r^{\lambda_r}\,.\nonumber \\
    \,
\end{eqnarray}
We restrict the sum over $n$ to $\max\{l,2\}$ because, in order for a term of
the form $z_1^{\alpha_1}\ldots z_{2k}^{\alpha_{2k}}$, with $\alpha_i>0$
for all $i \leq li \leq k$, we need to have at least $l$ individual $z_i$'s,
with $i\leq k$, appearing
in $S_1^{\lambda_1} S_2^{\lambda_2} \ldots S_r^{\lambda_r}$. But each term in 
the sum $S_j$ involves at most $j$ individual $z_i$'s, hence overall
we require $\sum_j=1^r j \lambda_j = n \geq l$.

Now, define
\begin{eqnarray}
    T:=\sum_{i=1}^{2k}p^{z_i}\,.
\end{eqnarray}
It is then not too hard
to see (e.g. by considering the number of ways in which $z_1^{\alpha_1}\ldots
z_l^{\alpha_l} z_{k+1}^{\alpha_{k+1}}\ldots z_{k+d}^{\alpha_{k+d}}$ can be
formed)  that $a_{\alpha,p}$ is bounded by the coefficient of
$z_1^{\alpha_1}\ldots z_l^{\alpha_l} z_{k+1}^{\alpha_{k+1}}\ldots
z_{k+d}^{\alpha_{k+d}}$ in
\begin{eqnarray}
    \sum_{n=\max\{l,2\}}^{\infty} \frac{T^{2n}}{p^n}
    \sum_{\substack{\lambda_1+2\lambda_2+\cdots+r\lambda_r=n\\\lambda_i\ge
    0\,,\,r\ge 1}} \frac{(\lambda_1+\cdots+\lambda_r)!}{\lambda_1 ! \cdots
    \lambda_r !}\,.
\end{eqnarray}

Also,
\begin{eqnarray}
    \sum_{\substack{\lambda_1+2\lambda_2+\cdots+r\lambda_r=n\\\lambda_i\ge
    0\,,\,r\ge 1}} \frac{(\lambda_1+\cdots+\lambda_r)!}{\lambda_1 ! \cdots
    \lambda_r !} \le 2^n
    \sum_{\substack{\lambda_1+2\lambda_2+\cdots+r\lambda_r=n\\\lambda_i\ge
    0\,,\,r\ge 1}} 1 \le 2^{2n}\,.
\end{eqnarray}
For the first step above use:
\begin{eqnarray}
    \frac{(\lambda_1+\cdots+\lambda_r)!}{\lambda_1 ! \cdots \lambda_r !}
    = \binom{\lambda_r}{\lambda_r}
    \binom{\lambda_{r-1}+\lambda_{r-2}}{\lambda_{r-1}}
      \ldots
      \binom{\lambda_1+\ldots\lambda_r}{\lambda_1}
\end{eqnarray}
and bound each binomial coefficient by: $\binom{m}{j} \leq 2^m$. For the second
step, the number of terms is bounded by the number of unordered partitions of
$n$, which is easily $\leq 2^{n-1}$, since the number of ordered partitions of
$n$ equals $2^{n-1}$.

Hence, $a_{\alpha,p}$ is more simply bounded by the coefficient of
$z_1^{\alpha_1}\ldots z_l^{\alpha_l} z_{k+1}^{\alpha_{k+1}}\ldots
z_{k+d}^{\alpha_{k+d}}$ in
\begin{eqnarray}
    \sum_{n=\max\{l,2\}}^{\infty} \frac{e^{2n}}{p^n}\, T^{2n} \,.
\end{eqnarray}

Let
\begin{eqnarray}
   [z_1^{\alpha_1} \ldots z_l^{\alpha_l} z_{k+1}^{\alpha_{k+1}}\ldots
   z_{k+d}^{\alpha_{k+d}}]_n:=\textrm{ Coefficient of $z_1^{\alpha_1}\ldots
   z_l^{\alpha_l} z_{k+1}^{\alpha_{k+1}}\ldots z_{k+d}^{\alpha_{k+d}}$ in
   $T^{2n}$}.
\end{eqnarray}
Setting $z_{l+1}=\ldots = z_k=0$, and $z_{k+d+1}=\ldots = z_{2k}=0$
in $T^{2n}$ gives
\begin{eqnarray}
    \left(\sum_{i=1}^{l}p^{z_i}+\sum_{i=1}^{d}p^{z_{k+i}}+(2k-l-d)\right)^{2n}
    = \sum_{j=0}^{2n} \binom{2n}{j} (2k-l-d)^{2n-j} \left(\sum_{i=1}^{l}p^{z_i}+\sum_{i=1}^{d}p^{z_{k+i}}\right)^j. \notag\\
    \,
\end{eqnarray}
Taking the multinomial expansion of the bracketed term, and
applying the operator
\begin{eqnarray}
    \left.\frac{\partial^{\alpha_1}}{\partial z_1^{\alpha_1}}\, \ldots\,
    \frac{\partial^{\alpha_l}}{\partial z_l^{\alpha_l}}
    \,\frac{\partial^{\alpha_{k+1}}}{\partial z_{k+1}^{\alpha_{k+1}}}\,
    \ldots\,\frac{\partial_{\alpha_{k+d}}}{\partial
    z_{k+d}^{\alpha_{k+d}}}\right|_{(z_1,\ldots,z_{2k})=0}\,, 
\end{eqnarray}
to $T^{2n}$, thus gives
\begin{eqnarray}
    [z_1^{\alpha_1} \ldots z_l^{\alpha_l} z_{k+1}^{\alpha_{k+1}}\ldots
    z_{k+d}^{\alpha_{k+d}}]_n=
    \qquad\qquad\qquad\qquad\qquad\qquad\qquad\qquad \\ \, \nonumber\\
    \left( \log p\right)^{|\alpha|}
    \sum_{\substack{\lambda=(\lambda_1,\ldots,\lambda_{l+d})\\ |\lambda|\le
    2n\,,\, \lambda_i\ge 1}} \binom{2n}{|\lambda|}\, (2k-l-d)^{2n-|\lambda|}
    \,\frac{\lambda_1^{\alpha_1}\ldots
    \lambda_l^{\alpha_l}\,\lambda_{l+1}^{\alpha_{k+1}}\ldots
    \lambda_{l+d}^{\alpha_{k+d}}}{\alpha_1 ! \ldots \alpha_l !\,\alpha_{k+1}
    !\ldots \alpha_{k+d}!} \,\frac{|\lambda|!}{\lambda_1 ! \ldots \lambda_{l+d}
    !}\,.\nonumber
\end{eqnarray}
Note that $0^0$ is defined to be 1 whenever it occurs. Thus,
\begin{eqnarray}
    [z_1^{\alpha_1} \ldots z_l^{\alpha_l} z_{k+1}^{\alpha_{k+1}}\ldots
    z_{k+d}^{\alpha_{k+d}}]_n \le \qquad\qquad\qquad\qquad\qquad\qquad
     \\ \,\nonumber \\ (\log p)^{|\alpha|} \sum_{j=l+d}^{2n}
    \binom{2n}{j} (2k)^{2n-j} \left(1-\frac{l+d}{2k}\right)^{2n-j} e^j
    \sum_{\substack{\lambda=(\lambda_1,\ldots,\lambda_{l+d})\\ |\lambda|=j\,,\,
    \lambda_i\ge 1}} \frac{j!}{\lambda_1!\ldots \lambda_{l+d}!}\,.\nonumber
\end{eqnarray}
The factor $e^j$ is accounted for by $e^{\lambda_1+\ldots \lambda_{l+d}}=e^j$, and
comparing to the terms obtained by multiplying out the Taylor series for each
$e^{\lambda_j}$.

By the multinomial theorem, interpreting $(l+d)^j$ to be $(1+1+\ldots+1)^j$, we
therefore get
\begin{eqnarray}
    [z_1^{\alpha_1} \ldots z_l^{\alpha_l} z_{k+1}^{\alpha_{k+1}}\ldots
    z_{k+d}^{\alpha_{k+d}}]_n \le (\log p)^{|\alpha|} \sum_{j=l+d}^{2n}
    \binom{2n}{j} (2k)^{2n-j} \left(1-\frac{l+d}{2k}\right)^{2n-j} e^j
    (l+d)^j\,.\nonumber \\
\end{eqnarray}
From this we deduce, using $\binom{2n}{j} \leq 2^{2n}$, and relabeling the sum to
start at $j=0$, that
\begin{eqnarray}
    [z_1^{\alpha_1} \ldots z_l^{\alpha_l} z_{k+1}^{\alpha_{k+1}}\ldots
    z_{k+d}^{\alpha_{k+d}}]_n  \le
    \qquad\qquad\qquad\qquad\qquad\qquad\qquad\qquad \nonumber \\ \, \\
    (\log p)^{|\alpha|}\, 4^n\, (2k)^{2n-l-d}\, e^{l+d}\, (l+d)^{l+d}
    \sum_{j=0}^{2n-l-d} \left(1-\frac{l+d}{2k}\right)^{2n-l-d-j} \left(\frac{e
    (l+d)}{2k}\right)^j\,.\nonumber
\end{eqnarray}
Hence,
\begin{eqnarray}
    [z_1^{\alpha_1} \ldots z_l^{\alpha_l} z_{k+1}^{\alpha_{k+1}}\ldots
    z_{k+d}^{\alpha_{k+d}}]_n  \le (\log p)^{|\alpha|}\, 8^n\, (2k)^{2n-l-d}\,
    e^{l+d}\, (l+d)^{l+d} \,.
\end{eqnarray}
And so
\begin{eqnarray}
    |a_{\alpha,p}|&\le& \sum_{n=\max\{l,2\}}^{\infty} \frac{1}{p^n} \,
    [z_1^{\alpha_1} \ldots z_l^{\alpha_l} z_{k+1}^{\alpha_{k+1}}\ldots
    z_{k+d}^{\alpha_{k+d}}]_n  \\ &\le& (\log p)^{|\alpha|}\,
    e^{l+d}\, (l+d)^{l+d} \sum_{n=\max\{l,2\}}^{\infty} \frac{e^{2n}\,32^n\,
    k^{2n-l-d}}{p^n}\,. 
\end{eqnarray}
Choose $c$ in $p>ck^2$ to be $c=64 e^2$ say, then
\begin{eqnarray}
    |a_{\alpha,p}| \le e^{l+d}\, (l+d)^{l+d}\, k^{2\max\{l,2\}-l-d}\,
    \frac{(\log p)^{|\alpha|}}{p^{\max\{l,2\}}}\,.
\end{eqnarray}
Finally,
\begin{eqnarray}
    \sum_{p>ck^2} |a_{\alpha,p}|&\le& e^{l+d}\, (l+d)^{l+d} \,
    k^{2\max\{l,2\}-l-d} \sum_{p>ck^2} \frac{(\log
    p)^{|\alpha|}}{p^{\max\{l,2\}}} \\ \, \nonumber \\ &\ll& e^{l+d}\,
    (l+d)^{l+d} \, k^{2\max\{l,2\}-l-d}\,\frac{|\alpha|!}{l!} \,\frac{(\log c
    k^2)^{|\alpha|-1}}{(ck^2)^{\max\{l,2\}-1}} \\ \,\nonumber \\ &\ll&
    (32 |\alpha|)^{|\alpha|}\, (\log k)^{|\alpha|-1}\, k^{2-l-d}\,. 
\end{eqnarray} 
In summary, the contribution to $a_{\alpha}$ of the combinatorial sum
corresponding to the ``the large primes'' is
\begin{eqnarray}
    \label{eq:large primes comb m alpha => 1}
    \ll (32 |\alpha|)^{|\alpha|}\, (\log k)^{|\alpha|-1}\,
    k^{2-m(\alpha)}\,.
\end{eqnarray}

\subsubsection{The convergence factor sum}

We wish to bound the Taylor coefficients (about zero) of
\begin{eqnarray}
    \sum_{p> ck^2} \left[
    \frac{S_{1,p}}{p}+\sum_{i,j=1}^k
    \log\left(1-\frac{p^{z_{k+j}-z_i}}{p}\right)\right]
    =:\sum_{p> ck^2} C_p\,.
\end{eqnarray}
Expand $\log(1-w) = - \sum_1^\infty w^m/m$, $w=p^{z_{k+j}-z_i-1}$
and cancel the $S_{1,p}/p$ term with the $m=1$ term to get
\begin{eqnarray}
    C_p=-\sum_{m=2}^{\infty} \frac{1}{m}\, \sum_{i,j=1}^k
    \frac{p^{m(z_{k+j}-z_i)}}{p^m}\,.
\end{eqnarray}

The Taylor coefficients of a local factor $C_p$ are zero except for the
coefficients of monomials of the type $z_i^{u}$, with $1\le i\le 2k$ (case
$m(\alpha)=1$), or $z_i^u z_{k+j}^v$, with $1\le i,j\le k$ (case
$m(\alpha)=2$). Here $u,v\in \mathbb{Z}_{\geq 0}$. So, by symmetry, it is
enough to consider the monomials $z_1^u$ and $z_1^u z_{k+1}^v$.

We deal with the case $m(\alpha)=1$ first. So, let $a_{\alpha,p}$ denote the
coefficient of $z_1^u$ in $C_p$, where
\begin{eqnarray}
    \alpha=(u,0,\ldots,0)\,,\qquad u\in \mathbb{Z}_{\geq 0}\,.
\end{eqnarray}
Then,
\begin{eqnarray}
    |a_{\alpha,p}|\le k\,(\log p)^u \, \sum_{m=2}^{\infty} \frac{m^u
    }{u!\,p^m}\le \frac{10\,k\,(\log p)^u}{p^2}\,. 
\end{eqnarray}
Therefore, when $m(\alpha)=1$, the contribution to $a_{\alpha}$ of the
convergence factor sum corresponding to the ``the large primes'' is
\begin{eqnarray}
    \label{eq:large primes conv m alpha = 1}
    \ll k\,\sum_{p> ck^2} \frac{(\log p)^{|\alpha|}}{p^2}\ll
    \frac{|\alpha|!\,(4\log k)^{|\alpha|-1}}{k}\,.
\end{eqnarray}
The latter inequality follows by comparing the sum to
$\int_{ck^2}^\infty \log(t)^{|\alpha|-1}/t^2 dt$ (with one less power in the exponent
to account for the density of primes), integrating by parts $|\alpha|$ times,
and using the assumption that $10 \leq c\leq 1000 \leq k$:
\begin{equation}
    \int_{ck^2}^\infty \log(t)^{|\alpha|-1}/t^2 dt =
    (|\alpha|-1)! \sum_{j=0}^{|\alpha|-1} \frac{(\log{ck^2})^j}{j!ck^2}
    \ll |\alpha|! \frac{(4\log k)^{|\alpha|-1}}{k^2} \,. 
\end{equation}

On the other hand, when $m(\alpha)=2$, the contribution to $a_{\alpha}$ is
\begin{eqnarray}
    \label{eq:large primes conv m alpha = 2}
    \ll \sum_{p> ck^2} \frac{(\log p)^{|\alpha|}}{p^2}\ll
    \frac{|\alpha|!\,(4\log k)^{|\alpha|-1}}{k^2}\,.
\end{eqnarray}

\subsection{Bounding the coefficients of the arithmetic factor}

We are now ready to state the main theorem of this section.

\begin{thm}\label{arithmfac}
The coefficients $a_{\alpha}$ in the Taylor expansion
\begin{eqnarray}
    \log A(z_1,\ldots,z_{2k})=:\log a_k+B_k\sum_{i=1}^k
    z_i-z_{k+i}+\sum_{|\alpha|>1} a_{\alpha}\, z_1^{\alpha_1}\ldots
    z_{2k}^{\alpha_{2k}}
\end{eqnarray}
satisfy
\begin{eqnarray}
    a_{\alpha}\ll \left\{\begin{array}{ll} \lambda_2^{|\alpha|}\,\left(\log
    k\right)^{|\alpha|} k+ \lambda_2^{|\alpha|}\, |\alpha|!\, \left(\log
    k\right)^{|\alpha|-1} k\,,  & \textrm{if  $m(\alpha)=1$}\\ \,&\nonumber\\
    \lambda_2^{|\alpha|}\, m(\alpha)^{|\alpha|} \left(\log
    k\right)^{|\alpha|} + \lambda_2^{|\alpha|}\, |\alpha|!\,\left(\log
    k\right)^{|\alpha|-1} k^{2-m(\alpha)}\,, & \textrm{if $m(\alpha)>1$}\\
    \end{array}\right. \\
    \label{eq:a alpha final bound}
\end{eqnarray}
as $k\to \infty$, and uniformly in $\alpha$, where $\lambda_2$ is some absolute constant. 
More simply, but slightly less precisely, 
\begin{eqnarray}
    a_{\alpha}\ll  \lambda_2^{|\alpha|}\,(\log k)^{|\alpha|}\,\left[
    m(\alpha)^{|\alpha|} \,k^{2-\min\{m(\alpha),2\}} +  |\alpha|!\,
    k^{2-m(\alpha)}\right]
\end{eqnarray}
as $k\to \infty$. Asymptotic constants are absolute.
\end{thm}

\begin{proof}
The terms $\lambda_2^{|\alpha|}\,\left(\log k\right)^{|\alpha|} k$ and
$\lambda_2^{|\alpha|}\, m(\alpha)^{|\alpha|} \left(\log k\right)^{|\alpha|}$
in \eqref{eq:a alpha final bound} come from the small primes, and arise by combining the 
contributions to $a_{\alpha}$ of:

\begin{itemize}
\item The combinatorial sum for the small primes when $m(\alpha) = 1$, \eqref{eq:small primes comb m alpha = 1}:
\begin{eqnarray} \label{eq:summarith1}
    \ll \eta_9^{|\alpha|}(\log k)^{|\alpha|-1} k\,. 
\end{eqnarray}
\item The combinatorial sum for the small primes when $m(\alpha)>1$, \eqref{eq:small primes comb m alpha > 1}: 
\begin{eqnarray} \label{eq:summarith2}
    \ll \eta_8^{|\alpha|} m(\alpha)^{|\alpha|} (\log k)^{|\alpha|}\,.
\end{eqnarray}
\item The convergence factor sum for the small primes when $m(\alpha) = 1$, \eqref{eq:small primes conv m alpha = 1}:
\begin{eqnarray} \label{eq:summarith3}
    \ll 50^{|\alpha|} (\log k)^{|\alpha|} k\,.
\end{eqnarray}
\item The convergence factor sum for the small primes when $m(\alpha) = 2$, \eqref{eq:small primes conv m alpha = 2}:
\begin{eqnarray} \label{eq:summarith4}
    \ll 50^{|\alpha|} (\log k)^{|\alpha|}\,.
\end{eqnarray}
\end{itemize}
While the terms $\lambda_2^{|\alpha|}\, |\alpha|!\, \left(\log k\right)^{|\alpha|-1} k$ and
$\lambda_2^{|\alpha|}\, |\alpha|!\,\left(\log k\right)^{|\alpha|-1} k^{2-m(\alpha)}$ 
in \eqref{eq:a alpha final bound} come from the
large primes, and arise by combining the contributions to $a_{\alpha}$ of:

\begin{itemize}
\item The combinatorial sum for the large primes when $m(\alpha) \ge 1$, \eqref{eq:large primes comb m alpha => 1}:
\begin{eqnarray}
    \label{eq:arith summary 0}
    \ll 32^{|\alpha|} |\alpha|^{|\alpha|} (\log k)^{|\alpha|-1} k^{2-m(\alpha)}\,.
\end{eqnarray}
\item The convergence factor sum for the large primes when $m(\alpha) = 1$, \eqref{eq:large primes conv m alpha = 1}:
\begin{eqnarray}
    \label{eq:arith summary 1}
    \ll 4^{|\alpha|}(|\alpha|!) (\log k)^{|\alpha|-1} /k\,.
\end{eqnarray}
\item The convergence factor sum for the large primes when $m(\alpha) = 2$, \eqref{eq:large primes conv m alpha = 2}:
\begin{eqnarray}
    \label{eq:arith summary 2}
    \ll 4^{|\alpha|}(|\alpha|!) (\log k)^{|\alpha|-1} /k^2\,.
\end{eqnarray}
\end{itemize}
The $\lambda_2^{|\alpha|} |\alpha|!$ in the statement
of the theorem accounts for both the $4^{|\alpha|}|\alpha|!$
in~\eqref{eq:large primes conv m alpha = 1} and~\eqref{eq:large primes conv m alpha = 2}, and, on using
Stirling's asymptotic, for the $(32 |\alpha|)^{|\alpha|}$
in~\eqref{eq:large primes comb m alpha => 1}.

\end{proof}

Remark: A review of the previous argument shows the 
statement of the theorem can be made more precise in the 
 case $m(\alpha)=1$: 

\begin{thm}
\label{thm:lincoeff}
For $\alpha$ satisfying $m(\alpha)=1$, define

\begin{displaymath}
    sgn(\alpha):=\left\{\begin{array}{cl}
    (-1)^{|\alpha|+1}\,,&\textrm{if $\alpha_{k+1}=\cdots=\alpha_{2k}= 0$}\,,\\
    \,\\
    -1\,,&\textrm{if $\alpha_1=\cdots=\alpha_k=0$}\,.
    \end{array}\right.
\end{displaymath}
Then, with $|\alpha|$ fixed, and as $k\to \infty$, the coefficients $a_{\alpha}$ satisfy
\begin{eqnarray}\label{eq:case m=1 asympt}
    a_{\alpha}&=& sgn(\alpha) \,\frac{k}{|\alpha|!}\,\left(\sum_{p\le k^2}
    (\log p)^{|\alpha|} \sum_{n=1}^{\infty} \frac{n^{|\alpha|-1}}{p^n}
    \right)\,   \left[1+O\left(\frac{1}{\log k}\right)\right] \nonumber \\
    \,\nonumber \\ 
    &=& sgn(\alpha) \,\frac{k}{|\alpha|!}\,\left(\sum_{p\le k^2}
    \frac{(\log p)^{|\alpha|}}{p} \right)\,   \left[1+O\left(\frac{1}{\log
    k}\right)\right]\,.
\end{eqnarray}
Asymptotic constants depend only on $|\alpha|$. In particular,
\begin{equation}
    B_k=a_{(1,0,\ldots,0)} \sim 2k \log k. 
\end{equation}
\end{thm}

\begin{proof}

Our plan is to show that, asymptotically as $k\to \infty$ and for
$|\alpha|$ fixed, the dominant contribution to the $a_{\alpha}$ when $m(\alpha)=1$
comes from the convergence factor sum corresponding to the small primes.
Notice this asymptotic is not uniform in $\alpha$, so it is not of immediate utility
in the proof of the main theorem, but it is included here because it might be
of independent interest.

To this end, by the symmetry of $A(z_1,\ldots,z_{2k})$ in the first half of the
variables $z_1,\ldots,z_k$ and, separately, in the second half $z_{k+1},\ldots,
z_{2k}$, we may assume $\alpha_1 \ge \alpha_2 \ge \cdots \ge \alpha_k$ and
$\alpha_{k+1} \ge \alpha_{k+2} \ge \cdots \ge \alpha_{2k}$. Thus, since
$m(\alpha) = 1$, then all the $\alpha_j$'s are zero except $\alpha_1$ or
$\alpha_{k+1}$, but not both.

Consider the case $\alpha_1 \ne 0$ first. Then
$\alpha = (|\alpha|,0,\ldots,0)$, and $a_{\alpha}$
is the coefficient of $z_1^{|\alpha|}$ in $A(z_1,\ldots,z_{2k})$.
By \eqref{eq:thm62proof0} and \eqref{eq:thm62proof1}, the contribution
of the convergence factor sum corresponding to the small primes to this
coefficient is

\begin{eqnarray} \label{eq:smallprimecontm1}
    \frac{k}{|\alpha|}\, \sum_{p\le c k^2} \frac{\log p}{p} \, \times
    \textrm{ Coefficient of $z_1^{|\alpha|-1}$ in }
    \frac{p^{-z_1}}{1-\frac{p^{-z_1}}{p}}\,.
\end{eqnarray}
where $10 < c < 1000$. Expanding, we obtain
\begin{eqnarray}
    \frac{p^{-z_1}}{1-\frac{p^{-z_1}}{p}} &=& \sum_{n=1}^{\infty} \frac{p^{-nz_1}}{p^{n-1}} 
    = \sum_{n=1}^{\infty} \frac{1}{p^{n-1}} \sum_{r=0}^{\infty} \frac{(-1)^r}{r!} n^r (\log p)^r z_1^r 
\end{eqnarray}
Singling out the case $r=|\alpha|-1$ above, we have
\begin{eqnarray} \label{eq:case m=1 asympt 1}
   \eqref{eq:smallprimecontm1} &=&  (-1)^{|\alpha| - 1} \frac{k}{|\alpha|!} 
   \sum_{p\le ck^2} (\log p)^{|\alpha|} \sum_{n=1}^{\infty} \frac{n^{|\alpha|-1}}{p^n} \nonumber \\
    &=& sgn(\alpha) \frac{k}{|\alpha|!} \sum_{p\le ck^2} \frac{(\log p)^{|\alpha|}}{p} \left[1 + O(1/\log k) \right]\,,
\end{eqnarray}
 where we used $(-1)^{|\alpha|-1} = (-1)^{|\alpha|+1}= sgn(\alpha)$, 
 $\sum_{p\le ck^2} (\log p)^{|\alpha|}/p \gg \log k$, and (hence) 
\begin{eqnarray}
    \sum_{p\le ck^2} (\log p)^{|\alpha|} \sum_{n=1}^{\infty} \frac{n^{|\alpha|-1}}{p^n}
    &=& \sum_{p\le ck^2} \frac{(\log p)^{|\alpha|}}{p} + O(1) \nonumber\\
    &=& \sum_{p\le ck^2} \frac{(\log p)^{|\alpha|}}{p} \left[1+O(1/\log k)\right]\,. 
\end{eqnarray}
Also, since $c$ is fixed, we may replace the range of summation $p \le ck^2$
 in \eqref{eq:case m=1 asympt 1} by $p\le k^2$ without affecting the 
 asymptotic.

The remaining contributions to $a_{\alpha}$ (which, recall, is the coefficient 
of $z_1^{|\alpha|}$) come from the combinatorial sum for the small primes, 
the combinatorial sum for the large primes, and the convergence factor sum 
for the large primes. But these contributions, which are bounded by 
\eqref{eq:summarith1}, \eqref{eq:arith summary 0}, and
  \eqref{eq:arith summary 1}, respectively, are asymptotically smaller
  than \eqref{eq:case m=1 asympt 1}, as $k\to \infty$ and for $|\alpha|$ fixed, 
  by at least a factor of $1/\log k$. 
 Put together, this yields the asymptotic \eqref{eq:case m=1 asympt} 
  in the case $\alpha_1\ne 0$.
  
Last, the analysis in the case $\alpha_{k+1}\ne 0$ is completely similar except the coefficient of 
 $z_1^{|\alpha|-1}$ in $p^{-z_1}/(1-p^{-z_1}/p)$ in \eqref{eq:smallprimecontm1} 
 is replaced by the coefficient of $z_{k+1}^{|\alpha|-1}$ in 
 $-p^{z_{k+1}}/(1-p^{z_{k+1}}/p)$, thereby changing $sgn(\alpha)$ to -1.
\end{proof}

\section{The product of zetas}
\label{sec: zeta product}

Finally, we bound the Taylor coefficients $b_{\alpha}$ of
\begin{eqnarray} \label{corr9}
    \log\left(\prod_{i,j=1}^k (z_i-z_{k+j})\zeta(1+z_i-z_{k+j})\right)=:\gamma
    k\sum_{i=1}^k z_i-z_{k+i}+\sum_{|\alpha|>1} b_{\alpha}\,
    z_1^{\alpha_1}\ldots z_{2k}^{\alpha_{2k}}\,.\nonumber\\
    \,
\end{eqnarray}
The Taylor coefficients are zero except for those of monomials of the type
$z_i^{u}$, with $1\le i\le 2k$ (case $m(\alpha)=1$), or $z_i^u z_{k+j}^v$, with
$1\le i,j\le k$ (case $m(\alpha)=2$). Here $u,v\in \mathbb{Z}_{\geq 0}$. By
symmetry, it is enough to consider the monomials $z_1^u$ and $z_1^u z_{k+1}^v$.

We deal with the case $m(\alpha)=1$ first. So, let $\alpha$ be of the form
\begin{eqnarray}
    \alpha=(u,0,\ldots,0)\,,\qquad u\in \mathbb{Z}_{\geq 0}\,.
\end{eqnarray}
Setting $z_2=\cdots=z_{2k}=0$, the lhs of \eqref{corr9} becomes
\begin{eqnarray}
    k\,\log\left[z_1\,\zeta(1+z_1)\right]=\gamma k \,z_1+\sum_{u=2}^{\infty}
    b_{(u,0,\ldots,0)}\,z_1^u\,.
\end{eqnarray}
Now, by the well-known Taylor expansion, we have
\begin{eqnarray}
    z\,\zeta(1+z)=1+\sum_{n=0}^{\infty}
    \frac{(-1)^n}{n!}\,\gamma_n\,z^{n+1}\,, 
\end{eqnarray}
where the $\gamma_n$'s are the generalized Euler constants
satisfying, $\gamma_0=\gamma=.577\ldots$, and, see Theorem 2 of \cite{B},
\begin{equation}
    \qquad |\gamma_n|\le
    4\,\frac{(n-1)!}{\pi^n}\, \quad n \geq 1 \,.
\end{equation}
Consider the derivative
\begin{eqnarray}
    \frac{d}{dz} \log \left[z \,\zeta(1+z )\right]=\frac{\sum_{n=0}^{\infty}
    \frac{(-1)^n\,(n+1)}{n!}\,\gamma_n\,z^{n}}{1+\sum_{n=0}^{\infty}
    \frac{(-1)^n}{n!}\,\gamma_n\,z^{n+1}}\,.
\end{eqnarray}
Note in particular, for $|z|<1/10$, we have
\begin{eqnarray}
    \left|\frac{d}{dz} \log \left[z \,\zeta(1+z
    )\right]\right|=\frac{8\,\sum_{n=0}^{\infty}
    \frac{1}{(10\pi)^n}}{1-\frac{4}{10}\,\sum_{n=0}^{\infty}\frac{1}{(10\pi)^n}}\le
    100\,.
\end{eqnarray}
So, by Cauchy's estimate, the coefficients $d_n$ in the expansion
\begin{eqnarray}
    \log \left[z\,\zeta(1+z)\right] =:\sum_{m=1}^{\infty} d_n\,z^n\,,
\end{eqnarray}
satisfy
\begin{eqnarray}
    |d_n|\le 100\,(10)^n\,.
\end{eqnarray}
From which it follows
\begin{eqnarray}
    |b_{\alpha}|\ll k\,(10)^{|\alpha|}\,,\qquad \textrm{when $m(\alpha)=1$}\,.
\end{eqnarray}
Analogous reasoning yields
\begin{eqnarray}
    |b_{\alpha}|\ll (100)^{|\alpha|}\,,\qquad \textrm{when $m(\alpha)=2$}\,.
\end{eqnarray}
Put together, we have
\begin{lem}\label{lem:prodz}
The coefficients $b_{\alpha}$ in the expansion
\begin{eqnarray}
    \log\left(\prod_{i,j=1}^k (z_i-z_{k+j})\zeta(1+z_i-z_{k+j})\right)=:\gamma
    k\sum_{i=1}^k z_i-z_{k+i}+\sum_{|\alpha|>1} b_{\alpha}\,
    z_1^{\alpha_1}\ldots z_{2k}^{\alpha_{2k}}\,,\nonumber\\
    \,
\end{eqnarray}
are zero when $m(\alpha)>2$, otherwise, as $k\to \infty$, and uniformly in $\alpha$, they satisfy
\begin{eqnarray}
    b_{\alpha}\ll \lambda_3^{|\alpha|}\, k^{2-m(\alpha)}\,,
\end{eqnarray}
where $\lambda_3$ is some absolute constant. Asymptotic constants are absolute.
\end{lem}

\printindex

\end{document}